\documentclass[11pt,leqno]{article}
\usepackage{amssymb}
\usepackage{mathrsfs}
\usepackage{cleveref}
\usepackage{amsfonts}
\usepackage{graphicx}
\usepackage{mathptmx}
\usepackage{latexsym,amsmath,amssymb,amsfonts,amsthm}

\allowdisplaybreaks
\usepackage{epsfig}
\usepackage{subfigure}
\usepackage{tikz}
\newcommand{\bcen}{\begin{center}}
\newcommand{\ecen}{\end{center}}

\begin{document}
\setcounter{page}{1}
\title{Determining anisotropic real-analytic metric from boundary electromagnetic information}
\author{Genqian Liu}

\date{}
\protect
\maketitle ~~~\\[-15mm]

\begin{center}
{\footnotesize   School of Mathematics and Statistics, Beijing Institute of Technology, Beijing 100081, China\\
 Emails:  liugqz@bit.edu.cn \\
 }
\end{center}


\begin{abstract}
\ For a compact, connected, oriented Riemannian $3$-manifold $(M, g)$ with smooth boundary $\partial M$, we explicitly give a local representation and a full symbol expression for the electromagnetic Dirichlet-to-Neumann map by factorizing Maxwell's equations and using an isometric transform. We prove that one can reconstruct a compact, connected, real-analytic Riemannian $3$-manifold $M$ with boundary from the set of tangential electric fields and tangential magnetic fields, given on a non-empty open subset $\Gamma$ of the boundary, of all electric and magnetic fields with tangential electric data supported in $\Gamma$. We note that for this result we need no assumption on the topology of the manifold other than compactness and connectedness, nor do we need a priori knowledge of all of $\partial M$. In addition, as a by-product of the explicit symbol expression of $\Lambda_{g,\Gamma}$, we show that for a given smooth Riemannian metric $g$, the electromagnetic Dirichlet-to-Neumann map $\Lambda_{g,\Gamma}$ uniquely determines all order tangential and normal derivatives of electromagnetic parameters $\mu$ and $\sigma$ on $\Gamma$. Therefore, $\mu$ and $\sigma$ are completely determined in $M$ by $\Lambda_{g,\Gamma}$ if these two parameter functions and  metric $g$ are all real analytic in $M$ up to $\Gamma$.
 \end{abstract}

\markright{\sl\hfill Genqian Liu  \hfill}

\section{Introduction}
\renewcommand{\thesection}{\arabic{section}}
\renewcommand{\theequation}{\thesection.\arabic{equation}}
\setcounter{equation}{0} \label{intro}

\vskip 0.45 true cm

 Let $(M, g)$ be a compact, connected, oriented Riemannian manifold with $C^1$-smooth boundary
$\partial M$, and assume that $\mbox{dim}\, M=3$.
 We consider the inverse problem of recovering metric of the medium $(M, g)$ by probing with time-harmonic electromagnetic fields from the measurements made on an open subset $\Gamma\subset \partial M$. The fields in $(M, g)$ are described by $E$ and $H$ (electric and magnetic fields), and the behavior of the fields is governed by the Maxwell equations in $M$,
 \begin{eqnarray} \label{19.7.6-1} \left\{ \!\!\begin{array}{ll} \mbox{curl}\,  E = i\omega \mu H,\\
 \mbox{curl}\, H =-i\omega(\epsilon +i\gamma/\omega) E=-i\omega \sigma E.\end{array}\right.\end{eqnarray}
 Here the constant $\omega>0$ is a fixed frequency; $\epsilon$,  $\gamma$ and $\mu$  are electric permittivity, conductivity and magnetic permeability, respectively. We assume that $\epsilon$, $\gamma$ and $\mu$ are all nonnegative functions in $M$, and $\epsilon$, $\mu$ are strictly positive on $\bar M$.
 The boundary condition is $\nu \times E\big|_{\partial M} =f$, the
tangential component of the electric field at the boundary, where $\nu$ is the unit outer normal to $\partial M$, and $A\times B$ denotes the vector product of the vectors $A$ and $B$ in the tangential space of Riemannian manifold $(M, g)$ (see section 2). The tangential
component of the magnetic field is then $\nu \times H$.
A number $\omega$ is said to be a resonant frequency if there exists a non-trivial $(E,H)$ such that  \begin{eqnarray}\label{2020.1.31-1}  \left\{\! \!\begin{array} {ll}
  \mbox{curl}\, E= i \omega \mu H \quad &\mbox{in}\,\; M,\\
 \mbox{curl}\, H=-i \omega \sigma E \quad &\mbox{in}\,\; M,\\
     \nu \times E = 0  \quad &\mbox{on}\,\; \partial M.\end{array}\right. \end{eqnarray}
It is well-known that there is a discrete set $\mathfrak{N}$ of resonant frequencies with accumulation point at infinity such that if $\omega$ is outside this set, then for any $f\in TH^{\frac{1}{2}}(\partial M)$ (the space of tangential vector fields on $\partial M$ with components in $H^{\frac{1}{2}}(\partial M)$), the system (\ref{19.7.6-1}) has a unique solution $(E,H)\in (L^2(M))^3 \times (L^2(M))^3$ (see section 7 of \cite{KeSU} or p.$\,$166 of \cite{Isak}) satisfying $\nu\times E\big|_{\partial M}= f\;$.
 We assume more regularity on $\Gamma$ and $g$. In fact we
assume $\Gamma$ is a real analytic piece of boundary, and that the metric tensor $g$ of
$M$ is real analytic up to $\Gamma$. Let $f \in C(\partial M)\cap TH^{\frac{1}{2}}(\partial M)$ with $\mbox{supp}\; f \subset \Gamma$. For $\omega\notin \mathfrak{N}$, let $(E,H)$ be the solution of (\ref{19.7.6-1}) with $\nu\times E=f$ on $\partial M$. The part boundary map, which we call the
 electromagnetic Dirichlet-to-Neumann, is
 \begin{eqnarray}\label{2020-3.28-1} \Lambda_{g,\Gamma}:  f  \mapsto \nu \times H\big|_{\Gamma}.\end{eqnarray}

Similar to one of conjectures for the conductivity equation (see \cite{LU}), an open problem for the electromagnetic field is: {\it Can one recover the real-analytic metric $g$ of $M$ uniquely up to isometry from the knowledge of the electromagnetic Dirichlet-to-Neumann map $\Lambda_{g,\Gamma}$?} Or roughly speaking, can one obtain metric information in a real-analytic, inaccessible medium $M$ when the data consist of electromagnetic field measurements on the part surface of the medium?
This open problem also arises naturally, e.g. in medical imaging applications, non-destructive material evaluation and geoelectrics (see \cite{Chew}, p$\,$166 of \cite{Isak}, or \cite{GKLU}); just as wrote by A. Kirsch (see, p.$\,$174 of \cite{Kir}) ``A fundamental question for every inverse problem is the question of {\it uniqueness}: is the information --- at least in principle --- sufficient to determine the unknown quantity?''

 The root of this problem is in Electrical Impedance Tomography (EIT). The
question in EIT is whether one can determine the (anisotropic) electrical
conductivity of a medium $\Omega$ in Euclidean space by making voltage and current
measurements at the boundary of the medium. Calder\'{o}n proposed this
problem \cite{Cald} motivated by geophysical prospection.
  The Calder\'{o}n problem has been studied extensively in past decades, and the problem of uniqueness was solved affirmatively in higher dimensional case ($\mbox{dim}\,M \ge 3$) for the smooth isotropic conductivity of a body by J. Sylvester and
G. Uhlmann \cite{SU1} and for Lipschitz conductivities on a bounded Lipschitz domain by Caro and Rogers in \cite{CR}, and in two dimensional case for $C^2$-conductivities by A. Nachman \cite{Nac1} and for Lipschitz conductivities by Brown and Uhlmann in \cite{BU}. We also refer the reader to \cite{Ale, AP, Bro, Br, BT, Cha, FKSU, GLU, KT, NSU, Nov, PPU, SU, SU2}.

Let us point out that for many inverse problems determining isotropic quantities, a key method is the construction of complex geometrical optics (CGO) solutions with a large parameter which was introduced by Sylvester and Uhlmann in \cite{SU}. Unfortunately, such a very useful method is not valid for any anisotropic inverse problem because one cannot find CGO solutions of the form $e^{\varsigma \cdot x} (1+w_{\varsigma})$ with $\varsigma\in {\mathbb{C}}^n$ and $\varsigma\cdot \varsigma=0$ for such an anisotropic problem (here the corresponding top-order differential equations have variables coefficients) unless some additional conditions are added. In a celebrated paper \cite{LU}, by considering the symbol of the Dirichlet-to-Neumann map, Lee and Uhlmann solved a series conjectures and answered the case for the anisotropic conductivity equation $\sum_{j,k=1}^n \frac{\partial }{\partial x_j}\big(\gamma^{jk} \frac{\partial u}{\partial x_k}\big)=0$. The same problem was further discussed by Lassas, Taylor and Uhlmann in \cite{LTU}, and the corresponding conjecture for $n\ge 3$ is completely solved without any topology assumption on manifold other than connectedness. The reason why the Riemannian manifolds appear naturally in the study of Electrical Impedance Tomography (EIT) (see [LeU]) is
that in dimension $n \ge  3$, the EIT problem is equivalent to the problem of
determining a Riemannian metric $g$ from Dirichlet-to-Neumann map  with
\begin{eqnarray*}  g_{ij} = ( \mbox{det}\,  \gamma^{kl})^{1/(n-2)} (\gamma^{ij} )^{-1}.\end{eqnarray*}
In addition, in order to reveal the true behavior for the Dirichlet-to-Neumann map on the boundary  (even for a bounded domain of the Euclidean space ${\mathbb{R}}^n$), the boundary normal coordinates have to be introduced which leads to Riemannian metric, symbol calculus and differential geometric techniques.

 Thus, it is our aim to discuss the more general case of recovering anisotropic metric of the electromagnetic medium by making part boundary measurements, and the equation modeling the problem is the full system of Maxwell's equations. As pointed out by V. Isakov in p.$\,$166 of \cite{Isak} ``$\cdot\cdot\cdot\cdot$ which is not a simple task $\cdot\cdot\cdot\cdot$. Maxwell's system is not elliptic, so there are additional difficulties in this case.''
 For the isotropic inverse problem, the first global uniqueness result for smooth electromagnetic parameters on a smooth domain is \cite{OPS}; another proof was given in \cite{OS}. More recently, in \cite{CZ} the case of continuously differentiable isotropic electromagnetic parameters on a domain with $C^1$ boundary was examined. The global uniqueness for isotropic Lipschitz electromagnetic parameters in a bounded Lipschitz domain was obtained in a constructive manner in \cite{Pich} by Pichler. Partial boundary data problems were studied in \cite{COR} using the reflection argument introduced in
\cite{Isako}, and in \cite{COST}, extending the ideas from \cite{ABU, KeSjU} to Maxwell's equations. The inverse parameter problem on a manifold
has been studied in \cite{OPS2}, and \cite{KeSU} considered the problem in a non-Euclidean (admissible manifold) setting.

In this paper, we give an affirmative answer to the above open problem by symbol calculus and a method of Green's function (the later was introduced by
 Lassas, Taylor and Uhlmann in \cite{LTU}). Our main results are the following:

\vskip 0.2 true cm

\noindent{\bf Theorem 1.1.} \ {\it Let $(M,g)$ be a compact, connected, oriented, smooth Riemannian $3$-manifold with $C^1$-smooth boundary $\partial M$, and let $\Gamma\subset \partial M$ be a piece of smooth surface on $\partial M$. Suppose electromagnetic parameters $\mu$, $\sigma$ and metric $g$ are smooth in $M$ up to $\Gamma$, and $\mu$} and $\,\mbox{Re}(\sigma)$ are strictly positive on $\bar M$. {\it Then the electromagnetic Dirichlet-to-Neumann map $\Lambda_{g,\Gamma}$ has a precise local expression (see Proposition 2.6). Furthermore, the full symbol matrix $\psi(x',\xi')$ of  $\Lambda_{g,\Gamma}$ \begin{eqnarray*}\psi(x',\xi')\sim \begin{bmatrix} \psi^{11}(x', \xi')  & \psi^{12} (x', \xi') \\
 \psi^{21}(x', \xi')  & \psi^{22} (x', \xi')\end{bmatrix} \end{eqnarray*}
 can explicitly be given in local boundary normal coordinates (see, the proof of Proposition 3.1), and
\begin{eqnarray*} \Lambda_{g,\Gamma} f(x')=\frac{1}{(2\pi)^2} \int_{{\mathbb{ R}}^2} e^{ix'\cdot \xi'} \psi(x',\xi') \hat{f}(\xi')\,d\xi' \quad \; \mbox{for any}\;\, f\in C(\partial M)\cap C_0^\infty(\Gamma),\end{eqnarray*}
where $\hat{f}(\xi')$ is the Fourier transform of $f$.}

\vskip 0.25 true cm

If $(M,g)$ and $(\tilde{M},\tilde{g})$ are both Riemannian manifolds, a smooth map $\varrho: M\to \tilde{M}$ is called a (Riemannian) isometry if it is a diffeomorphism that satisfies $\varrho^* \tilde{g} =g$.
If there exists a Riemannian isometry between $(M,g)$ and
$\tilde{M},\tilde{g})$, we say that they are isometric as Riemannian manifolds.

\vskip 0.25 true cm

\noindent{\bf Theorem 1.2.} \ {\it  Let $(M_1,g_1)$ and $(M_2,g_2)$ be compact, connected, oriented, real-analytic Riemannian $3$-manifolds with $C^1$-smooth boundaries $\partial M_1$ and $\partial M_2$. Assume that $\partial M_1$ and $\partial M_2$ contain a non-empty open set
$\Gamma_1 = \Gamma_2 = \Gamma$, on which each boundary is real analytic, with the metric tensors analytic up to $\Gamma_j$.
Furthermore, assume that the electromagnetic parameters $\mu$ and $\sigma$ are real analytic in $M$ up to $\Gamma$, and that $\mu$ and}  $\,\mbox{Re}(\sigma)$ {\it are strictly positive on $\bar M$. Finally, assume that the electromagnetic Dirichlet-to-Neumann maps $\Lambda_{g_1,\Gamma}$  and $\Lambda_{g_2,\Gamma}$ coincide. Then there exists a real-analytic diffeomorphism $\varrho: M_1\cup \Gamma\to  M_2\cup \Gamma$ with $\varrho\big|_{\Gamma}=\mbox{identity}$, such that $g_1=\varrho^*{g_2}$.}

\vskip 0.20 true cm

Strictly speaking in the statement above we mean by the set $\Gamma$ the sets
$\Gamma_1 \subset \partial M_1$ and $\Gamma_2\subset \partial M_2$,  which are identified by a diffeomorphism.

\vskip 0.2 true cm

When $M$ is a bounded domain in ${\mathbb{R}}^3$ and $g$ is the Euclidean metric, we immediately have the following corollary:
\vskip 0.2 true cm

 \noindent{\bf Corollary 1.3.} \ {\it Let $M \subset {\mathbb{R}}^3$ be a bounded domain
with $C^1$-smooth boundary, and let $\Gamma$ be a piece of real-analytic surface on $\partial M$. Suppose $\tilde{g}$ is a real-analytic metric in $M$ up to $\Gamma$ such that $\Lambda_{\tilde{g},\Gamma} = \Lambda_{g,\Gamma}$, where $g$ is the Euclidean metric. Assume that the electromagnetic parameters $\mu$ and $\sigma$ are real analytic in $M$ up to $\Gamma$, and that $\mu$ and}  $\,\mbox{Re}(\sigma)$ {\it are strictly positive on $\bar M$. Then there exists a real-analytic diffeomorphism $\varrho: M\cup \Gamma\to  \tilde{M}\cup \Gamma$ with $\varrho\big|_{\Gamma}=\mbox{identity}$, such that $g=\varrho^*{\tilde{g}}$.}

\vskip 0.38 true cm

In addition, as a by-product of our new method of full symbol we can show the following:

\vskip 0.2 true cm

\noindent{\bf Theorem 1.4.} \ {\it Let $(M, g)$ be a compact, connected, oriented, $C^\infty$-smooth Riemannian $3$-manifold with $C^1$-smooth boundary $\partial M$. Assume that $\,\Gamma$ is a piece of $C^\infty$-smooth surface on $\partial M$. Assume that electromagnetic parameters $\mu$ and $\sigma$ are $C^\infty$-smooth in $M$ up to $\Gamma$, and that $\mu$ and}  $\,\mbox{Re}(\sigma)$ {\it are strictly positive on $\bar M$. Then the electromagnetic Dirichlet-to-Neumann map $\Lambda_{g,\Gamma}$ uniquely determines $\frac{\partial^{|K|} \mu}{\partial x^K}$ and $\frac{\partial^{|K|} \sigma}{\partial x^K}$ on $\Gamma$ for all multi-indices $K=(k_1, k_2, k_3)$ with $|K|=k_1+k_2+k_3\ge 0$}.

\vskip 0.28 true cm

Theorem 1.4 generalizes  Joshi-McDowall's result (see \cite{JoMcD}) for electromagnetic field (in which they considered the case of a  bounded domain with $C^\infty$-smooth boundary in Euclidean space ${\mathbb{R}}^n$, and $\Gamma=\partial M$). This immediately leads to

\vskip 0.28 true cm

 \noindent{\bf Theorem 1.5.} \ {\it Let $(M,g)$ be a compact, connected, oriented, real analytic Riemannian $3$-manifold with $C^1$-smooth boundary. Let $\Gamma\subset \partial M$ be a piece of real-analytic surface on $\partial M$. Assume that the electromagnetic parameters $\mu$, $\sigma$ and the given metric $g$, are real-analytic in $M$ up to $\Gamma$, and that $\mu$ and}  $\,\mbox{Re}(\sigma)$ {\it are strictly positive on $\bar M$.
Assume also that $\omega>0$ is such that the system (\ref{19.7.6-1}) has a unique solution in $M$ with the prescribed (compact support) tangential component of the electric field on $\partial \Omega$, so that $\Lambda_{g,\Gamma}$ is well defined. Then the electromagnetic Dirichlet-to-Neumann map
$\Lambda_{g,\Gamma}$ uniquely determines the function $\mu$ and $\sigma$ in $M\cup \Gamma$.}

\vskip 0.28 true cm

Let us remark that Theorem 1.5 is not a trivial conclusion. To the best of my knowledge, this is the first novel result on determining the anisotropic parameters for the electromagnetic part boundary measurements.

\vskip 0.28 true cm

The main ideas of this paper are as follows. From the Maxwell equations we obtain a second-order partial differential equations $$\mbox{curl}\;\mbox{curl}\; E -  \big(\mbox{grad}\, (\log \mu)\big) \times \mbox{curl}\; E-\omega^2 \mu \sigma\, E=0\quad \mbox{in}\;\, M,$$ which can be rewritten as \begin{eqnarray}\label{19.9.14-1}\left\{\big(\frac{\partial^2}{\partial x_3^2} I_3\big)+ B \big(\frac{\partial }{\partial x_3}  I_3\big)+C\right\}E=0\end{eqnarray} in the boundary normal coordinates. Let us point out that it is a valuable work to obtain explicit expression for $B$ and $C$ because some new methods and calculations are needed. The second step is to use the factorization for the left-hand side of (\ref{19.9.14-1}) so that we obtain a pseudodifferential operator $\Phi$ of order one. This method  for the Laplace equation $\Delta_g u=0$ is well known (see, for example, \cite{Tre}, vol.$\,$1, pages 159-161, or \cite{Esk}) and have been well developed by Lee and Uhlmann in \cite{LU}. Since the third component of $\nu\times E$ and $\nu\times H$ are both vanish in local boundary normal coordinates, we introduce a linear isometric operator
$\wp: (C^\infty (M)))^3 \to (C^\infty (M))^2$ defined by  \begin{eqnarray*} \wp\left(\begin{bmatrix} a_1\\
 a_2\\ 0\end{bmatrix}\right)  \mapsto    \begin{bmatrix} a_1\\
 a_2\end{bmatrix} \quad \,\mbox{for any} \;\, (a_1, a_2, 0)\in (C^\infty (M))^3,\end{eqnarray*}
then we get an equivalent local representation formula for $\Lambda_{g, \Gamma}$ on $\Gamma \subset \partial M$:
 \begin{eqnarray} \label{20.4.5-0} \;\;\, \;\;\;  \; \wp(\nu \times H)=\left[\begin{matrix} \Lambda^{11} &\Lambda^{12}\\
 \Lambda^{21}& \Lambda^{22}\end{matrix} \right]\wp(\nu\times E)= \frac{1}{i\omega \mu\sqrt{|g|}}  \left[ \begin{matrix} L^{11}  &  L^{12}\\
   L^{21}  & L^{22}   \end{matrix} \right] \begin{bmatrix} -g_{12} & - g_{22}\\
 g_{11} & g_{12}  \end{bmatrix} \wp(\nu\times E),\end{eqnarray}
where $L^{jk}$ are given by (\ref{19.7.21-7}) (see section 2). This formula transforms a three-dimensional problem into a two-dimensional problem.
In order to get (\ref{20.4.5-0}) we also use the equation $\mbox{div}\, (\sigma E)=0$ in $M$, which is derived from Maxwell's equations (\ref{19.7.6-1}) (In fact, the Maxwell equations (\ref{19.7.6-1}) is equivalent to the following system of equations for electric field $E$ (see section 4):
\begin{eqnarray*} \left\{\!\!\begin{array}{ll} \mbox{curl}\; \mbox{curl}\, E
-  \big(\mbox{grad}\, (\log \mu)\big) \times \mbox{curl}\; E -\omega^2 \mu \sigma E=0 \;\, &\mbox{in}\;\,M,\\
\mbox{div}\, (\sigma E)=0 \,\; &\mbox{in}\;\, M.\,)\end{array}\right.\end{eqnarray*}
The final step is to determine $g_{\alpha\beta}$ and all their normal derivatives on $\Gamma$ by the full symbol of the  electromagnetic Dirichlet-to-Neumann  map $\Lambda_{g,\Gamma}$. This step plays a key role in this paper because some subtle new estimates are employed. More precisely, by discussing the principal symbol of pseudodifferential operator $\Lambda^{11}$ we can obtain $g_{\alpha\beta}$ on $\Gamma$ (Of course, by $\Lambda^{22}$ we can also obtain the same conclusion); and by discussing the symbol with homogeneous of degree $0$ in $\xi'$ for the operator $L^{11}+L^{22}$, we can get the first-order normal derivatives of $g^{\alpha\beta}$ on $\Gamma$; furthermore, from the symbol with homogeneous of degree $-1$ in $\xi'$ for operator $L^{11}+L^{22}$, we get the second-order normal derivatives of $g^{\alpha\beta}$ on $\Gamma$ (Let us remark that this argument cannot be omitted); generally, by using the induction we can obtain $(m+2)$-order normal derivatives on $\Gamma$ from the symbol with homogeneous of degree $-m-1$ in $\xi'$ for the operator $L^{11}+L^{22}$ for any $m\ge 0$. Combining the expansion of Taylor's series of $g_{jk}$ at every $x_0\in \Gamma$  and the analyticity of the Riemannian metric in $M$ up to $\Gamma$, we can prove Theorem 1.2 via a method of Green's function.
 The effective method of Green's function for the Laplace equation $\Delta_g u=0$ on a Riemannian manifold $M$ was introduced by Lassas, Taylor and Uhlmann in \cite{LTU}. Note that our method here is slightly different from that of \cite{LTU}. We attach a collar domain across $\Gamma$ for $M$ instead of a ``half-ball'' domain of some point $x_0\in \Gamma$ for $M$, so that we can analytically extend the Riemannian manifolds $M_j$ to the larger Riemaiannian manifolds ${\tilde{M}}_j$. The advantage of this new technique is that it can immediately ensure the identity of metric $g$ on $\Gamma$ when $M_1$ and $M_2$ are shown to be isometric. In addition, we will introduce (electric) dyadic Green's function, which leads to the uniqueness of metric. Finally, by calculating any order (normal and tangent) derivatives of $\mu$ and $\sigma$ on $\Gamma$, we get Theorem 1.4, and hence by using analyticity up to $\Gamma$ we can obtain Theorem 1.5.

This paper is organized as follows: in section 2, we give the local expression of the electromagnetic Dirichlet-to-Neumann map $\Lambda_{g,\Gamma}$. In section 3, we explicitly calculate the full symbol of $\Lambda_{g,\Gamma}$.
Therefore, Proposition 2.6 and the proof of Proposition 3.1 have provided the detail proof for Theorem 1.1. Furthermore, we show that the full symbol of $\Lambda_{g,\Gamma}$ uniquely determines $g_{\alpha\beta}$ and their all-order derivatives on $\Gamma$; hence, by using the (electric) dyadic Green's function we will prove Theorem 1.2 in section 4. Finally, for given metric $g$, Theorem 1.4 and the uniqueness of the anisotropic electromagnetic parameters $\mu$ and $\sigma$ from $\Lambda_{g,\Gamma}$ are proved in section 5.

\vskip 1.49 true cm

\section{A factorization of Maxwell's equations on Riemannian manifold}

\vskip 0.45 true cm

Let $M$ be an oriented compact Riemannian $3$-manifold endowed with a Riemannian metric $g$.
 By $TM$ (respectively, $T^*(M)$) we denote the tangent (respectively, cotangent) bundle on $M$, which are the disjoint union of the tangent (respectively, cotangent) spaces at all points of $M$:
    \begin{eqnarray*} TM=\Pi_{p\in M} T_p M \quad  \; \;\; (\mbox{respectively}, T^* M = \Pi_{p\in M} T^*_p M).\end{eqnarray*} A vector field on $M$ is a section of the map $\pi: TM\to M$.  More concretely, a vector field is a smooth map $X: M\to TM$, usually written $p\mapsto X_p$, with the property that
\begin{eqnarray*} \pi \circ X={\mbox{Id}}_M,\end{eqnarray*}
or equivalently, $X_p \in T_pM$ for each $p\in M$. If $(U; x_1,x_2,x_3)$ is any smooth coordinate chart for $M$, we can write the value of $X$ at any point $p\in U$ in terms of the coordinate basis vectors $\{\frac{\partial}{\partial x_j}\big|_{p}\}$ of $T_pM$:
\begin{eqnarray*}  X_p= \sum_{j=1}^3 X^j (p) \, \frac{ \partial}{\partial x_j}\bigg|_{p}.\end{eqnarray*}
This defines three functions $X^j:U\to {\mathbb{ R}}$, called the component functions of $X$ in the given chart.
It is standard to use the notation $\mathfrak{X}(M)$ to denote the set of all smooth vector fields on $M$.
Let $\nabla$ be the associated Levi-Civita connection. For a Riemannian manifold $(M,g)$, the connection coefficients (i.e., Christoffel symbols) associated with the metric $g$ are given by (see, for example, \cite{Ta2}) \begin{eqnarray} \label{19.8.3-11} \Gamma_{lk}^j= \frac{1}{2} \sum\limits_{m=1}^3 g^{jm} \bigg( \frac{\partial g_{km}}{\partial x_l} +\frac{\partial g_{lm}}{\partial x_k} -\frac{\partial g_{lk}}{\partial x_m}\bigg),\end{eqnarray}
 where $\big[g^{jk}\big]_{3\times 3}$ is the inverse of $g=\big[g_{jk}\big]_{3\times 3}$.
It is well-known that in a local coordinate system with the naturally associated frame field on the tangent bundle,
\begin{eqnarray*} \nabla_{\frac{\partial}{\partial x_k}} X = \sum\limits_{j=1}^n \big(\frac{\partial X^j}{\partial x_k} +\sum\limits_{l=1}^n \Gamma_{lk}^j X^l  \big)\frac{\partial }{\partial x_j}\quad \; \mbox{for}\;\; X=\sum\limits_{j=1}^n X^j \frac{\partial}{\partial x_j}. \end{eqnarray*}
  If we denote \begin{eqnarray*} {X^j}_{;k}= \frac{\partial X^j}{\partial x_k} +\sum\limits_{l=1}^n \Gamma_{lk}^j X^l,\end{eqnarray*}
 then     \begin{eqnarray*} \label{18/10/29}  \nabla_Y X = \sum\limits_{j,k=1}^{n} Y^k {X^j}_{;k} \,\frac{\partial}{\partial x_j} \;\; \mbox{for}\,\; Y=\sum\limits_{k=1}^n Y^{k} \frac{\partial}{\partial x_k}.\end{eqnarray*}
 The associated Riemannian curvature tensor is:
 $$R(X,Y)Z =[\nabla_X, \nabla_Y] Z- \nabla_{[X,Y]} Z.$$
 In a local coordinates, the expression for the Riemannian curvature tensor is \begin{eqnarray} \label{19.8.29-1} R^j_{klm} = \frac{\partial \Gamma_{km}^j}{\partial x_l} -\frac{\partial \Gamma_{kl}^j}{\partial x_m} + \sum_{s=1}^3\big(\Gamma_{sl}^j \Gamma_{km}^s -\Gamma_{sm}^j \Gamma_{kl}^s\big),\end{eqnarray}
 and the Ricci curvature tensor is defined by \begin{eqnarray} \label{19.8.29-2} R_{km}=\sum_{j=1}^3 R^j_{kjm}.\end{eqnarray}
 It follows from p.$\,$410 of \cite{Youn} that if $E$ and $F$ are two vectors, i.e.,  \begin{eqnarray*}E=\sum_{j=1}^3 E^j \frac{\partial}{\partial x_j}\; \;\, \mbox{and}\;\,\;  F= \sum_{j=1}^3 F^j \frac{\partial}{\partial x_j},\end{eqnarray*} then the vector product $E\times  F$ of $E$ and $F$ is
\begin{eqnarray} && E\times F= \sqrt{|g|} \left\{  \begin{vmatrix}
g^{11} & g^{21} & g^{31}\\
 E^1 & E^2 &E^3\\
 F^1&F^2&F^3\end{vmatrix}  \frac{\partial }{\partial x_1} +
 \begin{vmatrix}
g^{12} & g^{22} & g^{32}\\
 E^1 & E^2 &E^3\\
 F^1&F^2&F^3\end{vmatrix} \frac{\partial }{\partial x_2}
+ \begin{vmatrix}
g^{13} & g^{23} & g^{33}\\
 E^1 & E^2 &E^3\\
 F^1&F^2&F^3\end{vmatrix}  \frac{\partial }{\partial x_3}\right\},\end{eqnarray}
where $|a|$ denotes the determinant of the matrix \begin{eqnarray*}  a=\begin{bmatrix}a^{11} & a^{12} & a^{13}\\
a^{21} & a^{22}& a^{23}\\
a^{31} & a^{32} & a^{33}\end{bmatrix}. \end{eqnarray*}
That is,
\begin{eqnarray} \label{2020.4.12-1}&& E\times F= \sqrt{|g|}\, \begin{vmatrix}
 dx_1  & dx_2 & dx_3\\
 E^1 & E^2 &E^3\\
 F^1&F^2&F^3\end{vmatrix} \\
&&\qquad \quad  = \sqrt{|g|} \left\{  \begin{vmatrix}
 E^2 &E^3\\
 F^2&F^3
\end{vmatrix}  dx_1 +    \begin{vmatrix}
 E^3 &E^1\\
 F^3&F^1\end{vmatrix} dx_2 +   \begin{vmatrix}
 E^1 &E^2\\
 F^1& F^2\end{vmatrix} dx_3\right\}\nonumber
   \end{eqnarray}
with  $dx_j=\sum_{k=1}^3 g^{jk}\frac{\partial }{\partial x_k}, \;\, j=1,2,3$.
On the other hand, according to p.$\,$410 of \cite{Youn} we can also write $E\times F$ as the following form:
\begin{eqnarray} \label{20.4.11-1}\!\! \!\!\!\!\!&\!\!\!& \;\;\;\;\; \;\,\;\;\quad \;\quad  E\times F=\frac{1}{\sqrt{|g|}}   \begin{vmatrix}
 \frac{\partial }{\partial x_1}   &  \frac{\partial }{\partial x_2} &  \frac{\partial }{\partial x_3}\\
 \sum_{l=1}^3 g_{1l} E^l  & \sum_{l=1}^3 g_{2l} E^l &\sum_{l=1}^3 g_{3l} E^l\\
 \sum_{l=1}^3 g_{1l} F^l&\sum_{l=1}^3 g_{2l} F^l&\sum_{l=1}^3 g_{3l} F^l\end{vmatrix}\\
\!\!\!\!\!&\!\!\!&\!\!\!=\! \frac{1}{\!\sqrt{|g|}} \!\left\{ \!  \begin{vmatrix}
    \sum_{\!l=1}^3 g_{2l} E^l &\sum_{\
    \!l=1}^3 g_{3l} E^l\\
 \sum_{l=1}^3 g_{2l} F^l&\sum_{\!l=1}^3 g_{3l} F^l\end{vmatrix}  \!\frac{\partial }{\partial x_1} \!+\!   \begin{vmatrix}
    \sum_{\!l=1}^3 g_{3l} E^l &\sum_{\!l=1}^3 g_{1l} E^l\\
 \sum_{\!l=1}^3 g_{3l} F^l&\sum_{l=1}^3 g_{1l} F^l\end{vmatrix}\!  \frac{\partial }{\partial x_2} \!+\!\begin{vmatrix}
    \sum_{\!l=1}^3 g_{1l} E^l &\sum_{\!l=1}^3 g_{2l} E^l\\
 \sum_{\!l=1}^3 g_{1l} F^l&\sum_{\!l=1}^3 g_{2l} F^l\end{vmatrix}  \!\frac{\partial }{\partial x_3}\!\right\}\!. \nonumber \end{eqnarray}
In fact, from (\ref{2020.4.12-1}) we  have
\begin{eqnarray} \label{2020.4.12-2} \;\;\;\,\quad \, \;  E\times F=\!\!\!\!\!\!\!&&\!\!\! \sqrt{|g|}\, \begin{vmatrix}
 \sum_{k=1}^3 g^{1k} \frac{\partial }{\partial x_k}   & \sum_{k=1}^3 g^{2k} \frac{\partial }{\partial x_k} &\sum_{k=1}^3 g^{3k} \frac{\partial }{\partial x_k}\\
 E^1 & E^2 &E^3\\
 F^1&F^2&F^3\end{vmatrix} \\
 =\!\!\!\!\!\!\!&&\!\!\!  \sqrt{|g|}\begin{vmatrix}
 \sum_{k=1}^3 g^{1k} \frac{\partial }{\partial x_k}   & \sum_{k=1}^3 g^{2k} \frac{\partial }{\partial x_k} &\sum_{k=1}^3 g^{3k} \frac{\partial }{\partial x_k}\\
 \sum_{k=1}^3 g^{1k} \big(\sum_{l=1}^3 g_{kl} E^l\big)   &  \sum_{k=1}^3 g^{2k} \big(\sum_{l=1}^3 g_{kl} E^l\big)  & \sum_{k=1}^3 g^{3k} \big(\sum_{l=1}^3 g_{kl} E^l\big) \\
 \sum_{k=1}^3 g^{1k} \big(\sum_{l=1}^3 g_{kl} F^l\big)   &  \sum_{k=1}^3 g^{2k} \big(\sum_{l=1}^3 g_{kl} F^l\big)  & \sum_{k=1}^3 g^{3k} \big(\sum_{l=1}^3 g_{kl} F^l\big)\end{vmatrix} \nonumber \\
 =\!\!\!\!\!\!\!&&\!\!\!  \sqrt{|g|}  \begin{vmatrix}
 \sum_{k=1}^3 g^{1k} \frac{\partial }{\partial x_k}   &  \sum_{k=1}^3 g^{1k} \big(\sum_{l=1}^3 g_{kl} E^l\big) & \sum_{k=1}^3 g^{1k} \big(\sum_{l=1}^3 g_{kl} F^l\big) \\  \sum_{k=1}^3 g^{2k} \frac{\partial }{\partial x_k} &  \sum_{k=1}^3 g^{2k} \big(\sum_{l=1}^3 g_{kl} E^l\big)  &  \sum_{k=1}^3 g^{2k} \big(\sum_{l=1}^3 g_{kl} F^l\big)  \\  \sum_{k=1}^3 g^{3k} \frac{\partial }{\partial x_k}  &
  \sum_{k=1}^3 g^{3k} \big(\sum_{l=1}^3 g_{kl} E^l\big)   &   \sum_{k=1}^3 g^{3k} \big(\sum_{l=1}^3 g_{kl} F^l\big)\end{vmatrix} \nonumber\\
       =\!\!\!\!\!\!\!&&\!\!\!  \sqrt{|g|} \left| \begin{bmatrix}  g^{11} & g^{12} & g^{13}\\ g^{21} & g^{22} & g^{23}\\
    g^{31} & g^{32} & g^{33}\end{bmatrix} \begin{bmatrix} \frac{\partial }{\partial x_1}  & \sum_{l=1}^3 g_{1l} E^l & \sum_{l=1}^3 g_{1l} F^l\\
    \frac{\partial }{\partial x_2}  & \sum_{l=1}^3 g_{2l} E^l & \sum_{l=1}^3 g_{2l} F^l\\
    \frac{\partial }{\partial x_3}  & \sum_{l=1}^3 g_{3l} E^l & \sum_{l=1}^3 g_{3l} F^l
    \end{bmatrix}\right| \nonumber\\
     =\!\!\!\!\!\!\!&&\!\!\!  \frac{1}{\sqrt{|g|}}  \begin{vmatrix} \frac{\partial }{\partial x_1} &  \frac{\partial}{\partial x_2} & \frac{\partial} {\partial x_3}\\  \sum_{l=1}^3 g_{1l} E^l &  \sum_{l=1}^3 g_{2l} E^l &  \sum_{l=1}^3 g_{3l} E^l \\
     \sum_{l=1}^3 g_{1l} F^l &  \sum_{l=1}^3 g_{2l} F^l &  \sum_{l=1}^3 g_{3l} F^l \end{vmatrix}, \nonumber
    \end{eqnarray}
so (\ref{20.4.11-1}) is verified.

For any smooth real-valued function $f$ on a Riemannian manifold $(M,g)$, the {\it gradient} of $f$ is defined by
  \begin{eqnarray*} \mbox{grad}\, f= (df)^\#,\end{eqnarray*}
 and the {\it divergence operator} $\mbox{div}: \mathfrak{X}(M)\to C^\infty (M)$ is defined by
 \begin{eqnarray*} \mbox{div}\, X =  *d*X^\flat, \end{eqnarray*}
 where $*: \wedge^k T^*M \to \wedge^{3-k} T^* M$ is the Hodge star operator, and $\flat$ and $\#$ are {\it flat} and {\it sharp} operators by lowering index and raising index, respectively.
  In smooth local coordinates, the $\mbox{grad}\, f$ and $\mbox{div}\, X$  have the following expression
 \begin{eqnarray*} \mbox{grad}\, f= \sum\limits_{j,k=1}^3 g^{jk}\frac{\partial f}{\partial x_k}\,\frac{\partial }{\partial x_j},\end{eqnarray*}
\begin{eqnarray}\label{19.9.13-2} \mbox{div} \, X= \sum_{j=1}^3\frac{1}{\sqrt{|g|}} \frac{\partial}{\partial x_j} \big( \sqrt{|g|}\, X^j) \quad\, \mbox{for}\;\, X=\sum\limits_{j=1}^3 X^j \frac{\partial }{\partial x_j}\in \mathfrak{X}(M).\end{eqnarray}
  Accordingly, the Laplace-Beltrami
operator $\Delta_g$ is just given by
 \begin{eqnarray} \label{19.7.6-5} \Delta_g:= \mbox{div}\; \mbox{grad} = \frac{1}{\sqrt{|g|}} \sum_{j,k=1}^3 \frac{\partial}{\partial x_j} \bigg(\sqrt{g}\,g^{jk} \frac{\partial}{\partial x_k}\bigg). \end{eqnarray}
  The {\it curl operator}, denoted by
 $\mbox{curl}: \mathfrak{X}(M) \to     \mathfrak{X}(M)$, is defined by
\begin{eqnarray*} & \mbox{curl}\, X= (* d(X^\flat))^\# \;\;\; \mbox{for}\,\; X\in  \mathfrak{X}(M).\end{eqnarray*}
The curl operator defined only in dimension $3$ because it is only in that case that $\wedge^2 T^* M$ is isomorphic to $TM$ (via the Hodge star operator). Suppose that $X=\sum_{j=1}^3 X^j \frac{\partial}{\partial x_j}\in \mathfrak{X}(M)$ in smooth coordinates, then by definition of $\mbox{curl}\, X$, we have (see p.$\,$454 of \cite{Youn})
\begin{eqnarray}\label{19.7.28-1}  &&\quad\; \mbox{curl}\, X\!=\! \nabla\times X\!=\! \frac{1}{\sqrt{|g|}}\left\{\! \bigg(\frac{\partial}{\partial x_2} \big(\sum_{l=1}^3 g_{3l}X^l\big) \!-\!  \frac{\partial} {\partial x_3}\big(\sum_{l=1}^2 g_{2l}X^l\big)\bigg)\frac{\partial }{\partial x_1} \right.\\
&& \left. +\bigg(\!\frac{\partial}{\partial x_3}\! \big(\sum_{l=1}^3 g_{1l}X^l\big) \!- \! \frac{\partial }{\partial x_1}\big(\sum_{l=1}^3 g_{3l}X^l\big)\!\bigg)\frac{\partial }{\partial x_2}\!+\! \bigg(\!\frac{\partial }{\partial x_1}\big(\sum_{l=1}^3 g_{2l}X^l\big)\! -\!  \frac{\partial }{\partial x_2}\big(\sum_{l=1}^3 g_{1l}X^l\big)\!\bigg)\frac{\partial }{\partial x_3}\!\right\}\!\in \!\mathfrak{X}(M).\nonumber\end{eqnarray}
The operators $\mbox{div}$, $\mbox{grad}$ and $\mbox{curl}$ on an oriented three-dimensional Riemannian manifold  $M$ are related by the following commutative diagram:
  \[\begin{matrix}
C^\infty(M) &
\stackrel{\mathrm{grad}}{\longrightarrow} & \mathfrak{X}(M) & \stackrel{\mbox{curl}}{\longrightarrow} & \mathfrak{X}(M) & \stackrel{\mbox{div}}{\longrightarrow} &  C^\infty (M)\\
 \Big\downarrow\vcenter{%
\rlap{$\scriptstyle{\mbox{Id}}$}} &  &\Big\downarrow\vcenter{%
\rlap{$\flat$}} & &  \Big\downarrow\vcenter{%
\rlap{$*$}}& & \Big\downarrow\vcenter{%
\rlap{$*$}}
 \\
 \wedge^0 (M)& \stackrel{d}{\longrightarrow} &\wedge^1 (M) & \stackrel{d}{\longrightarrow} &
 \wedge^2 (M) & \stackrel{d}{\longrightarrow}& \wedge^3 (M).
\end{matrix}\]
The identities $\mbox{curl} \, \mbox{grad}\equiv 0$ and $\mbox{div}\, \mbox{curl}\equiv 0$ follow from $dd\equiv 0$ and $**\eta =(-1)^{k(n-k)} \eta$ if $\eta$ is a $k$-form on $n$-dimensional Riemannian manifolds.

\vskip 0.28 true cm

 \noindent{\bf Lemma 2.1.} \ {\it Let $M$ be an oriented compact Riemannian $3$-manifold with smooth metric tensor $g$. If  $X=\sum_{j=1}^3 X^j\frac{\partial }{\partial x_j}$ is a smooth vector field on $M$, then in local coordinates,}
\begin{eqnarray}
\label{19.8.24-1}  &&\mbox{curl}\, \mbox{curl}\, X= \mbox{grad}\, \mbox{div}\, X - \sum\limits_{j=1}^3 \bigg(\Delta_g X^j + 2 \sum\limits_{k,l,m=1}^3 g^{ml} \Gamma_{km}^j\frac{\partial X^k}{\partial x_l}   + \sum\limits_{k,l,m=1}^3 g^{ml} \frac{\partial \Gamma_{km}^j}{\partial x_l} X^k \\
 && \qquad \qquad \,\quad \;\; +\sum\limits_{k,l,m,h=1}^3 g^{ml} \Gamma_{hl}^j\Gamma_{km}^h X^k- \sum\limits_{k,l,m,h=1}^3g^{ml} \Gamma_{kh}^j \Gamma_{ml}^h X^k - \sum\limits_{k=1}^3 R_k^j X^k \bigg)\frac{\partial}{\partial  x_j},
\nonumber\end{eqnarray}
where  $R^j_k=\sum\limits_{m=1}^3 g^{jm}R_{mk}$.
\vskip 0.25 true cm

 \noindent  {\it Proof.} \  In local coordinates, let $X^\flat= \eta$ be the dual $1$-form of the vector field $X= \sum_{j=1}^3 X^j\frac{\partial }{\partial x_j}$, i.e., $\eta=\sum_{j=1}^3X _j dx_j$, where $X_j =\sum_{k=1}^3 g_{jk} X^k$.
Then for such a $1$-form $\eta$, we have
 \begin{eqnarray*} &(* d)(*d) \eta = (*d*) d \eta,\end{eqnarray*}
  i.e., \begin{eqnarray*}
   & (\mbox{curl}\, \mbox{curl}\, X)^\flat =  \delta d \eta,\end{eqnarray*}
    where $\delta: \wedge^k(M)\to \wedge^{k-1}(M)$ is defined by $\delta \zeta= (-1)^{3(k+1)+1}( * d * )\zeta$ for any $k$-form $\zeta$ with $0\le k\le 3=\mbox{dim}\, M$.
    Recalling that the Hodge Laplacian $\Delta_{{}_H} = d\delta +\delta d$, we find that
    \begin{align}  \label{19.9.20-1}  (\mbox{curl}\, \mbox{curl}\, X)^\flat &= - d\delta \eta + \Delta_{{}_H} \eta \\
    & = d(* d*) \eta +\Delta_{{}_H}\eta \nonumber\\
    &=  d(\mbox{div}\, X) +\Delta_{{}_H} \eta.\nonumber\end{align}
 It is easy to see that
\begin{align*} \Delta_{{}_H} \eta &=d\delta \eta +\delta d\eta = -\sum\limits_{j,m=1}^3 X^m_{\;\;\;;m;j} dx_j +\sum_{j,l,m=1}^3 \big(g^{ml} X_{m;j;l}-g^{ml} X_{j;m;l}\big)dx_j\\
&   = -\sum\limits_{j,l,m=1}^3 g^{ml} X_{j;m;l} dx_j -\sum\limits_{j,m=1}^3\big(X^m_{\;\;\;; m;j}-X^m_{\;\;\;; j;m}\big) dx_j\nonumber \\
&    = - \sum\limits_{j,l,m=1}^3 g^{ml} X_{j;m;l} dx_j+ \sum\limits_{j,k,m=1}^3 R^{m}_{kmj}X^{k}dx_j\nonumber\\
&   = -\sum\limits_{j,l,m=1}^3g^{ml} X_{j;m;l} dx_j +\sum\limits_{j,k=1}^3R_{kj}X^k dx_j,\nonumber\end{align*}
where $X^m_{\;\;\,\, ; m; j}= (X^m_{\,\;\; ;m})_{;j}$, $X_{m; j; l }= (X_{m; j})_{;l}$ and $X_{m;j}= \frac{\partial X_m}{\partial x_j} -\sum_{l=1}^3
\Gamma_{mj}^l X_l$.
By using the sharp operator $\#$ (i.e., by raising an index) and by noting that $g^{js}_{\;\;\,\,,m}=0$ for Riemannian metric $g$,
we get
  \begin{align*} \big(\Delta_{{}_H} X^\flat\big)^{\#} &=-\sum\limits_{j,l,m=1}^3 g^{ml} X^j_{\;\; ;m;l} \frac{\partial }{\partial x_j} +\sum\limits_{j,k,m=1}^3 g^{jm} R_{km}X^k \frac{\partial }{\partial x_j}\\
  &= -\sum\limits_{j,l,m=1}^3 g^{ml} X^j_{\;\; ; m;l} \frac{\partial }{\partial x_j} +\sum\limits_{j,k=1}^3 R^j_k X^k \frac{\partial }{\partial x_j}.\end{align*}
   In view of \begin{align*}  X^j_{\;\,;m;l}& = \frac{\partial X^j_{\;\;;m}}{\partial x_l}+ \sum\limits_{k=1}^3 X^k_{\;\;;m} \Gamma^j_{kl} -\sum\limits_{k=1}^3 X^j_{\;\;;k} \Gamma_{ml}^k\\
  & =\frac{\partial }{\partial x_l} \bigg( \frac{\partial X^j}{\partial x_m}
+\sum\limits_{k=1}^3 \Gamma^j_{km}  X^k\bigg) +\sum\limits_{k=1}^3\bigg( \frac{\partial X^k}{\partial x_m} +\sum\limits_{h=1}^3\Gamma_{hm}^kX^h\bigg) \Gamma_{kl}^j -\sum\limits_{k=1}^3 \bigg(\frac{\partial X^j}{\partial x_k}+\sum\limits_{h=1}^3 \Gamma_{hk}^jX^h \bigg) \Gamma_{ml}^k\\
&= \frac{\partial^2 X^j}{\partial x_m\partial x_l} \!+ \!\sum\limits_{k=1}^3\!\Gamma_{km}^j \frac{\partial X^k}{\partial x_l} \!+\!\sum\limits_{k=1}^3\!\Gamma_{kl}^j \frac{\partial X^k}{\partial x_m} \!- \!\sum\limits_{k=1}^3\! \Gamma_{ml}^{k} \frac{\partial X^j}{\partial x_k}\!+\!\sum\limits_{k=1}^3\! \bigg( \!\frac{\partial \Gamma_{km}^j}{\partial x_l}
\!+\!\sum\limits_{h=1}^3\!\Gamma_{hl}^j\Gamma_{km}^h\!-\! \sum\limits_{h=1}^3\! \Gamma_{kh}^j \Gamma_{ml}^h \bigg)X^k,\end{align*}
we have \begin{align*}\sum\limits_{l,m=1}^3 g^{ml}X^{j}_{\;\; ;m;l}&  =\sum\limits_{l,m=1}^3\bigg(g^{ml} \frac{\partial^2 X^j}{\partial x_m\partial x_l}
 - \sum\limits_{k=1}^3g^{ml} \Gamma_{ml}^k \frac{\partial X^j}{\partial X_k} \bigg)  + \sum\limits_{k,l,m=1}^3\bigg(g^{ml}\Gamma_{km}^j \frac{\partial X^k}{\partial x_l} +g^{ml}\Gamma_{kl}^j \frac{\partial X^k}{\partial x_m}\bigg) \\
 &\quad \, +\sum\limits_{k,l,m=1}^3 \big(g^{ml} \frac{\partial \Gamma_{km}^j}{\partial x_l}  +\sum\limits_{h=1}^3g^{ml}\Gamma_{hl}^j\Gamma_{km}^h   -\sum\limits_{h=1}^3g^{ml} \Gamma_{kh}^j \Gamma_{ml}^h\big) X^k\\
& =   \Delta_g X^j + 2 \sum\limits_{k,l,m=1}^3g^{ml} \Gamma_{km}^j\frac{\partial X^k}{\partial x_l} +\sum\limits_{k,l,m=1}^3 \big(g^{ml} \frac{\partial \Gamma_{km}^j}{\partial x_l}  +\sum\limits_{h=1}^3g^{ml}\Gamma_{hl}^j\Gamma_{km}^h   -\sum\limits_{h=1}^3g^{ml} \Gamma_{kh}^j \Gamma_{ml}^h\big) X^k  \end{align*}
because of $\sum_{l,m=1}^3 \big( g^{ml} \frac{\partial^2 X^j}{\partial x_m\partial x_l} -\sum_{k=1}^3 g^{ml} \Gamma_{ml}^k \frac{\partial X^j}{\partial x_k} \big) = \Delta_g X^j$.
 It follows that
 \begin{eqnarray}\label{18.8.24-7}  && \big(\Delta_{{}_H} X^\flat\big)^{\#}=\sum\limits_{j=1}^3 \bigg\{ -  \Delta_g X^j - 2\sum\limits_{k,l,m=1}^3 g^{ml} \Gamma_{km}^j\frac{\partial X^k}{\partial x_l} - \sum\limits_{k,l,m=1}^3\bigg(g^{ml} \frac{\partial \Gamma_{km}^j}{\partial x_l}  \\
 && \qquad \qquad \quad \;\; +\sum\limits_{h=1}^3g^{ml} \Gamma_{hl}^j\Gamma_{km}^h -\sum\limits_{h=1}^3g^{ml} \Gamma_{kh}^j \Gamma_{ml}^h\bigg) X^k + \sum\limits_{k=1}^3 R^j_{k} X^k\bigg\} \frac{\partial}{\partial x_j}.\nonumber\end{eqnarray}
Noting that $\mbox{grad}\, \mbox{div}\, X =\big(d \,(\mbox{div}\, X)\big)^\sharp$, we find by (\ref{19.9.20-1}) that
\begin{eqnarray*} &&\mbox{curl}\, \mbox{curl}\, X= \mbox{grad}\, \mbox{div}\, X - \sum\limits_{j=1}^3 \bigg(\Delta_g X^j + 2 \sum\limits_{k,l,m=1}^3g^{ml} \Gamma_{km}^j\frac{\partial X^k}{\partial x_l}   + \sum\limits_{m,k,l=1}^3g^{ml} \frac{\partial \Gamma_{km}^j}{\partial x_l} X^k \\
 && \qquad \qquad \,\quad \;\; +\sum\limits_{k,l,m,h=1}^3g^{ml}\Gamma_{hl}^j \Gamma_{km}^h X^k-\sum\limits_{k,l,m,h=1}^3 g^{ml} \Gamma_{kh}^j \Gamma_{ml}^h X^k - \sum\limits_{k=1}^3R^j_{k} X^k \bigg)\frac{\partial}{\partial  x_j}.\nonumber\end{eqnarray*} \qed

\vskip 0.25 true cm

\noindent{\bf Remark 2.2.} \  {\it The formula (\ref{19.8.24-1}) can also be obtained by a direct calculation from the definition of $\mbox{curl}$ if we using the relations $gg^{-1}=I_3$, $\frac{\partial g_{jk}}{\partial x_l}= \sum\limits_{m=1}^3\big(g_{jm}\Gamma_{kl}^m +g_{km}\Gamma^m_{jl}\big)$ and $\frac{\partial g^{jk}}{\partial x_l}= -\sum\limits_{m=1}^3\big(g^{jm}\Gamma^k_{ml}+g^{km}\Gamma_{ml}^j\big).$}

\vskip 0.28 true cm

 \noindent{\bf Lemma 2.3.} \ {\it Let $(E,H)$ be a solution of the Maxwell equations (\ref{19.7.6-1}). Then  \begin{eqnarray}\label{19.9.6-1}
 \mbox{div}\,(\sigma E)=\mbox{div}\, (\mu H)=0 \quad \mbox{in} \,\, M.\end{eqnarray}}

 \vskip 0.25 true cm

 \noindent  {\it Proof.} \  As pointed out before, $\mbox{div}\, \mbox{curl}\, X=0$ for any vector field $X\in \mathfrak{X}(M)$. In fact, this also can directly be got as follows.  Applying (\ref{19.7.28-1}) and (\ref{19.9.13-2}) we have
  \begin{align*}   &            \mbox{div}\, \mbox{curl}\, X =
   \frac{1}{\sqrt{|g|}}  \bigg\{\big( \frac{\partial^2 (g_{3l} X^l)}{\partial x_1\partial x_2} -   \frac{\partial^2(g_{2l} X^l)}{\partial x_1\partial x_3}
\big) +
  \big( \frac{\partial^2(g_{1l} X^l)}{\partial x_2\partial x_3} -   \frac{\partial^2(g_{3l} X^l)}{\partial x_2\partial x_1}
\big)  +  \big( \frac{\partial^2(g_{2l} X^l)}{\partial x_1\partial x_3} -   \frac{\partial^2(g_{1l} X^l)}{\partial x_2\partial x_3} \big)\bigg\}
 \\
 &       \;\;  -\frac{1}{2\sqrt{|g|^3}} \bigg\{\! \big(\frac{\partial |g|}{\partial x_1}\big)  \big( \frac{\partial(g_{3l} X^l)}{\partial x_2} \!-   \frac{\partial(g_{2l} X^l)}{\partial x_3} \big) + \big(\frac{\partial |g|}{\partial x_2}\big)  \big( \frac{\partial(g_{1l} X^l)}{\partial x_3} -   \frac{\partial(g_{3l} X^l)}{\partial x_1} \big) \!  + \big(\frac{\partial |g|}{\partial x_3}\big)  \big( \frac{\partial(g_{2l} X^l)}{\partial x_1} -   \frac{\partial(g_{1l} X^l)}{\partial x_2} \big) \!\bigg\}\\
 &         \;\; + \Gamma_{1k}^k  \big( \frac{\partial(g_{3l} X^l)}{\partial x_2} -   \frac{\partial(g_{2l} X^l)}{\partial x_3} \big)  + \Gamma_{2k}^k  \big( \frac{\partial(g_{1l} X^l)}{\partial x_3} -   \frac{\partial(g_{3l} X^l)}{\partial x_1} \big) +  \Gamma_{3k}^k \big( \frac{\partial(g_{2l} X^l)}{\partial x_1} -   \frac{\partial(g_{1l} X^l)}{\partial x_2} \big)=0\end{align*}
because of $ \frac{1}{2|g|} \, \frac{\partial |g|}{\partial x_k} =\sum_{l=1}^3\Gamma_{kl}^l$ for $k=1,2,3$.
Taking $\mbox{div}$ to Maxwell's equations and applying  the above conclusion, we immediately get the desired result. \qed

\vskip 0.22 true cm

By Lemma 2.3, we have
 \begin{eqnarray} \label{2020.3.29-2} 0= \mbox{div}\, (\sigma E) =\sigma\, \mbox{div}\, E +\sum\limits_{l=1}^3 \frac{\partial \sigma}{\partial x_l} E^l,\end{eqnarray} so that $$\mbox{div}\, E =- \frac{1}{\sigma} \sum\limits_{l=1}^3 \frac{\partial \sigma}{\partial x_l} E^l.$$
Thus \begin{eqnarray}\label{19.12.10-1}  \mbox{grad} \, \mbox{div}\, E = - \sum\limits_{j=1}^3 \bigg(\sum\limits_{m=1}^3 g^{jm} \frac{\partial }{\partial x_m} \Big(\frac{1}{\sigma} \sum\limits_{l=1}^3 \frac{\partial \sigma}{\partial x_l} E^l\Big)\bigg) \frac{\partial }{\partial x_j}. \end{eqnarray}

Applying Maxwell's equations (\ref{19.7.6-1}) again, we obtain  \begin{eqnarray} \label{19.8.12-1} \mbox{curl curl} \,E-i\omega \,\mbox{curl}\,(\mu H)=0 \quad \mbox{in}\;\; M.\end{eqnarray}
Since $H(x)= (H_1(x), H_2(x), H_3(x))$, we get
\begin{align} \label{2020.4.11-2} & \mbox{curl}\, (\mu H) = \mbox{curl}\,( \mu H_1(x),  \mu H_2(x),  \mu H_3(x))= \frac{1}{\sqrt{|g|} } \bigg\{ \bigg(\frac{\partial}{\partial x_2}\big( \sum\limits_{l=1}^3 g_{3l} (\mu H^l)\big) - \frac{\partial}{\partial x_3}\big( \sum\limits_{l=1}^3 g_{2l} (\mu H^l)\big)\bigg)\frac{\partial}{\partial x_1}\\
 & + \bigg(\frac{\partial}{\partial x_3}\big( \sum\limits_{l=1}^3 g_{1l} ( \mu H^l)\big) - \frac{\partial}{\partial x_1}\big( \sum\limits_{l=1}^3 g_{3l} (\mu H^l)\big)\bigg)\frac{\partial}{\partial x_2}+  \bigg(\frac{\partial}{\partial x_1}\big( \sum\limits_{l=1}^3 g_{2l} (\mu H^l)\big) - \frac{\partial}{\partial x_2}\big( \sum\limits_{l=1}^3 g_{1l} (\mu H^l)\big)\bigg)\frac{\partial}{\partial x_3}\bigg\}\nonumber\\
 & =\mu \, \mbox{curl}\, H +\frac{1}{\sqrt{|g|} } \bigg\{   \bigg(\frac{\partial \mu}{\partial x_2} \sum\limits_{l=1}^3 g_{3l}  H^l - \frac{\partial \mu}{\partial x_3}\sum\limits_{l=1}^3 g_{2l} H^l\bigg)\frac{\partial}{\partial x_1}\nonumber\\
 & + \bigg(\frac{\partial \mu }{\partial x_3} \sum\limits_{l=1}^3 g_{1l}  H^l - \frac{\partial \mu}{\partial x_1} \sum\limits_{l=1}^3 g_{3l}  H^l\bigg)\frac{\partial}{\partial x_2}+  \bigg(\frac{\partial \mu }{\partial x_1}\sum\limits_{l=1}^3 g_{2l} H^l - \frac{\partial \mu}{\partial x_2} \sum\limits_{l=1}^3 g_{1l}  H^l\bigg)\frac{\partial}{\partial x_3}\bigg\}\nonumber\\
 &=\mu \, \mbox{curl}\, H + (\mbox{grad}\, \mu) \times H.\nonumber
\end{align}
Here, the second term of the last equality follows from the definition of vector product (see (\ref{20.4.11-1})). Combining (\ref{19.8.12-1}), (\ref{2020.4.11-2}) and Maxwell's equations, we get
\begin{eqnarray} \label{2020.4.11-4} \mbox{curl}\;\mbox{curl}\; E -\omega^2 \mu \sigma\, E -  \big(\mbox{grad}\, (\log \mu)\big) \times \mbox{curl}\; E=0 \quad \mbox{in}\;\, M.\end{eqnarray}
Because \begin{align*} & H= \frac{1}{i\omega \mu} \mbox{curl}\, E =\frac{1}{i\omega\mu\sqrt{|g|} } \bigg\{
  \bigg(\frac{\partial}{\partial x_2}\big( \sum\limits_{s=1}^3 g_{3s}  E^s\big) - \frac{\partial}{\partial x_3}\big( \sum\limits_{s=1}^3 g_{2s} E^s\big)\bigg)\frac{\partial}{\partial x_1}+ \bigg(\frac{\partial}{\partial x_3}\big( \sum\limits_{s=1}^3 g_{1s}  E^s\big) \\
  & \qquad - \frac{\partial}{\partial x_1}\big( \sum\limits_{s=1}^3 g_{3s}  E^s\big)\bigg)\frac{\partial}{\partial x_2} +  \bigg(\frac{\partial}{\partial x_1}\big( \sum\limits_{s=1}^3 g_{2s} E^s\big) - \frac{\partial}{\partial x_2}\big( \sum\limits_{s=1}^3 g_{1s}  E^s\big)\bigg)\frac{\partial}{\partial x_3}\bigg\},\end{align*}
  we have  \begin{align*} & \frac{1}{\sqrt{|g|} } \bigg\{   \bigg(\frac{\partial \mu}{\partial x_2} \sum\limits_{l=1}^3 g_{3l}  H^l - \frac{\partial \mu}{\partial x_3}\sum\limits_{l=1}^3 g_{2l} H^l\bigg)\frac{\partial}{\partial x_1}+ \bigg(\frac{\partial \mu }{\partial x_3} \sum\limits_{l=1}^3 g_{1l}  H^l - \frac{\partial \mu}{\partial x_1} \sum\limits_{l=1}^3 g_{3l}  H^l\bigg)\frac{\partial}{\partial x_2} \\
 & \quad +  \bigg(\frac{\partial \mu }{\partial x_1}\sum\limits_{l=1}^3 g_{2l} H^l - \frac{\partial \mu}{\partial x_2} \sum\limits_{l=1}^3 g_{1l}  H^l\bigg)\frac{\partial}{\partial x_3}\bigg\}\\
 & = \frac{1}{i\omega\mu|g|} \Bigg( \bigg\{ \frac{\partial \mu}{\partial x_2} \bigg[ g_{31} \bigg( \frac{\partial}{\partial x_2}\big(\sum\limits_{s=1}^3 g_{3s} E^s) -\frac{\partial }{\partial x_3} \big(\sum\limits_{s=1}^3 g_{2s} E^s\big) \bigg) +g_{32}  \bigg( \frac{\partial}{\partial x_3}\big(\sum\limits_{s=1}^3 g_{1s} E^s) -\frac{\partial }{\partial x_1} \big(\sum\limits_{s=1}^3 g_{3s} E^s\big) \bigg)
\\ &\quad  +g_{33}  \bigg( \frac{\partial}{\partial x_1}\big(\sum\limits_{s=1}^3 g_{2s} E^s) -\frac{\partial }{\partial x_2} \big(\sum\limits_{s=1}^3 g_{1s} E^s\big) \bigg) \bigg] - \frac{\partial \mu}{\partial x_3} \bigg[ g_{21} \bigg( \frac{\partial}{\partial x_2}\big(\sum\limits_{s=1}^3 g_{3s} E^s) -\frac{\partial }{\partial x_3} \big(\sum\limits_{s=1}^3 g_{2s} E^s\big) \bigg) \\
& \quad  +g_{22}  \bigg( \frac{\partial}{\partial x_3}\big(\sum\limits_{s=1}^3 g_{1s} E^s) -\frac{\partial }{\partial x_1} \big(\sum\limits_{s=1}^3 g_{3s} E^s\big) \bigg)
 +g_{23}  \bigg( \frac{\partial}{\partial x_1}\big(\sum\limits_{s=1}^3 g_{2s} E^s) -\frac{\partial }{\partial x_2} \big(\sum\limits_{s=1}^3 g_{1s} E^s\big) \bigg) \bigg] \bigg\}\frac{\partial}{\partial x_1} \\
& \quad + \bigg\{ \frac{\partial \mu}{\partial x_3} \bigg[ g_{11} \bigg( \frac{\partial}{\partial x_2}\big(\sum\limits_{s=1}^3 g_{3s} E^s) -\frac{\partial }{\partial x_3} \big(\sum\limits_{s=1}^3 g_{2s} E^s\big) \bigg) +g_{12}  \bigg( \frac{\partial}{\partial x_3}\big(\sum\limits_{s=1}^3 g_{1s} E^s) -\frac{\partial }{\partial x_1} \big(\sum\limits_{s=1}^3 g_{3s} E^s\big) \bigg)
\\ &\quad  +g_{13}  \bigg( \frac{\partial}{\partial x_1}\big(\sum\limits_{s=1}^3 g_{2s} E^s) -\frac{\partial }{\partial x_2} \big(\sum\limits_{s=1}^3 g_{1s} E^s\big) \bigg) \bigg] - \frac{\partial \mu}{\partial x_1} \bigg[ g_{31} \bigg( \frac{\partial}{\partial x_2}\big(\sum\limits_{s=1}^3 g_{3s} E^s) -\frac{\partial }{\partial x_3} \big(\sum\limits_{s=1}^3 g_{2s} E^s\big) \bigg) \\
& \quad  +g_{32}  \bigg( \frac{\partial}{\partial x_3}\big(\sum\limits_{s=1}^3 g_{1s} E^s) -\frac{\partial }{\partial x_1} \big(\sum\limits_{s=1}^3 g_{3s} E^s\big) \bigg)
 +g_{33}  \bigg( \frac{\partial}{\partial x_1}\big(\sum\limits_{s=1}^3 g_{2s} E^s) -\frac{\partial }{\partial x_2} \big(\sum\limits_{s=1}^3 g_{1s} E^s\big) \bigg) \bigg] \bigg\}\frac{\partial}{\partial x_2}\\
  & \quad + \bigg\{ \frac{\partial \mu}{\partial x_1} \bigg[ g_{21} \bigg( \frac{\partial}{\partial x_2}\big(\sum\limits_{s=1}^3 g_{3s} E^s) -\frac{\partial }{\partial x_3} \big(\sum\limits_{s=1}^3 g_{2s} E^s\big) \bigg) +g_{22}  \bigg( \frac{\partial}{\partial x_3}\big(\sum\limits_{s=1}^3 g_{1s} E^s) -\frac{\partial }{\partial x_1} \big(\sum\limits_{s=1}^3 g_{3s} E^s\big) \bigg)
\\ & \quad +g_{23}  \bigg( \frac{\partial}{\partial x_1}\big(\sum\limits_{s=1}^3 g_{2s} E^s) -\frac{\partial }{\partial x_2} \big(\sum\limits_{s=1}^3 g_{1s} E^s\big) \bigg) \bigg] - \frac{\partial \mu}{\partial x_2} \bigg[ g_{11} \bigg( \frac{\partial}{\partial x_2}\big(\sum\limits_{s=1}^3 g_{3s} E^s) -\frac{\partial }{\partial x_2} \big(\sum\limits_{s=1}^3 g_{2s} E^s\big) \bigg) \\
& \quad +g_{12}  \bigg( \frac{\partial}{\partial x_3}\big(\sum\limits_{s=1}^3 g_{1s} E^s) -\frac{\partial }{\partial x_1} \big(\sum\limits_{s=1}^3 g_{3s} E^s\big) \bigg)
 +g_{13}  \bigg( \frac{\partial}{\partial x_1}\big(\sum\limits_{s=1}^3 g_{2s} E^s) -\frac{\partial }{\partial x_2} \big(\sum\limits_{s=1}^3 g_{1s} E^s\big) \bigg) \bigg] \bigg\}\frac{\partial}{\partial x_3}  \Bigg)\\
       & = \frac{1}{i\omega\mu|g|} \Bigg( \bigg\{ \big( \frac{\partial \mu}{\partial x_2} \, g_{31} - \frac{\partial \mu}{\partial x_3} \, g_{21}\big)
  \sum\limits_{s=1}^3 \Big(g_{3s} \frac{\partial}{\partial x_2} - g_{2s} \frac{\partial}{\partial x_3} +  \frac{\partial  g_{3s}} {\partial x_2}  - \frac{\partial g_{2s}}{\partial x_3} \Big)E^s\\
 & \quad \,  + \big( \frac{\partial \mu}{\partial x_2} \, g_{32} - \frac{\partial \mu}{\partial x_3} \, g_{22}\big)
  \sum\limits_{s=1}^3 \Big(g_{1s} \frac{\partial}{\partial x_3} - g_{3s} \frac{\partial}{\partial x_1} + \frac{\partial  g_{1s}} {\partial x_3}  - \frac{\partial g_{3s}}{\partial x_1} \Big)E^s
   \\
     & \quad \,  + \big( \frac{\partial \mu}{\partial x_2} \, g_{33} - \frac{\partial \mu}{\partial x_3} \, g_{23}\big)
  \sum\limits_{s=1}^3 \Big(g_{2s} \frac{\partial}{\partial x_1} - g_{1s} \frac{\partial}{\partial x_2} + \frac{\partial  g_{2s}} {\partial x_1}  - \frac{\partial g_{1s}}{\partial x_2} \Big)E^s \bigg\}\frac{\partial}{\partial x_1} \\
     & \quad \,  +\bigg\{  \big( \frac{\partial \mu}{\partial x_3} \, g_{11} - \frac{\partial \mu}{\partial x_1} \, g_{31}\big)
  \sum\limits_{s=1}^3 \Big(g_{3s} \frac{\partial}{\partial x_2} - g_{2s} \frac{\partial}{\partial x_3} + \frac{\partial  g_{3s}} {\partial x_2}  - \frac{\partial g_{2s}}{\partial x_3} \Big)E^s\\
   & \quad \,  +  \big( \frac{\partial \mu}{\partial x_3} \, g_{12} - \frac{\partial \mu}{\partial x_1} \, g_{32}\big)
  \sum\limits_{s=1}^3 \Big(g_{1s} \frac{\partial}{\partial x_3} - g_{3s} \frac{\partial}{\partial x_1} + \frac{\partial  g_{1s}} {\partial x_3}  - \frac{\partial g_{3s}}{\partial x_1} \Big)E^s\\
   & \quad \,  +  \big( \frac{\partial \mu}{\partial x_3} \, g_{13} - \frac{\partial \mu}{\partial x_1} \, g_{33}\big)
  \sum\limits_{s=1}^3 \Big(g_{2s} \frac{\partial}{\partial x_1} - g_{1s} \frac{\partial}{\partial x_2} + \frac{\partial  g_{2s}} {\partial x_1}  - \frac{\partial g_{1s}}{\partial x_2} \Big)E^s\bigg\}\frac{\partial}{\partial x_2}\\
      & \quad \,  +\bigg\{  \big( \frac{\partial \mu}{\partial x_1} \, g_{21} - \frac{\partial \mu}{\partial x_2} \, g_{11}\big)
  \sum\limits_{s=1}^3 \Big(g_{3s} \frac{\partial}{\partial x_2} - g_{2s} \frac{\partial}{\partial x_3} + \frac{\partial  g_{3s}} {\partial x_2}  - \frac{\partial g_{2s}}{\partial x_3} \Big)E^s\\
   & \quad \,  +  \big( \frac{\partial \mu}{\partial x_1} \, g_{22} - \frac{\partial \mu}{\partial x_2} \, g_{12}\big)
  \sum\limits_{s=1}^3 \Big(g_{1s} \frac{\partial}{\partial x_3} - g_{3s} \frac{\partial}{\partial x_1} + \frac{\partial  g_{1s}} {\partial x_3}  - \frac{\partial g_{3s}}{\partial x_1} \Big)E^s\\
   & \quad \,  +  \big( \frac{\partial \mu}{\partial x_1} \, g_{23} - \frac{\partial \mu}{\partial x_2} \, g_{13}\big)
  \sum\limits_{s=1}^3 \Big(g_{2s} \frac{\partial}{\partial x_1} - g_{1s} \frac{\partial}{\partial x_2} + \frac{\partial  g_{2s}} {\partial x_1}  - \frac{\partial g_{1s}}{\partial x_2} \Big)E^s\bigg\}\frac{\partial}{\partial x_3}\Bigg).  \end{align*}
 Combining Lemma 2.1, Lemma 2.3, (\ref{19.12.10-1}) and (\ref{19.8.12-1}), we obtain
 \begin{eqnarray} \label{2020.2.3-1} &&  \;\;\, \;\;\; -\!\sum\limits_{\!j=1}^3\!\! \bigg( \!\!\Delta_g E^j \!+\!2 \!\sum\limits_{\!k,l,m\!=\!1}^3\! \!\!g^{ml} \Gamma_{\!km}^j \!\frac{\partial E^k}{\partial x_l}\!+\!\!\sum\limits_{\!k,l,m\!=\!1}^3\!\! \!g^{ml} \frac{\partial \Gamma^j_{\!km}}{\partial x_l} E^k \!\!+\!\!\sum\limits_{\!k,l,m,h\!=\!1}^3\! \!\! g^{ml} \Gamma_{\!hl}^j \Gamma_{km}^h    E^k \! \!-\!\! \sum\limits_{\!k,l,m,h\!=\!1}^3 \!\!\!g^{ml} \Gamma^j_{\!kh}\Gamma_{\!ml}^h E^k   \\
   &&  \;  -R^j_k E^k
   +\omega^2 \mu \sigma E^j \bigg) \frac{\partial}{\partial x_j}
       - \sum\limits_{j=1}^3 \bigg\{ \sum\limits_{m=1}^3 g^{jm} \bigg[ \sum\limits_{l=1}^3 \bigg( \frac{\partial}{\partial x_m}
  \big( \frac{1}{\sigma} \,\frac{\partial \sigma}{\partial x_l} \big) \bigg) E^l   +\frac{1}{\sigma}\sum\limits_{l=1}^3 \frac{\partial \sigma}{\partial x_l} \frac{\partial E^l}{\partial x_m}\bigg]\bigg\} \frac{\partial }{\partial  x_j}\nonumber\\
   &&  \;  -\frac{1}{\mu |g|} \Bigg( \bigg\{ \big( \frac{\partial \mu}{\partial x_2} \, g_{31} - \frac{\partial \mu}{\partial x_3} \, g_{21}\big)
  \sum\limits_{s=1}^3 \Big(g_{3s} \frac{\partial}{\partial x_2} - g_{2s} \frac{\partial}{\partial x_3} +  \frac{\partial  g_{3s}} {\partial x_2}  - \frac{\partial g_{2s}}{\partial x_3} \Big)E^s   \nonumber\\
 &&  \; + \big( \frac{\partial \mu}{\partial x_2} \, g_{32} - \frac{\partial \mu}{\partial x_3} \, g_{22}\big)
  \sum\limits_{s=1}^3 \Big(g_{1s} \frac{\partial}{\partial x_3} - g_{3s} \frac{\partial}{\partial x_1} + \frac{\partial  g_{1s}} {\partial x_3}  - \frac{\partial g_{3s}}{\partial x_1} \Big)E^s  \nonumber
   \\
     &&  \; + \big( \frac{\partial \mu}{\partial x_2} \, g_{33} - \frac{\partial \mu}{\partial x_3} \, g_{23}\big)
  \sum\limits_{s=1}^3 \Big(g_{2s} \frac{\partial}{\partial x_1} - g_{1s} \frac{\partial}{\partial x_2} + \frac{\partial  g_{2s}} {\partial x_1}  - \frac{\partial g_{1s}}{\partial x_2} \Big)E^s \bigg\} \frac{\partial}{\partial x_1}  \nonumber\\
     &&  \; +\bigg\{  \big( \frac{\partial \mu}{\partial x_3} \, g_{11} - \frac{\partial \mu}{\partial x_1} \, g_{31}\big)
  \sum\limits_{s=1}^3 \Big(g_{3s} \frac{\partial}{\partial x_2} - g_{2s} \frac{\partial}{\partial x_3} + \frac{\partial  g_{3s}} {\partial x_2}  - \frac{\partial g_{2s}}{\partial x_3} \Big)E^s   \nonumber\\
   && \;  +  \big( \frac{\partial \mu}{\partial x_3} \, g_{12} - \frac{\partial \mu}{\partial x_1} \, g_{32}\big)
  \sum\limits_{s=1}^3 \Big(g_{1s} \frac{\partial}{\partial x_3} - g_{3s} \frac{\partial}{\partial x_1} + \frac{\partial  g_{1s}} {\partial x_3}  - \frac{\partial g_{3s}}{\partial x_1} \Big)E^s   \nonumber\\
   && \;+  \big( \frac{\partial \mu}{\partial x_3} \, g_{13} - \frac{\partial \mu}{\partial x_1} \, g_{33}\big)
  \sum\limits_{s=1}^3 \Big(g_{2s} \frac{\partial}{\partial x_1} - g_{1s} \frac{\partial}{\partial x_2} + \frac{\partial  g_{2s}} {\partial x_1}  - \frac{\partial g_{1s}}{\partial x_2} \Big)E^s\bigg\}\frac{\partial}{\partial x_2} \nonumber\\
      && \; +\bigg\{  \big( \frac{\partial \mu}{\partial x_1} \, g_{21} - \frac{\partial \mu}{\partial x_2} \, g_{11}\big)
  \sum\limits_{s=1}^3 \Big(g_{3s} \frac{\partial}{\partial x_2} - g_{2s} \frac{\partial}{\partial x_3} + \frac{\partial  g_{3s}} {\partial x_2}  - \frac{\partial g_{2s}}{\partial x_3} \Big)E^s \nonumber\\
   &&  \; +  \big( \frac{\partial \mu}{\partial x_1} \, g_{22} - \frac{\partial \mu}{\partial x_2} \, g_{12}\big)
  \sum\limits_{s=1}^3 \Big(g_{1s} \frac{\partial}{\partial x_3} - g_{3s} \frac{\partial}{\partial x_1} + \frac{\partial  g_{1s}} {\partial x_3}  - \frac{\partial g_{3s}}{\partial x_1} \Big)E^s \nonumber\\
   && \;  +  \big( \frac{\partial \mu}{\partial x_1} \, g_{23} - \frac{\partial \mu}{\partial x_2} \, g_{13}\big)
  \sum\limits_{s=1}^3 \Big(g_{2s} \frac{\partial}{\partial x_1} - g_{1s} \frac{\partial}{\partial x_2} + \frac{\partial  g_{2s}} {\partial x_1}  - \frac{\partial g_{1s}}{\partial x_2} \Big)E^s\bigg\}\frac{\partial}{\partial x_3}\Bigg)=0.  \nonumber \end{eqnarray}
\noindent Write the above equation as the form of components relative to coordinates:
\begin{eqnarray}  \label{19.11.26-1} \end{eqnarray}
\begin{align}
 &  \!   {\mathcal{M}}_g E= \left\{ \Delta_g I_3 +  \begin{bmatrix} 2\sum\limits_{l,m=1}^3 g^{ml} {\Gamma_{1m}^1 \frac{\partial }{\partial x_l}}_{{}_{{}_{}}} &  2\sum\limits_{l,m=1}^3 g^{ml} \Gamma_{2m}^1 \frac{\partial }{\partial x_l} &  2\sum\limits_{l,m=1}^3  g^{ml} \Gamma_{3m}^1\frac{\partial }{\partial x_l}\\
    2\sum\limits_{l,m=1}^3 g^{ml} \Gamma_{1m}^2\frac{\partial }{\partial x_l} &  2\sum\limits_{l,m=1}^3 g^{ml} \Gamma_{2m}^2\frac{\partial }{\partial x_l} &  2\sum\limits_{l,m=1}^3 g^{ml} {\Gamma_{3m}^2\frac{\partial }{\partial x_l}}_{{}_{{}_{}}}\\
    2\sum\limits_{l,m=1}^3 g^{ml} \Gamma_{1m}^3 \frac{\partial }{\partial x_l}&  2\sum\limits_{l,m=1}^3 g^{ml} {\Gamma_{2m}^3\frac{\partial }{\partial x_l}}_{{}_{{}_{}}} & 2\sum\limits_{l,m=1}^3 g^{ml} {\Gamma_{3m}^3\frac{\partial }{\partial x_l}}_{{}_{{}_{}}}\end{bmatrix}\right. \nonumber
    \\
 &    \!     +  \begin{bmatrix} a_{11}+\omega^2 \mu\sigma -R_1^{\,1} & a_{12}-R_2^{\,1} & a_{13}-R_3^{\,1}\\
 a_{21}- R_1^{\,2} &a_{22}+\omega^2 \mu\sigma -R_2^{\,2} & a_{23}-R_3^{\,2} \\
    a_{31}-R_1^{\,3} & a_{32}-R_2^{\,3} & a_{33}+ \omega^2 \mu\sigma - R_3^{\,3}\end{bmatrix}\nonumber\\
   & \!   + \!  \begin{bmatrix} \sum\limits_{m\!=\!1}^3 \!g^{1m}\! \Big(\! \frac{\partial}{\partial x_m} \!\big(\frac{1}{\sigma}\frac{\partial \sigma}{\partial x_1}\!\big) \!+\!\frac{1}{\sigma}\! \frac{\partial \sigma}{\partial x_1}\frac{\partial}{\partial x_m}\!\Big)  & \sum\limits_{m\!=\!1}^3 \!g^{1m} \!\Big( \!\frac{\partial}{\partial x_m} \!\big(\!\frac{1}{\sigma}\frac{\partial \sigma}{\partial x_2}\!\big) \!+\!\frac{1}{\sigma} \!\frac{\partial \sigma}{\partial x_2}\frac{\partial}{\partial x_m}\!\Big) &\sum\limits_{m\!=\!1}^3 \!g^{1m} \!\Big(\! \frac{\partial}{\partial x_m}\! \big(\!\frac{1}{\sigma}\frac{\partial \sigma}{\partial x_3}\big)\! +\!\frac{1}{\sigma} \frac{\partial \sigma}{\partial x_3}\frac{\partial}{\partial x_m}\!\Big)\\
    \sum\limits_{m\!=\!1}^3 \!g^{2m} \!\Big( \!\frac{\partial}{\partial x_m} \!\big(\!\frac{1}{\sigma}\frac{\partial \sigma}{\partial x_1}\!\big) \! +\!\frac{1}{\sigma} \!\frac{\partial \sigma}{\partial x_1}\frac{\partial}{\partial x_m}\Big)  & \sum\limits_{m\!=\!1}^3 \!g^{2m} \!\Big(\! \frac{\partial}{\partial x_m} \!\big(\frac{1}{\sigma}\frac{\partial \sigma}{\partial x_2}\big)\! +\!\frac{1}{\sigma} \!\frac{\partial \sigma}{\partial x_2}\frac{\partial}{\partial x_m}\Big) &\sum\limits_{m\!=\!1}^3 \!g^{2m} \!\Big( \!\frac{\partial}{\partial x_m} \big(\frac{1}{\sigma}\frac{\partial \sigma}{\partial x_3}\big) \!+\!\frac{1}{\sigma}\! \frac{\partial \sigma}{\partial x_3}\frac{\partial}{\partial x_m}\Big)\\
     \sum\limits_{m\!=\!1}^3\! g^{3m}\! \Big( \!\frac{\partial}{\partial x_m}\! \big(\frac{1}{\sigma}\frac{\partial \sigma}{\partial x_1}\big) \!+\!\frac{1}{\sigma} \frac{\partial \sigma}{\partial x_1}\frac{\partial}{\partial x_m}\Big)  & \sum\limits_{m\!=\!1}^3 \!g^{3m}\! \Big(\! \frac{\partial}{\partial x_m} \big(\frac{1}{\sigma}\frac{\partial \sigma}{\partial x_2}\big) \!+\!\frac{1}{\sigma} \frac{\partial \sigma}{\partial x_2}\frac{\partial}{\partial x_m}\Big) &\sum\limits_{m\!=\!1}^3 \!g^{3m} \!\Big( \frac{\partial}{\partial x_m} \big(\frac{1}{\sigma}\frac{\partial \sigma}{\partial x_3}\big) \!+\!\frac{1}{\sigma} \frac{\partial \sigma}{\partial x_3}\frac{\partial}{\partial x_m}\Big)\end{bmatrix}\nonumber \\
    \!\!&\! + \!\frac{1}{\mu|g|}\!\! \begin{bmatrix}\! g_{31} \frac{\partial \mu}{\partial x_2} \!-\!\!g_{21} \frac{\partial \mu}{\partial x_3}  & g_{32} \frac{\partial \mu}{\partial x_2} \!-\!\!g_{22} \frac{\partial \mu}{\partial x_3}& g_{33} \frac{\partial \mu}{\partial x_2} \!-\!\!g_{23} \frac{\partial \mu}{\partial x_3} \\  \!
     g_{11} \frac{\partial \mu}{\partial x_3} \!-\!\!g_{31} \frac{\partial \mu}{\partial x_1}  & g_{12} \frac{\partial \mu}{\partial x_3} \!-\!\!g_{32} \frac{\partial \mu}{\partial x_1}& g_{13} \frac{\partial \mu}{\partial x_3} \!-\!g_{33} \frac{\partial \mu}{\partial x_1} \\
     \! g_{21} \frac{\partial \mu}{\partial x_1} \!-\!\!g_{11} \frac{\partial \mu}{\partial x_2}  & g_{22} \frac{\partial \mu}{\partial x_1}\! -\!\!g_{12} \frac{\partial \mu}{\partial x_2}& g_{23} \frac{\partial \mu}{\partial x_1} \!-\!\!g_{13} \frac{\partial \mu}{\partial x_2} \end{bmatrix}\!\!
     \begin{bmatrix} \!g_{31}\! \frac{\partial }{\partial x_2} \!-\!\!g_{21} \frac{\partial }{\partial x_3}  & g_{32} \frac{\partial }{\partial x_2} \!-\!\!g_{22} \frac{\partial }{\partial x_3}& g_{33} \frac{\partial }{\partial x_2} \!-\!\!g_{23} \frac{\partial }{\partial x_3} \!\\
    \! g_{11} \frac{\partial }{\partial x_3} \!-\!g_{31} \frac{\partial }{\partial x_1}  & g_{12} \frac{\partial }{\partial x_3} \!-\!g_{32} \frac{\partial }{\partial x_1}& g_{13} \frac{\partial }{\partial x_3} \!-\!\!g_{33} \frac{\partial }{\partial x_1} \!\\
    \! g_{21} \frac{\partial }{\partial x_1} \!-\!\!g_{11} \frac{\partial }{\partial x_2}  & g_{22} \frac{\partial }{\partial x_1}\! -\!\!g_{12} \frac{\partial}{\partial x_2}& g_{23} \frac{\partial }{\partial x_1} \!-\!\!g_{13} \frac{\partial }{\partial x_2} \!\end{bmatrix} \nonumber\\    &    \! \left. + \!\frac{1}{\mu |g|} \!  \begin{bmatrix}\! g_{31} \frac{\partial \mu}{\partial x_2} \!-\!g_{21} \frac{\partial \mu}{\partial x_3}  \!& g_{32} \frac{\partial \mu}{\partial x_2} \!-\! g_{22} \frac{\partial \mu}{\partial x_3}  \!& g_{33} \frac{\partial \mu}{\partial x_2} \!-\!g_{23} \frac{\partial \mu}{\partial x_3} \!\\
     \! g_{11} \frac{\partial \mu}{\partial x_3} \!-\!g_{31} \frac{\partial \mu}{\partial x_1} \! & g_{12} \frac{\partial \mu}{\partial x_3} \!-\!g_{32} \frac{\partial \mu}{\partial x_1}   \!&  \!g_{13} \frac{\partial \mu}{\partial x_3} \!-\!g_{33} \frac{\partial \mu}{\partial x_1} \!\\
    \! g_{21} \frac{\partial \mu}{\partial x_1} \!-\!g_{11} \frac{\partial \mu}{\partial x_2}  \!&  \!  g_{22} \frac{\partial \mu}{\partial x_1} \!-\!g_{12} \frac{\partial \mu}{\partial x_2}   \!&  \! g_{23} \frac{\partial \mu}{\partial x_1} \!-\!g_{13} \frac{\partial \mu}{\partial x_2} \! \end{bmatrix}               \!\! \begin{bmatrix} \! \frac{\partial g_{31}}{\partial x_2} \!-\! \frac{\partial g_{21}}{\partial x_3}  \! &\frac{\partial  g_{32} }{\partial x_2} \!-\! \frac{\partial g_{22}}{\partial x_3}\! &  \frac{\partial g_{33}}{\partial x_2}\! -\! \frac{\partial g_{23}}{\partial x_3} \\
      \!  \frac{\partial g_{11}}{\partial x_3} \!-\! \frac{\partial g_{31}}{\partial x_1}  \! &  \frac{\partial g_{12}}{\partial x_3} \!-\! \frac{\partial g_{32}}{\partial x_1}   \!  &  \frac{\partial g_{13}}{\partial x_3}\! -\! \frac{\partial g_{33}}{\partial x_1} \\
    \!  \frac{\partial  g_{21}}{\partial x_1}\! - \!\frac{\partial g_{11}}{\partial x_2} \!  &  \frac{\partial g_{22}}{\partial x_1} \!-\! \frac{\partial g_{12}}{\partial x_2}   \! &  \frac{\partial g_{23}}{\partial x_1}\! -\! \frac{\partial g_{13}}{\partial x_2} \!\end{bmatrix}
              \! \right\}  \! \! \begin{bmatrix}\! E^1\\ \! E^2\\ \! E^{3} \!\end{bmatrix}\!=\!0, \nonumber \end{align}
  where \begin{eqnarray} \label{19.9.13-1} a_{jk}= \sum\limits_{l,m}g^{ml} \big(\frac{\partial \Gamma_{km}^j}{\partial x_l} + \sum\limits_{h}\Gamma_{hl}^j \Gamma_{km}^h -\sum\limits_h \Gamma_{kh}^j \Gamma_{ml}^h \big)\quad \mbox{for all}\;\,j,k=1,2,3.\end{eqnarray}
  Clearly, ${\mathcal{M}}_g E=0$ in $M$ is a system of second-order linear elliptic equations.
\vskip 0.18 true cm

 If $x':=(x_1, x_2)$  are any local coordinates for
$\partial M$ near $x_0=0\in \Gamma\subset\partial M$, then one can obtain the boundary normal coordinates for $M$ in some neighborhood of $x_0$ determined by $(x_1, x_2)$. In these coordinates $x_3>0$ in
 $M$, and $\partial M$ is locally characterized by $x_3= 0$. A standard computation shows
that the metric $g$ on $\bar M$ then has the form
 (see p.$\,$532 of \cite{Ta2})
\begin{eqnarray} \label{18/a-1} \; \quad\;\big[g_{jk} (x',x_{3}) \big]_{3\times 3} = \begin{bmatrix}
 g_{11} (x',x_{3}) & g_{12} (x',x_{3}) &   0\\
 g_{21} (x',x_{3})  & g_{22} (x',x_{3}) &  0\\
 0& 0& 1 \end{bmatrix}.  \end{eqnarray}
 By (\ref{18/a-1}) we immediately see that the inverse of metric tensor $g$ in the boundary normal coordinates has form:
    \begin{eqnarray*} g^{-1}(x',x_n) =\begin{bmatrix} g^{11}(x', x_3) &  g^{12} (x', x_3)& 0 \\
        g^{21}(x',x_3) &  g^{22}(x',x_3)& 0\\
    0&0 &1\end{bmatrix}. \end{eqnarray*}
      Under this normal coordinates, we take the outward unit normal vector ${\nu}(x)=[0, 0,-1]^t$, where $A^t$ denotes the transpose of a vector $A$.

 In what follows, we will let Greek indices run from $1$ to $2$ and Roman indices from $1$ to $3$. Let $I_3$ denote the $3 \times  3$ identity
matrix. In boundary normal coordinates, we have that $g^{k3}=g^{3k}=0$ for any $1\le k\le 2$,
 and \begin{eqnarray} \label{19.9.6-2} \quad\;\; \left. \begin{array} {ll} &\qquad \Gamma_{3k}^3 =\frac{1}{2} \sum_{m=1}^3 g^{3m}\bigg(\frac{\partial g_{3m}}{\partial x_k} +\frac{\partial g_{km}}{\partial x_3} -\frac{\partial g_{3k}}{\partial x_m}\bigg)=\frac{1}{2} \bigg(\frac{\partial g_{33}}{\partial x_k}+\frac{\partial g_{k3}}{\partial x_3} -\frac{\partial g_{3k}}{\partial x_n}\bigg)=0,\\
  &\qquad \Gamma_{33}^k =\frac{1}{2}\sum\limits_{m=1}^3 g^{km}\bigg(\frac{\partial g_{3m}}{\partial x_3} +\frac{\partial g_{3m}}{\partial x_3} -\frac{\partial g_{33}}{\partial x_m}\bigg) =0\end{array}\right. \quad \, k=1,2,3.\end{eqnarray}
   Thus, in local boundary normal coordinates, we can rewrite (\ref{19.11.26-1}) as
\begin{align*}
 \!\!&  \!   {\mathcal{M}}_g E= \left\{ \bigg( \frac{\partial^2}{\partial x_3^2} +\big(\frac{1}{2} \sum_{\alpha,\beta} g^{\alpha\beta} \frac{\partial g_{\alpha\beta}}{\partial x_3} \big) \frac{\partial }{\partial x_3} +\sum\limits_{\alpha, \beta} g^{\alpha \beta} \frac{\partial^2}{\partial x_\alpha\partial x_\beta} +\sum\limits_{\alpha, \beta} \big(g^{\alpha\beta}\sum_{\gamma} \Gamma_{\alpha \gamma}^\gamma  +\frac{\partial g^{\alpha\beta}}{\partial x_\alpha}\big)\frac{\partial}{\partial x_\beta}\bigg)I_3  \right.\\
\!\!&  \!      + \begin{bmatrix} 2\Gamma_{13}^1 & {2\Gamma_{23}^1}_{{}_{{}_{}}} &0\\
   2\Gamma_{13}^2 & {2\Gamma_{23}^2}_{{}_{{}_{}}}  & 0\\
   0& 0& 0\end{bmatrix} \big(\frac{\partial }{\partial x_3}\, I_3\big)+ \begin{bmatrix} 2\sum\limits_{\alpha,\beta}g^{\alpha\beta} {\Gamma_{1\alpha}^1 \frac{\partial }{\partial x_\beta}}_{{}_{{}_{}}} &  2\sum\limits_{\alpha,\beta}g^{\alpha\beta} \Gamma_{2\alpha}^1 \frac{\partial }{\partial x_\beta}&  2\sum\limits_{\alpha,\beta}g^{\alpha\beta} \Gamma_{3\alpha}^1\frac{\partial }{\partial x_\beta}\\
    2\sum\limits_{\alpha,\beta}g^{\alpha\beta} \Gamma_{1\alpha}^2\frac{\partial }{\partial x_\beta} &  2\sum\limits_{\alpha,\beta}g^{\alpha\beta} \Gamma_{2\alpha}^2\frac{\partial }{\partial x_\beta} &  2\sum\limits_{\alpha,\beta}g^{\alpha\beta} {\Gamma_{3\alpha}^2\frac{\partial }{\partial x_\beta}}_{{}_{{}_{}}}\\
    2\sum\limits_{\alpha,\beta}g^{\alpha\beta} \Gamma_{1\alpha}^3 \frac{\partial }{\partial x_\beta}&  2\sum\limits_{\alpha,\beta}g^{\alpha\beta} {\Gamma_{2\alpha}^3\frac{\partial }{\partial x_\beta}}_{{}_{{}_{}}} &  0\end{bmatrix} \nonumber
    \\
  \!\!&    \!     +  \begin{bmatrix} a_{11}+\omega^2 \mu\sigma -R_1^{\,1} & a_{12}-R_2^{\,1} & a_{13}-R_3^{\,1}\\
 a_{21}- R_1^{\,2} &a_{22}+\omega^2 \mu\sigma -R_2^{\,2} & a_{23}-R_3^{\,2} \\
    a_{31}-R_1^{\,3} & a_{32}-R_2^{\,3} & a_{33}+ \omega^2 \mu\sigma - R_3^{\,3}\end{bmatrix}+ \begin{bmatrix} 0& 0& 0 \\
    0& 0& 0 \\
    \frac{1}{\sigma} \frac{\partial \sigma}{\partial x_1} & \frac{1}{\sigma} \frac{\partial \sigma}{\partial x_2} &\frac{1}{\sigma} \frac{\partial \sigma}{\partial x_3}\end{bmatrix} \big(\frac{\partial }{\partial x_3} I_3\big)\nonumber\\
  \!\! & \!   + \!  \begin{bmatrix} \sum\limits_{\beta} \!g^{1\beta}\! \Big(\! \frac{\partial}{\partial x_\beta} \!\big(\frac{1}{\sigma}\frac{\partial \sigma}{\partial x_1}\!\big) \!+\!\frac{1}{\sigma}\! \frac{\partial \sigma}{\partial x_1}\frac{\partial}{\partial x_\beta}\!\Big)  & \sum\limits_{\beta} \!g^{1\beta} \!\Big( \!\frac{\partial}{\partial x_\beta} \!\big(\!\frac{1}{\sigma}\frac{\partial \sigma}{\partial x_2}\!\big) \!+\!\frac{1}{\sigma} \!\frac{\partial \sigma}{\partial x_2}\frac{\partial}{\partial x_\beta}\!\Big) &\sum\limits_{\beta} \!g^{1\beta} \!\Big(\! \frac{\partial}{\partial x_\beta}\! \big(\!\frac{1}{\sigma}\frac{\partial \sigma}{\partial x_3}\big)\! +\!\frac{1}{\sigma} \frac{\partial \sigma}{\partial x_3}\frac{\partial}{\partial x_\beta}\!\Big)\\
    \sum\limits_{\beta} \!g^{2\beta} \!\Big( \!\frac{\partial}{\partial x_\beta} \!\big(\!\frac{1}{\sigma}\frac{\partial \sigma}{\partial x_1}\!\big) \! +\!\frac{1}{\sigma} \!\frac{\partial \sigma}{\partial x_1}\frac{\partial}{\partial x_\beta}\Big)  & \sum\limits_{\beta} \!g^{2\beta} \!\Big(\! \frac{\partial}{\partial x_\beta} \!\big(\frac{1}{\sigma}\frac{\partial \sigma}{\partial x_2}\big)\! +\!\frac{1}{\sigma} \!\frac{\partial \sigma}{\partial x_2}\frac{\partial}{\partial x_\beta}\Big) &\sum\limits_{\beta} \!g^{2\beta} \!\Big( \!\frac{\partial}{\partial x_\beta} \big(\frac{1}{\sigma}\frac{\partial \sigma}{\partial x_3}\big) \!+\!\frac{1}{\sigma}\! \frac{\partial \sigma}{\partial x_3}\frac{\partial}{\partial x_\beta}\Big)\\
       \!\frac{\partial}{\partial x_3}\! \big(\frac{1}{\sigma}\frac{\partial \sigma}{\partial x_1}\big)   &  \frac{\partial}{\partial x_3} \big(\frac{1}{\sigma}\frac{\partial \sigma}{\partial x_2}\big)  &\frac{\partial}{\partial x_3} \big(\frac{1}{\sigma}\frac{\partial \sigma}{\partial x_3}\big)\end{bmatrix}\nonumber \\
       \!\!&\! + \!\frac{1}{\mu|g|}\!\! \begin{bmatrix} \frac{\partial \mu}{\partial x_3} (g_{21} g_{21} -g_{22} g_{11}) &
   0 & 0\\
   0 & \frac{\partial \mu}{\partial x_3} (-g_{11} g_{22} +g_{12}g_{12}) & 0\\
    \frac{\partial \mu}{\partial x_1} (\!-g_{21} g_{21} +g_{11}g_{22}) &   \frac{\partial \mu}{\partial x_2} (g_{11}g_{22} -g_{12} g_{12}) & 0\end{bmatrix} \! \big(\!\frac{\partial }{\partial x_3}I_3\big) \\
           \!\!&\! + \!\frac{1}{\mu|g|}\!\! \begin{bmatrix} \frac{\partial \mu}{\partial x_2} \big( g_{21}\frac{\partial }{\partial x_1} -g_{11}
       \frac{\partial}{\partial x_2}\big)  &  \frac{\partial \mu}{\partial x_2} \big( g_{22}\frac{\partial }{\partial x_1} -g_{12}
       \frac{\partial}{\partial x_2}\big)  &  \frac{\partial\mu}{\partial x_3} ( g_{22} \frac{\partial }{\partial x_1}  -g_{21} \frac{\partial}{\partial x_2}) \\
       - \frac{\partial \mu}{\partial x_1} \big( g_{21}\frac{\partial }{\partial x_1} -g_{11}
       \frac{\partial}{\partial x_2}\big)  &  -\frac{\partial \mu}{\partial x_1} \big( g_{22}\frac{\partial }{\partial x_1} -g_{12}
       \frac{\partial}{\partial x_2}\big)  &    \frac{\partial \mu}{\partial x_3} \big( \!-  g_{12}\frac{\partial}{\partial x_1} + g_{11} \frac{\partial \mu}{\partial x_2} \big) \\
0& 0& \big(\!\!-\! g_{22}\frac{\partial \mu}{\partial x_1}+ \frac{\partial \mu}{\partial x_2} \big)\frac{\partial}{\partial x_1} +  \big( g_{21}\frac{\partial \mu}{\partial x_1}\!-\! g_{11} \frac{\partial \mu}{\partial x_2} \big)\frac{\partial}{\partial x_2}\end{bmatrix}
           \nonumber
          \\     \!\! & \! \left. + \!\frac{1}{\mu |g|} \!  \begin{bmatrix}- g_{21} \frac{\partial \mu}{\partial x_3}  \!& - g_{22} \frac{\partial \mu}{\partial x_3}  \!& \frac{\partial \mu}{\partial x_2} \!\\
     \! g_{11} \frac{\partial \mu}{\partial x_3}  & g_{12} \frac{\partial \mu}{\partial x_3}    \!&  - \frac{\partial \mu}{\partial x_1} \!\\
    \! g_{21} \frac{\partial \mu}{\partial x_1} -g_{11} \frac{\partial \mu}{\partial x_2}  \!&  \!  g_{22} \frac{\partial \mu}{\partial x_1} -g_{12} \frac{\partial \mu}{\partial x_2}   \!&   0 \! \end{bmatrix}
              \!\! \begin{bmatrix} \!- \frac{\partial g_{21}}{\partial x_3}  \! & - \frac{\partial g_{22}}{\partial x_3}\! &  0 \\
      \!  \frac{\partial g_{11}}{\partial x_3} \! &  \frac{\partial g_{12}}{\partial x_3}  \!  &  0 \\
    \!  \frac{\partial  g_{21}}{\partial x_1} - \frac{\partial g_{11}}{\partial x_2} \!  &  \frac{\partial g_{22}}{\partial x_1} - \frac{\partial g_{12}}{\partial x_2}   \! &  0 \!\end{bmatrix}
            \!  \!\right\}  \! \! \begin{bmatrix}\! E^1\\ \!E^2\\ \! E^{3} \!\end{bmatrix}\!=\!0. \nonumber \end{align*}
   Here, for the sake of simplicity, we have used the relationships that \begin{eqnarray} \label{19.9.1-1} \frac{1}{\sqrt{|g|}} \frac{\partial \sqrt{|g|}}{\partial x_k} =\frac{1}{2} \sum\limits_{\alpha, \beta} g^{\alpha\beta} \frac{\partial g_{\alpha\beta}}{\partial x_k}=\sum\limits_{\gamma} \Gamma_{k \gamma}^\gamma.\end{eqnarray}
Consequently, we have
\begin{eqnarray} \label{19.7.7-2,1}   {\mathcal{M}}_g E=\left\{\left(\frac{\partial^2}{\partial x_3^2}\, I_3\right) +B \left(\frac{\partial}{\partial x_3}I_3\right)+C\right\}
  \begin{bmatrix}E^1\\
   E^2\\
 E^3\end{bmatrix}=0, \end{eqnarray}
  where \begin{align*} &         B:=\frac{1}{2} \sum\limits_{\alpha,\beta} g^{\alpha\beta} \frac{\partial g_{\alpha\beta}}{\partial x_3}\,I_3+
  \begin{bmatrix} 2\Gamma_{13}^1 & 2{\Gamma_{23}^1}_{{}_{{}_{}}} &0\\
   2\Gamma_{13}^2 & 2{\Gamma_{23}^2}_{{}_{{}_{}}}  & 0\\
   0& 0& 0\end{bmatrix} +
         \begin{bmatrix} 0 & 0 &0\\
   0 & 0 & 0\\
   \frac{1}{\sigma} \frac{\partial \sigma}{\partial x_1} & \frac{1}{\sigma} \frac{\partial \sigma}{\partial x_2}& \frac{1}{\sigma} \frac{\partial \sigma}{\partial x_3}\end{bmatrix}+ \frac{1}{\mu} \begin{bmatrix} -\frac{\partial \mu}{\partial x_3} & 0& 0\\
   0 & -\frac{\partial \mu}{\partial x_3} & 0\\
    \frac{\partial \mu}{\partial x_1} & \frac{\partial \mu}{\partial x_2} & 0\end{bmatrix},\\
   &      C:= \!\sum\limits_{\alpha, \beta}\! \bigg(\! g^{\alpha\beta} \frac{\partial^2}{\partial x_\alpha \partial x_\beta}\!+ \!\big(g^{\alpha\beta} \sum\limits_{\gamma} \Gamma_{\alpha \gamma}^{\gamma} \!+\!\frac{\partial g_{\alpha\beta}}{\partial x_\alpha}\big)\frac{\partial }{\partial x_\beta}\bigg)I_3   \! +\! \begin{bmatrix} 2\!\sum\limits_{\alpha,\beta}\!g^{\alpha\beta} \Gamma_{1\alpha}^1{\frac{\partial}{\partial x_\beta}}_{{}_{{}_{}}} &  2\!\sum\limits_{\alpha,\beta}\!g^{\alpha\beta} \Gamma_{2\alpha}^1\frac{\partial}{\partial x_\beta} &  2\!\sum\limits_{\alpha,\beta}\!g^{\alpha\beta} {\Gamma_{3\alpha}^1\frac{\partial}{\partial x_\beta}}_{{}_{{}_{}}}\!\\
    2\!\sum\limits_{\alpha,\beta}\!g^{\alpha\beta} \Gamma_{1\alpha}^2\frac{\partial}{\partial x_\beta} &  2\!\sum\limits_{\alpha,\beta}\!g^{\alpha\beta} \Gamma_{2\alpha}^2 \frac{\partial}{\partial x_\beta}&  2\!\sum\limits_{\alpha,\beta}\!g^{\alpha\beta} {\Gamma_{3\alpha}^2\frac{\partial}{\partial x_\beta}}_{{}_{{}_{}}}\!\\
    2\!\sum\limits_{\alpha,\beta}\!g^{\alpha\beta} \Gamma_{1\alpha}^3\frac{\partial}{\partial x_\beta} &  2\!\sum\limits_{\alpha,\beta}\!g^{\alpha\beta} {\Gamma_{2\alpha}^3\frac{\partial}{\partial x_\beta}}_{{}_{{}_{}}} &  0\!\end{bmatrix}\\
    \\
     &  \qquad   \!+\!\begin{bmatrix} \sum\limits_{\beta}\! g^{1\beta} \!\big(\frac{1}{\sigma} \!\frac{\partial \sigma}{\partial x_1}\big)\! \frac{\partial }{\partial x_\beta} \! &  \sum\limits_{\beta}\! g^{1\beta} \!\big(\frac{1}{\sigma} \!\frac{\partial \sigma}{\partial x_2}\big)\! \frac{\partial }{\partial x_\beta}
\!& \sum\limits_{\beta} \!g^{1\beta} \!\big(\frac{1}{\sigma} \!\frac{\partial \sigma}{\partial x_3}\big)\! \frac{\partial }{\partial x_\beta}\\
\sum\limits_{\beta} \!g^{2\beta}\! \big(\frac{1}{\sigma}\! \frac{\partial \sigma}{\partial x_1}\big)\! \frac{\partial }{\partial x_\beta} \! &  \sum\limits_{\beta} g^{2\beta} \!\big(\frac{1}{\sigma} \!\frac{\partial \sigma}{\partial x_2}\big)\! \frac{\partial }{\partial x_\beta}
\! & \sum\limits_{\beta} \!g^{2\beta}\! \big(\frac{1}{\sigma}\! \frac{\partial \sigma}{\partial x_3}\big)\! \frac{\partial }{\partial x_\beta}\\
0 & 0
& 0\end{bmatrix} \\
\\
&      \!+\! \frac{1}{\mu|g|}\! \begin{bmatrix} \frac{\partial \mu}{\partial x_2} \big( g_{21} \frac{\partial }{\partial x_1} \!-\!g_{11} \frac{\partial }{\partial x_2}\big) & \frac{\partial \mu}{\partial x_2} \big( g_{22} \frac{\partial }{\partial x_1}\! -\!g_{21} \frac{\partial }{\partial x_2}\big)  & \!\frac{\partial \mu}{\partial x_3} \big( g_{22} \frac{\partial }{\partial x_1}\! -\!g_{21} \frac{\partial }{\partial x_2}\big)\\
     -\frac{\partial \mu}{\partial x_1} \big( g_{21} \frac{\partial }{\partial x_1} \!-\!g_{11} \frac{\partial }{\partial x_2}\big) & -\frac{\partial \mu}{\partial x_1} \!\big( g_{22} \frac{\partial }{\partial x_1}\! -\!g_{12} \frac{\partial }{\partial x_2}\big) & \frac{\partial \mu}{\partial x_3} \!\big(\!- g_{12} \frac{\partial }{\partial x_1} \!+\!g_{11} \frac{\partial }{\partial x_2}\big)\\
     0& 0& \!\big(\!\!-\!g_{22}\frac{\partial \mu}{\partial x_1}\! + \!\frac{\partial \mu}{\partial x_2}\! \big) \frac{\partial}{\partial x_1}\! +\!\big( g_{21}\frac{\partial \mu}{\partial x_1}\!-\! g_{11} \frac{\partial \mu}{\partial x_2}\big) \frac{\partial}{\partial x_2} \end{bmatrix}\\
     \\
   &  \! + \! \begin{bmatrix} \!a_{11}\!+\!\omega^2 \mu\sigma \!-\!R_1^{1} \!\!&\! a_{12}\!-\!R_2^{1} \!\!&\! a_{13}\!-\!R_3^{1}\\
\! a_{21}\!- \!R_1^{2}\! \!&\!a_{22}\!+\!\omega^2 \mu\sigma\!-\!R_2^{2}\! \!& \!a_{23}\!-\!R_3^{2} \\
  \!  a_{31}\!-\!R_1^{3}\! \!&\! a_{32}\!-\!R_2^{3} \!\! &\! a_{33}\!+\! \omega^2 \mu\sigma \!- \!R_3^{3}\end{bmatrix}
  \!+\!\begin{bmatrix} \!\sum\limits_{\beta}\! g^{1\beta} \!\frac{\partial}{\partial x_\beta} \!\big(\frac{1}{\sigma} \!\frac{\partial \sigma}{\partial x_1}\big) \! \!&  \sum\limits_{\beta}\! g^{1\beta} \! \frac{\partial }{\partial x_\beta}\big(\frac{1}{\sigma} \!\frac{\partial \sigma}{\partial x_2}\big)
\!\!& \sum\limits_{\beta} \!g^{1\beta} \! \frac{\partial }{\partial x_\beta}\!\big(\frac{1}{\sigma} \!\frac{\partial \sigma}{\partial x_3}\big)\! \\
\!\sum\limits_{\beta} \!g^{2\beta}\! \frac{\partial }{\partial x_\beta} \!\big(\frac{1}{\sigma}\! \frac{\partial \sigma}{\partial x_1}\big) \!\! &  \sum\limits_{\beta}\! g^{2\beta} \!\frac{\partial }{\partial x_\beta}\!\big(\frac{1}{\sigma} \!\frac{\partial \sigma}{\partial x_2}\big)\!
\! \!& \sum\limits_{\beta} \!g^{2\beta}\! \frac{\partial }{\partial x_\beta}\!\big(\frac{1}{\sigma}\! \frac{\partial \sigma}{\partial x_3}\big)\! \\
 \! \frac{\partial }{\partial x_3} \!\big(\frac{1}{\sigma}\! \frac{\partial \sigma}{\partial x_1}\big) \! &  \frac{\partial }{\partial x_3}\!\big(\frac{1}{\sigma} \!\frac{\partial \sigma}{\partial x_2}\big)
\! &  \frac{\partial }{\partial x_3}\!\big(\frac{1}{\sigma}\! \frac{\partial \sigma}{\partial x_3}\big)\!\end{bmatrix}\\
\\
& +\frac{1}{\mu|g|}\! \begin{bmatrix}- g_{21} \frac{\partial \mu}{\partial x_3} & -g_{22} \frac{\partial \mu}{\partial x_3} & \frac{\partial \mu}{\partial x_2} \\
g_{11} \frac{\partial \mu}{\partial x_3} & g_{12} \frac{\partial \mu}{\partial x_3} & -\frac{\partial \mu}{\partial x_1}\\
 g_{21} \frac{\partial \mu}{\partial x_1}\!-\!g_{11} \frac{\partial \mu}{\partial x_2}  & g_{22} \frac{\partial \mu}{\partial x_1}\!-\! g_{12}
 \frac{\partial \mu}{\partial x_2} & 0\end{bmatrix} \begin{bmatrix}
 -\frac{\partial g_{21}}{\partial x_3} &   -\frac{\partial g_{22}}{\partial x_3} & 0\\
 \frac{\partial g_{11}}{\partial x_3} &   \frac{\partial g_{12}}{\partial x_3} & 0\\
 \frac{\partial g_{21}}{\partial x_1}\!-\! \frac{\partial g_{11}}{\partial x_2} &    \frac{\partial g_{22}}{\partial x_1}\!-\! \frac{\partial g_{12}}{\partial x_2} & 0\end{bmatrix}. \nonumber  \end{align*}

\noindent Throughout this paper, we denote $\sqrt{-1}=i$.
\vskip 0.28 true cm

 \noindent{\bf Proposition 2.4.} \ {\it There exists a pseudodifferential operator $\Phi(x, D_{x'})$ of order one in $x'$ depending smoothly on $x_3$ such that \begin{eqnarray} \label{19.3.19-1:}  {\mathcal{M}}_g = \left( \frac{\partial}{\partial x_3} \,I_3 + B -\Phi \right) \left(\frac{\partial }{\partial x_3}\,I_3  +\Phi\right), \end{eqnarray}
 modulo a smoothing operator, where $D_{x'}=(D_{x_1}, \cdots, D_{x_{n-1}})$, $\,D_{x_j}=\frac{1}{i}\, \frac{\partial }{\partial x_j}$.}

 \vskip 0.25 true cm

 \noindent  {\it Proof.} \    Let us assume that we have a factorization \begin{eqnarray*} {\mathcal{M}}_g =\left(\frac{\partial}{\partial x_3} \,I_3 +B - \Phi\right)\left( \frac{\partial}{\partial x_3}\,I_3+\Phi\right).\end{eqnarray*}
  Then \begin{eqnarray*} {\mathcal{M}}_g E=\left(\frac{\partial}{\partial x_3} \,I_3 +B - \Phi\right)\left( \frac{\partial}{\partial x_3}\,I_3+\Phi\right)E=0\quad \mbox{for every solution} \;\, E \,\, \mbox{of Maxwell's equations}.\end{eqnarray*}
 It follows this and (\ref{19.7.7-2,1})  that
   \begin{eqnarray} \label{19.3.19-3}  0&=&  {\mathcal{M}}_g- \left(\frac{\partial}{\partial x_3} I_3 +B - \Phi\right)\left( \frac{\partial}{\partial x_3}I_3+\Phi\right)\\
     &=&\Phi^2+ \left(  \Phi\big(\frac{\partial}{\partial x_3}I_3 \big)-   \big(\frac{\partial}{\partial x_3} I_3\big)\Phi  \right) -B\Phi+C.\nonumber \end{eqnarray}
   Let $\phi(x, \xi')$, $b(x,\xi')$ and $c(x,\xi')$ be the full symbols of $\Phi$ and $B$ and $C$, respectively. Clearly,  $\phi(x, \xi') \sim \sum_{j\ge 0} \phi_{1-j} (x, \xi')$,   $\;b(x,\xi')=b_0(x,\xi')$ and $c(x, \xi')= c_2(x,\xi') +c_1(x, \xi') +c_0(x, \xi')$, where
  \begin{align}   &              \label{19.6.10-2}   b_0(x, \xi') =\frac{1}{2} \sum\limits_{\alpha,\beta} g^{\alpha\beta} \frac{\partial g_{\alpha\beta}}{\partial x_3}\,I_3+
  \begin{bmatrix} 2\Gamma_{13}^1 & 2{\Gamma_{23}^1}_{{}_{{}_{}}} &0\\
   2\Gamma_{13}^2 & 2{\Gamma_{23}^2}_{{}_{{}_{}}}  & 0\\
   0& 0& 0\end{bmatrix} +
         \begin{bmatrix} 0 & 0 &0\\
   0 & 0 & 0\\
   \frac{1}{\sigma} \frac{\partial \sigma}{\partial x_1} & \frac{1}{\sigma} \frac{\partial \sigma}{\partial x_2}& \frac{1}{\sigma} \frac{\partial \sigma}{\partial x_3}\end{bmatrix}+ \frac{1}{\mu} \begin{bmatrix} -\frac{\partial \mu}{\partial x_3} & 0& 0\\
   0 & -\frac{\partial \mu}{\partial x_3} & 0\\
    \frac{\partial \mu}{\partial x_1} & \frac{\partial \mu}{\partial x_2} & 0\end{bmatrix},\\ \nonumber \\
    &  c_2(x, \xi') = -\sum_{\alpha, \beta} g^{\alpha\beta} \xi_\alpha\xi_\beta\, I_3, \;\; \nonumber \\
   \nonumber \\
      & c_1(x, \xi') =i \sum\limits_{\alpha,\beta} \big( g^{\alpha\beta}\sum\limits_\gamma \Gamma_{\alpha \gamma}^\gamma +\frac{\partial g^{\alpha\beta}}{\partial x_\alpha} \big) \xi_\beta \, I_3 + \begin{bmatrix} 2i\sum\limits_{\alpha,\beta}g^{\alpha\beta} \Gamma_{1\alpha}^1{\xi_\beta}_{{}_{{}_{}}}&  2i\sum\limits_{\alpha,\beta}g^{\alpha\beta} \Gamma_{2\alpha}^1{\xi_\beta}_{{}_{{}_{}}}  &  2i\sum\limits_{\alpha,\beta}g^{\alpha\beta} \Gamma_{3\alpha}^1{\xi_\beta}_{{}_{}} \\
    2i\sum\limits_{\alpha,\beta}g^{\alpha\beta} {\Gamma_{1\alpha}^2{\xi_\beta}}_{{}_{{}_{}}}  &  2i\sum\limits_{\alpha,\beta}g^{\alpha\beta} \Gamma_{2\alpha}^2 {\xi_\beta}_{{}_{}} &  2i\sum\limits_{\alpha,\beta}g^{\alpha\beta} \Gamma_{3\alpha}^2{\xi_\beta}_{{}_{}} \\
    2i\sum\limits_{\alpha,\beta}g^{\alpha\beta} {\Gamma_{1\alpha}^3{\xi_\beta}}_{{}_{{}_{}}}  &  2i\sum\limits_{\alpha,\beta}g^{\alpha\beta} \Gamma_{2\alpha}^3 {\xi_\beta}_{{}_{}}  &  0\end{bmatrix},\nonumber\\
    \nonumber \\
           &    \!+\!\begin{bmatrix} i \sum\limits_{\beta}\! g^{1\beta} \!\big(\frac{1}{\sigma} \!\frac{\partial \sigma}{\partial x_1}\big) \xi_\beta \! & i\sum\limits_{\beta}\! g^{1\beta} \!\big(\frac{1}{\sigma} \!\frac{\partial \sigma}{\partial x_2}\big) \xi_\beta
\!& i\sum\limits_{\beta} \!g^{1\beta} \!\big(\frac{1}{\sigma} \!\frac{\partial \sigma}{\partial x_3}\big)  \xi_\beta\\
i\sum\limits_{\beta} \!g^{2\beta}\! \big(\frac{1}{\sigma}\! \frac{\partial \sigma}{\partial x_1}\big) \xi_\beta \! &  i\sum\limits_{\beta} g^{2\beta} \!\big(\frac{1}{\sigma} \!\frac{\partial \sigma}{\partial x_2}\big) \xi_\beta
\! & i\sum\limits_{\beta} \!g^{2\beta}\! \big(\frac{1}{\sigma}\! \frac{\partial \sigma}{\partial x_3}\big) \xi_\beta\\
0 & 0
& 0\end{bmatrix} \nonumber\\
\nonumber\\
&      \!+\! \frac{1}{\mu|g|}\! \begin{bmatrix} i\frac{\partial \mu}{\partial x_2} \big( g_{21} \xi_1 \!-\!g_{11} \xi_2\big) & i\frac{\partial \mu}{\partial x_2} \big( g_{22} \xi_1 -g_{12} \xi_2\big)  & i\frac{\partial \mu}{\partial x_3} \big( g_{22} \xi_1\! -\!g_{21} \xi_2\big)\\
     -i\frac{\partial \mu}{\partial x_1} \big( g_{21} \xi_1 \!-\!g_{11} \xi_2\big) & -i\frac{\partial \mu}{\partial x_1} \!\big( g_{22} \xi_1\! -\!g_{12} \xi_2\big) & i\frac{\partial \mu}{\partial x_3} \!\big(\!- g_{12} \xi_1 \!+\!g_{11} \xi_2\big)\\
     0& 0& \!i\big(\!\!-\!g_{22}\frac{\partial \mu}{\partial x_1} \!+\!\frac{\partial \mu}{\partial x_2} \big) \xi_1\! +\!i\big( g_{21}\frac{\partial \mu}{\partial x_1}\!-\! g_{11} \frac{\partial \mu}{\partial x_2}\big) \xi_2 \end{bmatrix},\nonumber\\
     \nonumber\\
        &  \!c_0(x, \xi')\!= \!\! \begin{bmatrix} \!a_{11}\!+\!\omega^2 \mu\sigma \!-\!R_1^{1} \!\!&\! a_{12}\!-\!R_2^{1} \!\!&\! a_{13}\!-\!R_3^{1}\\
\! a_{21}\!- \!R_1^{2}\! \!&\!a_{22}\!+\!\omega^2 \mu\sigma\!-\!R_2^{2}\! \!& \!a_{23}\!-\!R_3^{2} \\
  \!  a_{31}\!-\!R_1^{3}\! \!&\! a_{32}\!-\!R_2^{3} \!\! &\! a_{33}\!+\! \omega^2 \mu\sigma \!- \!R_3^{3}\end{bmatrix}
  \!+\!\begin{bmatrix} \!\sum\limits_{\beta}\! g^{1\beta} \!\frac{\partial}{\partial x_\beta} \!\big(\!\frac{1}{\sigma} \!\frac{\partial \sigma}{\partial x_1}\!\big) \! \!&  \sum\limits_{\beta}\! g^{1\beta} \! \frac{\partial }{\partial x_\beta}\big(\!\frac{1}{\sigma} \!\frac{\partial \sigma}{\partial x_2}\!\big)
\!\!& \sum\limits_{\beta} \!g^{1\beta} \! \frac{\partial }{\partial x_\beta}\!\big(\!\frac{1}{\sigma} \!\frac{\partial \sigma}{\partial x_3}\!\big)\! \\
\!\sum\limits_{\beta} \!g^{2\beta}\! \frac{\partial }{\partial x_\beta} \!\big(\!\frac{1}{\sigma}\! \frac{\partial \sigma}{\partial x_1}\!\big) \!\! &  \sum\limits_{\beta}\! g^{2\beta} \!\frac{\partial }{\partial x_\beta}\!\big(\!\frac{1}{\sigma} \!\frac{\partial \sigma}{\partial x_2}\!\big)\!
\! \!& \sum\limits_{\beta} \!g^{2\beta}\! \frac{\partial }{\partial x_\beta}\!\big(\!\frac{1}{\sigma}\! \frac{\partial \sigma}{\partial x_3}\!\big)\! \\
 \! \frac{\partial }{\partial x_3} \!\big(\frac{1}{\sigma}\! \frac{\partial \sigma}{\partial x_1}\!\big) \! &  \frac{\partial }{\partial x_3}\!\big(\!\frac{1}{\sigma} \!\frac{\partial \sigma}{\partial x_2}\!\big)
\! &  \frac{\partial }{\partial x_3}\!\big(\!\frac{1}{\sigma}\! \frac{\partial \sigma}{\partial x_3}\!\big)\!\end{bmatrix}\nonumber
\\
\nonumber \\
& +\frac{1}{\mu|g|}\! \begin{bmatrix}- g_{21} \frac{\partial \mu}{\partial x_3} & -g_{22} \frac{\partial \mu}{\partial x_3} & \frac{\partial \mu}{\partial x_2} \\
g_{11} \frac{\partial \mu}{\partial x_3} & g_{12} \frac{\partial \mu}{\partial x_3} & -\frac{\partial \mu}{\partial x_1}\\
 g_{21} \frac{\partial \mu}{\partial x_1}\!-\!g_{11} \frac{\partial \mu}{\partial x_2}  & g_{22} \frac{\partial \mu}{\partial x_1}\!-\! g_{12}
 \frac{\partial \mu}{\partial x_2} & 0\end{bmatrix} \begin{bmatrix}
 -\frac{\partial g_{21}}{\partial x_3} &   -\frac{\partial g_{22}}{\partial x_3} & 0\\
 \frac{\partial g_{11}}{\partial x_3} &   \frac{\partial g_{12}}{\partial x_3} & 0\\
 \frac{\partial g_{21}}{\partial x_1}\!-\! \frac{\partial g_{11}}{\partial x_2} &    \frac{\partial g_{22}}{\partial x_1}\!-\! \frac{\partial g_{12}}{\partial x_2} & 0\end{bmatrix}. \nonumber
              \end{align}
        Note that for any smooth vector-valued function $v$,  \begin{eqnarray*} \left( \phi \big(\frac{\partial}{\partial x_3}I_3 \big) - \big(\frac{\partial}{\partial x_3}I_3\big)\phi \right)v&=&
 \phi \big(\frac{\partial }{\partial x_3}I_3 \big)v- \big(\frac{\partial}{\partial x_3}I_3\big) (\phi v)\\
 &=& \phi \big(\frac{\partial }{\partial x_3} I_3\big)
 v- \bigg(\frac{\partial \phi}{\partial x_3}\bigg) v - \phi \big(\frac{\partial }{\partial x_3} I_3\big) v  = -\bigg(\frac{\partial \phi}{\partial x_3}\bigg) v.\nonumber\end{eqnarray*}
 It follows that   \begin{eqnarray} \label{19.3.19-8}\left( \phi\big(\frac{\partial}{\partial x_3}I_3\big)- \big(\frac{\partial}{\partial x_3} I_3\big)\phi  \right)
= -\frac{\partial \phi}{\partial x_3}.\end{eqnarray}
     Combining this, the right-hand side of (\ref{19.3.19-3}) and symbol formula for product of two pseudodifferential operators (see, for example,  p.$\,37$ of \cite{Tre}, p.$\,$13 of \cite{Ta2}, Theorem 18.1.8 of \cite{Ho3} or \cite{KN})\, we get the full symbol equation for (\ref{19.3.19-3}):
   \begin{eqnarray} \label{19.3.19-4} \sum_{|\vartheta|\ge 0} \frac{(-i)^{|\vartheta|}}{\vartheta!} \big(\partial^{\vartheta}_{\xi'} \phi\big)\big(\partial^\vartheta_{x'}\phi\big)
 -b_0 \phi  -\frac{\partial \phi}{\partial x_3} +c=0,\end{eqnarray}
where $\partial_{x'}^{\vartheta}= \frac{\partial^{|\vartheta|}}{\partial x_1^{\vartheta_1} \partial x_2^{\vartheta_2}}$, $\partial_{\xi'}^{\vartheta}= \frac{\partial^{|\vartheta|}}{\partial \xi_1^{\vartheta_1} \partial\xi_2^{\vartheta_2}}$, and $\vartheta=(\vartheta_1,\vartheta_2)$ is a $2$-tuple of nonnegative integers.

 Group the homogeneous terms of degree two in (\ref{19.3.19-4}) we get
\begin{eqnarray*} \label{19.3.19-16} \phi_1^2 +c_2=0, \end{eqnarray*}
so    \begin{align*}  \phi_1(x,\xi') = \pm  \sqrt{\sum_{\alpha, \beta} g^{\alpha\beta} \xi_\alpha\xi_\beta} \; I_3.
\end{align*}
Because $\nu$ is the outward normal vector of $\partial M$, we may take
\begin{eqnarray} \label{19.3.21-1}   \phi_1(x,\xi') =  \sqrt{\sum_{\alpha, \beta} g^{\alpha\beta} \xi_\alpha\xi_\beta} \; I_3.
\end{eqnarray}
The terms of degree one in (\ref{19.3.19-4}) give
\begin{eqnarray*} \!\!\!\!\!\!\!&& \phi_1 \phi_0 +\phi_0 \phi_1  -i \sum\limits_{m=1}^2 \frac{\partial \phi_1}{\partial \xi_m}\, \frac {\partial \phi_1}{\partial x_m} -b_0\phi_1 -\frac{\partial \phi_1}{\partial x_3} +c_1=0.\end{eqnarray*}
By (\ref{19.3.21-1}) we see that $\phi_1$ has commutative property with any $3\times 3$ matrix, and \begin{eqnarray} \label{19.9.9-1} \phi_1^{-1} (x, \xi')=
\frac{1}{ \sqrt{\sum\limits_{\alpha, \beta} g^{\alpha\beta} \xi_\alpha\xi_\beta}} \; I_3.\end{eqnarray}
From this, we immediately get
\begin{eqnarray} \label{19.8.2-9}  &&\;\,\;\;\quad\,\; \phi_0\! = \! \frac{1}{2\sqrt{\sum\limits_{\alpha,\beta}\! g^{\alpha \beta}\xi_\alpha \xi_\beta}} \left\{ i \sum\limits_{m=1}^2\! \frac{\partial \phi_1}{\partial \xi_m} \frac{\partial \phi_1}{\partial x_m} \!+\!\left(\frac{1}{2}\! \sum\limits_{\alpha, \beta}\! g^{\alpha\beta}  \frac{\partial g_{\alpha\beta}}{\partial x_3}I_3+ \begin{bmatrix} 2\Gamma_{13}^1 & 2\Gamma_{23}^1 &0\\
   2\Gamma_{13}^2 & 2\Gamma_{23}^2  & 0\\
   0& 0& 0\end{bmatrix} \right. \right.\\
   &&  \quad  \quad \;           \left.\left. +\begin{bmatrix} 0& 0& 0\\
   0& 0& 0\\
    \frac{1}{\sigma} \frac{\partial \sigma}{\partial x_1} & \frac{1}{\sigma} \frac{\partial \sigma}{\partial x_2} & \frac{1}{\sigma} \frac{\partial \sigma}{\partial x_3}\end{bmatrix}
    +\begin{bmatrix} -\frac{1}{\mu} \frac{\partial \mu}{\partial x_3} & 0& 0\\
    0 & -\frac{1}{\mu} \frac{\partial \mu}{\partial x_3} & 0 \\
    \frac{1}{\mu}\frac{\partial \mu}{\partial x_1} & \frac{1}{\mu} \frac{\partial \mu}{\partial x_2} & 0\end{bmatrix}
      \right)\phi_1  +\frac{\partial \phi_1}{\partial x_3} \! \right.\nonumber\\
        && \quad\quad \; \left. \, -i \!\sum\limits_{\alpha,\beta} \!\big(g^{\alpha\beta}\!\sum\limits_{\gamma}\! \Gamma_{\alpha\gamma}^\gamma\!+\!\frac{\partial g^{\alpha\beta}}{\partial x_\alpha}\!\big) \xi_\beta\, I_3- \begin{bmatrix} 2i\!\sum\limits_{\alpha,\beta}g^{\alpha\beta} \Gamma_{1\alpha}^1{\xi_\beta}_{{}_{{}_{}}}&  2i\!\sum\limits_{\alpha,\beta}g^{\alpha\beta} \Gamma_{2\alpha}^1{\xi_\beta}_{{}_{{}_{}}}  &  2i\!\sum\limits_{\alpha,\beta}g^{\alpha\beta} \Gamma_{3\alpha}^1{\xi_\beta}_{{}_{}} \\
    2i\!\sum\limits_{\alpha,\beta}g^{\alpha\beta} {\Gamma_{1\alpha}^2{\xi_\beta}}_{{}_{{}_{}}}  &  2i\!\sum\limits_{\alpha,\beta}g^{\alpha\beta} \Gamma_{2\alpha}^2 {\xi_\beta}_{{}_{}} &  2i\!\sum\limits_{\alpha,\beta}g^{\alpha\beta} \Gamma_{3\alpha}^2{\xi_\beta}_{{}_{}} \\
    2i\!\sum\limits_{\alpha,\beta}g^{\alpha\beta} {\Gamma_{1\alpha}^3{\xi_\beta}}_{{}_{{}_{}}}  &  2i\!\sum\limits_{\alpha,\beta}g^{\alpha\beta} \Gamma_{2\alpha}^3 {\xi_\beta}_{{}_{}}  &  0\end{bmatrix}\right. \nonumber\\
    && \left.   \qquad \;    -\begin{bmatrix} i \sum\limits_{\beta}\! g^{1\beta} \!\big(\frac{1}{\sigma} \!\frac{\partial \sigma}{\partial x_1}\big) \xi_\beta \! & i\sum\limits_{\beta}\! g^{1\beta} \!\big(\frac{1}{\sigma} \!\frac{\partial \sigma}{\partial x_2}\big) \xi_\beta
\!& i\sum\limits_{\beta} \!g^{1\beta} \!\big(\frac{1}{\sigma} \!\frac{\partial \sigma}{\partial x_3}\big)  \xi_\beta\\
i\sum\limits_{\beta} \!g^{2\beta}\! \big(\frac{1}{\sigma}\! \frac{\partial \sigma}{\partial x_1}\big) \xi_\beta \! &  i\sum\limits_{\beta} g^{2\beta} \!\big(\frac{1}{\sigma} \!\frac{\partial \sigma}{\partial x_2}\big) \xi_\beta
\! & i\sum\limits_{\beta} \!g^{2\beta}\! \big(\frac{1}{\sigma}\! \frac{\partial \sigma}{\partial x_3}\big) \xi_\beta\\
0 & 0
& 0\end{bmatrix} \right.\nonumber\\
\!\!\!\!\!\!\!&& \left. - \frac{1}{\mu|g|}\! \!\begin{bmatrix}\! i\frac{\partial \mu}{\partial x_2} \big( g_{21} \xi_1 \!-\!g_{11} \xi_2\big) & i\frac{\partial \mu}{\partial x_2} \big( g_{22} \xi_1 -g_{12} \xi_2\big) \! & i\frac{\partial \mu}{\partial x_3} \big( g_{22} \xi_1\! -\!g_{21} \xi_2\big)\\
    \! -i\frac{\partial \mu}{\partial x_1} \big( g_{21} \xi_1 \!-\!g_{11} \xi_2\big) & -i\frac{\partial \mu}{\partial x_1} \!\big( g_{22} \xi_1\! -\!g_{12} \xi_2\big) \!& i\frac{\partial \mu}{\partial x_3} \!\big(\!\!- \!g_{12} \xi_1 \!+\!g_{11} \xi_2\big)\\
 \!    0& 0 \!& \!i\big(\!\!-\!g_{22}\frac{\partial \mu}{\partial x_1} \!+\!\frac{\partial \mu}{\partial x_2} \big) \xi_1\! +\!i\big( g_{21}\frac{\partial \mu}{\partial x_1}\!-\! g_{11} \frac{\partial \mu}{\partial x_2}\big) \xi_2 \end{bmatrix}
    \!\right\}.\nonumber\end{eqnarray}
The terms of degree zero are
\begin{align*} &\phi_1 \phi_{-1} + \phi_{-1}\phi_1  +\phi_0^2 -i\sum\limits_{m=1}^2 \big(\frac{\partial \phi_1}{\partial \xi_m}\, \frac{\partial \phi_0}{\partial x_m} +\frac{\partial \phi_0}{\partial \xi_m}\, \frac{\partial \phi_1}{\partial x_m} \big) -\frac{1}{2}\! \sum\limits_{m,l=1}^2
\frac{\partial^2 \phi_1}{\partial \xi_m \partial \xi_l}\, \frac{\partial^2 \phi_1}{\partial x_m\partial x_l}   - b_0 \phi_0\!-\!\frac{\partial \phi_0}{\partial x_3} \!+\! c_0\!=\!0.\end{align*}
It follows from this, (\ref{19.3.21-1}) and (\ref{19.9.9-1}) that
\begin{align} \label{19.8.2-10} & \phi_{-1} = \frac{1}{2\sqrt{\sum\limits_{\alpha,\beta} g^{\alpha \beta}\xi_\alpha \xi_\beta}} \left\{ -\phi_0^2 +i\sum\limits_{m=1}^2 \big(\frac{\partial \phi_1}{\partial \xi_m}\, \frac{\partial \phi_0}{\partial x_m} +\frac{\partial \phi_0}{\partial \xi_m}\, \frac{\partial \phi_1}{\partial x_m} \big) +\frac{1}{2} \sum\limits_{m,l=1}^2
\frac{\partial^2 \phi_1}{\partial \xi_m \partial \xi_l}\, \frac{\partial^2 \phi_1}{\partial x_m\partial x_l} \right. \\
&      \left.+\!\left(\!\frac{1}{2}\! \sum\limits_{\alpha, \beta}\! g^{\alpha\beta}  \frac{\partial g_{\alpha\beta}}{\partial x_3}I_3+ \begin{bmatrix} 2\Gamma_{13}^1 & 2\Gamma_{23}^1 &0\\
   2\Gamma_{13}^2 & 2\Gamma_{23}^2  & 0\\
   0& 0& 0\end{bmatrix}  +\begin{bmatrix} 0& 0& 0\\
   0& 0& 0\\
    \frac{1}{\sigma} \frac{\partial \sigma}{\partial x_1} & \frac{1}{\sigma} \frac{\partial \sigma}{\partial x_2} & \frac{1}{\sigma} \frac{\partial \sigma}{\partial x_3}\end{bmatrix}
    + \begin{bmatrix} -\frac{1}{\mu} \frac{\partial \mu}{\partial x_3} & 0& 0\\
    0 & -\frac{1}{\mu} \frac{\partial \mu}{\partial x_3} & 0 \\
    \frac{1}{\mu}\frac{\partial \mu}{\partial x_1} & \frac{1}{\mu} \frac{\partial \mu}{\partial x_2} & 0\end{bmatrix}
      \right)\phi_0 \! \right.\nonumber\\
      &      +\!\frac{\partial \phi_0}{\partial x_3} \!-\! \! \begin{bmatrix} \!a_{11}\!+\!\omega^2 \mu\sigma \!-\!\!R_1^{1} \!\!\!&\! a_{12}\!-\!\!R_2^{1} \!\!\!&\!\! a_{13}\!-\!\!R_3^{1}\\
\! a_{21}\!- \!\!R_1^{2}\! \!\!&\!a_{22}\!+\!\omega^2 \mu\sigma\!-\!\!R_2^{2}\! \!\!& \!\!a_{23}\!-\!\!R_3^{2} \\
  \!  a_{31}\!-\!\! R_1^{3}\!\! \!&\!\! a_{32}\!-\!\!R_2^{3} \!\! \!&\! \!a_{33}\!+\! \omega^2 \mu\sigma \!-\! \!R_3^{3}\end{bmatrix}  \!\!-\!\!\begin{bmatrix} \!\sum\limits_{\beta}\! g^{1\beta} \!\frac{\partial}{\partial x_\beta} \!\big(\!\frac{1}{\sigma} \!\frac{\partial \sigma}{\partial x_1}\!\big) \! \!&  \sum\limits_{\beta}\! g^{1\beta} \! \frac{\partial }{\partial x_\beta}\big(\!\frac{1}{\sigma} \!\frac{\partial \sigma}{\partial x_2}\!\big)
\!\!& \sum\limits_{\beta} \!g^{1\beta} \! \frac{\partial }{\partial x_\beta}\!\big(\!\frac{1}{\sigma} \!\frac{\partial \sigma}{\partial x_3}\!\big)\! \\
\!\sum\limits_{\beta} \!g^{2\beta}\! \frac{\partial }{\partial x_\beta} \!\big(\!\frac{1}{\sigma}\! \frac{\partial \sigma}{\partial x_1}\!\big) \!\! &  \sum\limits_{\beta}\! g^{2\beta} \!\frac{\partial }{\partial x_\beta}\!\big(\!\frac{1}{\sigma} \!\frac{\partial \sigma}{\partial x_2}\!\big)\!
\! \!& \sum\limits_{\beta} \!g^{2\beta}\! \frac{\partial }{\partial x_\beta}\!\big(\frac{1}{\sigma}\! \frac{\partial \sigma}{\partial x_3}\!\big)\! \\
 \! \frac{\partial }{\partial x_3} \!\big(\!\frac{1}{\sigma}\! \frac{\partial \sigma}{\partial x_1}\!\big) \! &  \frac{\partial }{\partial x_3}\!\big(\!\frac{1}{\sigma} \!\frac{\partial \sigma}{\partial x_2}\!\big)
\! &  \frac{\partial }{\partial x_3}\!\big(\!\frac{1}{\sigma}\! \frac{\partial \sigma}{\partial x_3}\!\big)\!\end{bmatrix}\nonumber
\\
& \left. -\frac{1}{\mu|g|}\! \begin{bmatrix}- g_{21} \frac{\partial \mu}{\partial x_3} & -g_{22} \frac{\partial \mu}{\partial x_3} & \frac{\partial \mu}{\partial x_2} \\
g_{11} \frac{\partial \mu}{\partial x_3} & g_{12} \frac{\partial \mu}{\partial x_3} & -\frac{\partial \mu}{\partial x_1}\\
 g_{21} \frac{\partial \mu}{\partial x_1}\!-\!g_{11} \frac{\partial \mu}{\partial x_2}  & g_{22} \frac{\partial \mu}{\partial x_1}\!-\! g_{12}
 \frac{\partial \mu}{\partial x_2} & 0\end{bmatrix} \begin{bmatrix}
 -\frac{\partial g_{21}}{\partial x_3} &   -\frac{\partial g_{22}}{\partial x_3} & 0\\
 \frac{\partial g_{11}}{\partial x_3} &   \frac{\partial g_{12}}{\partial x_3} & 0\\
 \frac{\partial g_{21}}{\partial x_1}\!-\! \frac{\partial g_{11}}{\partial x_2} &    \frac{\partial g_{22}}{\partial x_1}\!-\! \frac{\partial g_{12}}{\partial x_2} & 0\end{bmatrix}\right\}.\nonumber
\end{align}
For general $m\ge 1$, by proceeding recursively we get
\begin{eqnarray*} && \phi_1 \phi_{-m-1} + \phi_{-m-1} \phi_1  +  \sum_{\substack{j,k,\vartheta\\ -m=j+k-|\vartheta|\\
-m\le j,k\le 1}} \frac{(-i)^{|\vartheta|}}{\vartheta!} (\partial^\vartheta_{\xi'} \phi_j) (\partial^\vartheta_{x'} \phi_k)  -b_0 \phi_{-m} -\frac{\partial \phi_{-m}}{\partial x_3}=0,\end{eqnarray*}
which implies
\begin{eqnarray} \label{19.8.3-1} && \;\;\;\;\;\phi_{\!-m\!-\!1}\! =\! \frac{1}{2\sqrt{\sum\limits_{\alpha,\beta}\! g^{\alpha \beta}\xi_\alpha \!\xi_\beta}}\! \left\{\! -\!\sum_{\substack{j,k,\vartheta\\ -m=j+k-|\vartheta|\\
-m\le j,k\le 1}}\!\! \!\frac{(\!-i)^{|\vartheta|}}{\vartheta!} \!(\partial^\vartheta_{\xi'} \phi_j) (\partial^\vartheta_{x'} \phi_k)\! \!+\!\!\left(\!\!\frac{1}{2}\! \sum\limits_{\alpha\!, \beta}\! g^{\alpha\!\beta}  \frac{\partial g_{\alpha\beta}}{\partial x_3}\!I_3\!+\! \!\begin{bmatrix} 2\Gamma_{13}^1 & 2\Gamma_{23}^1 &0\\
   2\Gamma_{13}^2 & 2\Gamma_{23}^2  & 0\\
   0& 0& 0\end{bmatrix} \right.\right. \\
  && \left. \left.\quad \;\quad \quad +\begin{bmatrix} 0& 0& 0\\
   0& 0& 0\\
    \frac{1}{\sigma} \frac{\partial \sigma}{\partial x_1} & \frac{1}{\sigma} \frac{\partial \sigma}{\partial x_2} & \frac{1}{\sigma} \frac{\partial \sigma}{\partial x_3}\end{bmatrix}
    +\begin{bmatrix} -\frac{1}{\mu} \frac{\partial \mu}{\partial x_3} & 0& 0\\
    0 & -\frac{1}{\mu} \frac{\partial \mu}{\partial x_3} & 0 \\
    \frac{1}{\mu}\frac{\partial \mu}{\partial x_1} & \frac{1}{\mu} \frac{\partial \mu}{\partial x_2} & 0\end{bmatrix}
      \right)  \phi_{-m}+\frac{\partial \phi_{-m}}{\partial x_3}\!\right\}, \;\;\;  m\!\ge \!1.\nonumber\end{eqnarray}
The proof is completed.  $\qquad \qquad \square$

 \vskip 0.20 true cm

We have obtained the full symbol $\phi(x, \xi')\sim \sum_{l\le 1} \phi_l(x, \xi')$ of the pseudodifferential operator $\Phi$ from above Proposition 2.4. This implies that modulo a smoothing operator, the pseudodifferential operator $\Phi$ has been determined on $\partial M$.

Recall that if $\omega$ is not a resonant frequency,  then for $f\in TH^{\frac{1}{2}} (\partial M)\cap C(\partial M)$ with $\mbox{supp}\,f\subset \Gamma$,
there exists a unique solution  $(E, H)\in ({\mathcal{D'}}(M))^3 \times ({\mathcal{D'}}(M))^3$ of Maxwell's equations
\begin{eqnarray} \label{19.7.20-4} \left\{ \begin{array}{ll}
\mbox{curl}\, E =i\omega \mu H \quad &\mbox{in}\;\, M,\\
\mbox{curl}\, H =-i\omega \sigma E \quad &\mbox{in}\;\, M,\\
\nu\times E=f\quad &\mbox{on}\;\, \partial M.\end{array} \right.\end{eqnarray}

\vskip 0.25 true cm

\noindent{\bf Proposition 2.5.} \ {\it If $E$ solves (\ref{19.7.20-4}), then
\begin{eqnarray}\label{19.7.20-5} \left(\frac{\partial}{\partial x_3}\, I_3\right) E = \Phi E\big|_{\partial M}\end{eqnarray}
modulo a smoothing operator. }

\vskip 0.26 true cm

 \noindent  {\it Proof.} \  Clearly, $\mathcal{M}(x,D) E=0$. It is well-known (see \cite{Tre}, Ch. III, Remark 1.2 $\&$ 4.1) that
  the solution $(E, H)$ of Maxwell's equations (\ref{19.7.20-4}) is smooth in normal variable, i.e., in boundary normal coordinates $(x',x_3)$ with $x_3\in [0,T]$, $E\in (C^\infty([0,T]; {\mathfrak{D}}' ({\mathbb{R}}^{2})))^3$ locally.
       Actually,  $E$ is smooth in the interior of $M$ by interior regularity for elliptic equation system ${\mathcal{M}}_g$ (see, for example, \cite{ADN} or \cite{GT}). From Proposition 2.4, we see that (\ref{19.7.20-4}) is locally equivalent to the following system of equations for $E,U\in (C^\infty([0,T]; {\mathfrak{D}}' ({\mathbb{R}}^{2})))^{3}$: \begin{eqnarray*} &&\bigg(\frac{\partial }{\partial x_3}\, I_3 + \Phi\bigg) E=U, \quad \; \nu\times E\big|_{x_3=0}=f,\\
 &&\bigg(\frac{\partial }{\partial x_3} I_3 +B-\Phi\bigg)U =W\in (C^\infty ([0,T]\times R^{2}))^3. \end{eqnarray*}
 Making the substitution $t=T-x_3$ for the second equation mentioned above (as done in \cite{LU}), we get a backwards generalized heat equation system:
 \begin{eqnarray*} \bigg(\frac{\partial}{\partial t} I_3\bigg) U -(-\Phi+B)U =-W.\end{eqnarray*}
  In view of the principal symbol $\phi_1(x,\xi')$ of $\Phi$ is strictly positive for any $\xi\ne 0$, we get that the solution operator for this heat equation system is smooth for $t>0$ (see p.$\,$134 of \cite{Tre} or \cite{Frie}). This implies that
    $U$ is smooth in the interior of $M$, and hence
    $U\big|_{x_3=T}$ is smooth.
    Therefore,
    \begin{eqnarray*} \left(\frac{\partial }{\partial x_3}I_3\right) E +\Phi E= U\in (C^\infty ([0,T]\times R^{2}))^3 \end{eqnarray*}
    locally.  Setting $Rf= U\big|_{\partial M} $, we immediately see that $R$ is a smoothing operator and
    \begin{eqnarray*} \bigg(\big(\frac{\partial}{\partial x_3} I_3 \big)E\bigg)\bigg|_{\partial M} = -\Phi E\big|_{\partial M} +Rf. \end{eqnarray*}
    The desired result is obtained.
$  \ \  \square$

\vskip 0.35 true cm

Recall that if $(E,H)$ solves (\ref{19.7.20-4})  then the electromagnetic Dirichlet-to-Neumann map $\Lambda_{g,\Gamma}$ is the map
$\Lambda_{g,\Gamma}: \nu\times E\big|_{\Gamma} \to  \nu \times H\big|_{\Gamma}$,  where $\nu = (0, 0, -1)$ is the outward unit normal to the boundary $\Gamma\subset\partial M$ in the boundary normal coordinates. Then we have
\begin{equation*}   \begin{split} &  \nu\times E =\sqrt{|g|}\left\{ \left|\begin{matrix} g^{11} & g^{21} & g^{31}\\
0 & 0& -1\\
 E^1 & E^2 & E^3 \end{matrix}\right|  \frac{\partial }{\partial x_1} + \left|\begin{matrix} g^{12} & g^{22} & g^{32}\\
0 & 0& -1\\
 E^1 & E^2 & E^3 \end{matrix}\right|  \frac{\partial }{\partial x_2}+ \left|\begin{matrix} g^{13} & g^{23} & g^{33}\\
0 & 0& -1\\
 E^1 & E^2 & E^3 \end{matrix}\right|  \frac{\partial }{\partial x_3}\right\}\\
 &      \qquad \;\;=\sqrt{|g|} \left\{ (-g^{21} E^1 +g^{11} E^2 ) \frac{\partial }{\partial x_1} + (-g^{22} E^1 +g^{12} E^2 ) \frac{\partial }{\partial x_2}
+ (-g^{23} E^1 +g^{13} E^2 ) \frac{\partial }{\partial x_3}\right\},\end{split}\end{equation*}
which can be rewritten as the components of vector field with respect to coordinates basis $\frac{\partial}{\partial x_1}$, $\frac{\partial}{\partial x_2}$, $\frac{\partial}{\partial x_3}$:
\begin{eqnarray*}\nu\times E
=\left[\begin{array} {ll} \sqrt{|g|} (-g^{21} E^1 +g^{11} E^2)\\
\sqrt{|g|} (-g^{22} E^1 +g^{12} E^2)\\
\sqrt{|g|} (-g^{23}E^1 +g^{13}E^2)\end{array}\right]= \left[\begin{array} {ll} \sqrt{|g|} (-g^{21} E^1 +g^{11} E^2)\\
\sqrt{|g|} (-g^{22} E^1 +g^{12} E^2)\\
\qquad \qquad \;\;\, 0\end{array}\right]
\end{eqnarray*}
because of $g^{13}=g^{23}=0$ in boundary normal coordinates.
From (\ref{19.7.28-1}) and (\ref{19.7.20-4}), we have
\begin{align*}   &\!\nu \times H= \frac{1}{i\omega \mu} \, \nu\times (\mbox{curl}\, E) \\
 &     \;\quad\;\;\;\;\, = \frac{1}{i\omega \mu} \, \nu \times \frac{1}{\sqrt{|g|}}\left\{ \!\bigg(\frac{\partial}{\partial x_2} (\sum_{l=1}^3 g_{3l}E^l) -  \frac{\partial} {\partial x_3}(\sum_{l=1}^2 g_{2l}E^l)\bigg)\frac{\partial }{\partial x_1} \!+\!\bigg(\frac{\partial}{\partial x_3} (\sum_{l=1}^3 g_{1l}E^l) -  \frac{\partial }{\partial x_1}(\sum_{l=1}^3 g_{3l}E^l)\bigg)\frac{\partial }{\partial x_2}\right.\\
& \left.\qquad \;\;\;\;\,+ \bigg(\frac{\partial }{\partial x_1}(\sum_{l=1}^3 g_{2l}E^l) -  \frac{\partial }{\partial x_2}(\sum_{l=1}^3 g_{1l}E^l)\bigg)\frac{\partial }{\partial x_3}\right\}\\
 &    \; \;\quad\;\;\,\; = \frac{1}{i\omega \mu} \left\{ \!\bigg(\!\!-g^{21}  \sum\limits_{l=1}^3 \big(\frac{\partial(g_{3l}E^l)}{\partial x_2} -\frac{\partial (g_{2l}E^l)}{\partial x_3} \big) + g^{11} \sum_{l=1}^3 \big(\frac{\partial (g_{1l}E^l)}{\partial x_3} -\frac{\partial (g_{3l}E^l)}{\partial x_1}\big)\bigg)\frac{\partial }{\partial x_1}\right.\\
 &    \; \;\quad\;\;\,\;\;\;\;\; +\bigg(\!\!-g^{22}  \sum\limits_{l=1}^3 \big(\frac{\partial(g_{3l}E^l)}{\partial x_2} -\frac{\partial (g_{2l}E^l)}{\partial x_3} \big) + g^{12} \sum_{l=1}^3 \big(\frac{\partial (g_{1l}E^l)}{\partial x_3} -\frac{\partial (g_{3l}E^l)}{\partial x_1}\big)\bigg)\frac{\partial }{\partial x_2}\\
 &  \; \;\quad\;\;\,\;\;\;\;\;\left. +\bigg(\!\!-g^{23}  \sum\limits_{l=1}^3 \big(\frac{\partial(g_{3l}E^l)}{\partial x_2} -\frac{\partial (g_{2l}E^l)}{\partial x_3} \big) + g^{13} \sum_{l=1}^3 \big(\frac{\partial (g_{1l}E^l)}{\partial x_3} -\frac{\partial (g_{3l}E^l)}{\partial x_1}\big)\bigg)\frac{\partial }{\partial x_3}
 \right\}\\
 &     \; \quad \;\;\;\; \;= \frac{1}{i\omega \mu}  \bigg(\! g^{11} \sum\limits_{l=1}^3 \big(\frac{\partial \big( g_{1l} E^l\big)}{\partial x_3} -\frac{\partial (g_{3l} E^l)}{\partial x_1} \big)
+ g^{21} \sum\limits_{l=1}^3 \big(\frac{\partial( g_{2l} E^l)}{\partial x_3} -\frac{\partial (g_{3l} E^l)}{\partial x_2} \big) \bigg) \frac{\partial}{\partial x_1}  \\
  &            \left.\; \; \,\;\quad  \;\;\quad \, + \bigg( \!g^{12} \sum\limits_{l=1}^3\big(\frac{\partial( g_{1l} E^l)}{\partial x_3} -\frac{\partial (g_{3l} E^l)}{\partial x_1} \big) + g^{22} \sum\limits_{l=1}^3\big(\frac{\partial( g_{2l} E^l)}{\partial x_3} -\frac{\partial (g_{3l} E^l)}{\partial x_2} \big) \bigg) \frac{\partial}{\partial x_2}  \right\}.\end{align*}
Thus, by rewriting the components of vector fields with respect to the coordinates basis $\frac{\partial }{\partial x_1}, \frac{\partial }{\partial x_2}, \frac{\partial }{\partial x_3}$,  the map $\,\Lambda: \nu\times E\to \nu\times H$ becomes   \begin{align*} &      \Lambda:  \left[ \begin{array}{ll} \sqrt{|g|} (-g^{21} E^1 +g^{11} E^2 )\\
\sqrt{|g|} (-g^{22} E^1 +g^{12} E^2 )\\
\qquad\qquad \;\; 0\end{array} \right] \mapsto \frac{1}{i\omega \mu}\left[ \begin{array}{ll}  g^{11}\sum\limits_{l=1}^3  \big(\frac{\partial( g_{1l} E^l)}{\partial x_3} -\frac{\partial (g_{3l} E^l)}{\partial x_1} \big)+ g^{21}\sum\limits_{l=1}^3  \big(\frac{\partial( g_{2l} E^l)}{\partial x_3} -\frac{\partial (g_{3l} E^l)}{\partial x_2} \big) \\
   g^{12} \sum\limits_{l=1}^3\big(\frac{\partial( g_{1l} E^l)}{\partial x_3} -\frac{\partial (g_{3l} E^l)}{\partial x_1} \big)+ g^{22} \sum\limits_{l=1}^3\big(\frac{\partial( g_{2l} E^l)}{\partial x_3} -\frac{\partial (g_{3l} E^l)}{\partial x_2} \big)  \\
 \qquad \qquad \qquad \quad\;\; \qquad 0\end{array}\right]. \end{align*}
Note that  \begin{eqnarray} \label{19.8.3;1} \left.\begin{array}{ll} g^{11}g_{11} +g^{12} g_{21}=1, \quad  g^{11} g_{12} +g^{12} g_{22} =0,\\
 g^{21}    g_{11}+g^{22}g_{21} =0, \quad g^{21} g_{12} +g^{22} g_{22}=1.\end{array} \right.\end{eqnarray}
Therefore,  in boundary normal coordinates, we find by (\ref{19.8.3;1}) that
\begin{align} \label{19.7.21-5} &      \nu\times H= \frac{1}{i\omega \mu} \begin{bmatrix}  g^{11} \big(\sum\limits_{\alpha} \frac{\partial (g_{1\alpha} E^\alpha)}{\partial x_3} -\frac{\partial E^3}{\partial x_1} \big)+g^{21}\big( \sum\limits_{\alpha} \frac{\partial (g_{2\alpha} E^\alpha)}{\partial x_3} -\frac{\partial E^3}{\partial x_2} \big) \\
  g^{12} \big(\sum\limits_{\alpha} \frac{\partial (g_{1\alpha} E^\alpha)}{\partial x_3} -\frac{\partial E^3}{\partial x_1} \big)+ g^{22} \big(\sum\limits_{\alpha} \frac{\partial (g_{2\alpha} E^\alpha)}{\partial x_3} -\frac{\partial E^3}{\partial x_2} \big)  \\
   0\end{bmatrix} \\
   &       \quad \quad \;\;\,    = \frac{1}{i\omega \mu} \begin{bmatrix}  \sum\limits_{\alpha} \left( g^{11} g_{1\alpha} \frac{\partial E^\alpha}{\partial x_3} +g^{11} \frac{\partial g_{1\alpha}}{\partial x_3} E^\alpha  + g^{21} g_{2\alpha} \frac{\partial E^\alpha}{\partial x_3} +g^{21} \frac{\partial g_{2\alpha}}{\partial x_3} E^\alpha\right)  - g^{11} \frac{\partial E^3}{\partial x_1} - g^{21} \frac{\partial E^3}{\partial x_2}\\ \\
       \sum\limits_{\alpha} \left( g^{12} g_{1\alpha} \frac{\partial E^\alpha}{\partial x_3} +g^{12} \frac{\partial g_{1\alpha}}{\partial x_3} E^\alpha  + g^{22} g_{2\alpha} \frac{\partial E^\alpha}{\partial x_3} +g^{22} \frac{\partial g_{2\alpha}}{\partial x_3} E^\alpha\right)  - g^{12} \frac{\partial E^3}{\partial x_1} - g^{22} \frac{\partial E^3}{\partial x_2}\\  \\
       0 \end{bmatrix}  \nonumber \\
  &              \qquad \;\,\;  = \frac{1}{i\omega \mu} \left[ \begin{matrix}\sum\limits_{\alpha} (g^{11} g_{1\alpha} +g^{21} g_{2\alpha} ) \frac{\partial E^\alpha}{\partial x_3} +\sum\limits_{\alpha} (g^{11} \frac{\partial g_{1\alpha}}{\partial x_3}
   +g^{21} \frac{\partial g_{2\alpha}}{\partial x_3} ) E^\alpha -\sum\limits_{\alpha} g^{\alpha 1} \frac{\partial E^3}{\partial x_\alpha} \\
   \\
   \sum\limits_{\alpha} (g^{12} g_{1\alpha} +g^{22} g_{2\alpha} ) \frac{\partial E^\alpha}{\partial x_3} +\sum\limits_{\alpha} (g^{12} \frac{\partial g_{1\alpha}}{\partial x_3}
   +g^{22} \frac{\partial g_{2\alpha}}{\partial x_3} ) E^\alpha -\sum\limits_{\alpha} g^{\alpha 2} \frac{\partial E^3}{\partial x_\alpha} \\
   \\
      0 \end{matrix}\right]\nonumber\\
   &      \qquad \;\,\; =\frac{1}{i\omega \mu} \left[ \begin{matrix} \frac{\partial E^1}{\partial x_3}   + \sum\limits_{\alpha,\beta} g^{\beta 1} \,\frac{\partial g_{\beta \alpha }}{\partial x_3} \, E^\alpha   -\sum\limits_{\alpha} g^{\alpha 1} \,\frac{\partial E^3}{\partial x_\alpha}   \\
   \\
    \frac{\partial E^2}{\partial x_3}    + \sum\limits_{\alpha,\beta} g^{\beta 2}\, \frac{\partial g_{\beta \alpha}}{\partial x_3} \, E^\alpha -\sum\limits_{\alpha} g^{\alpha 2}\, \frac{\partial E^3}{\partial x_\alpha} \\ \\
   0 \end{matrix}\right]
   . \nonumber\end{align}
According to (\ref{2020.3.29-2}), we have
\begin{eqnarray*}&& 0= \sigma \,\mbox{div}\, E +\sum\limits_{k=1}^3 \frac{\partial \sigma}{\partial x_k} E^k\\
&& \;\;\,  = \sigma ( \sum\limits_{k=1}^3 \frac{\partial E^k}{\partial x_k} +\sum\limits_{k,l=1}^3 \Gamma_{kl}^l E^k \big) + \sum\limits_{k=1}^3 \frac{\partial \sigma}{\partial x_k} E^k. \end{eqnarray*}
This implies
 \begin{eqnarray*} & \sigma \big( \frac{\partial E^3}{\partial x_3} + \sum\limits_{l=1}^3 \Gamma_{3l}^l E^3 \big) + \frac{\partial \sigma}{\partial x_3} E^3 =- \sigma \big( \sum\limits_{k=1}^2 \frac{\partial E^k}{\partial x_k} +\sum\limits_{k=1}^2 \sum\limits_{l=1}^3 \Gamma_{kl}^l E^k\big) - \sum\limits_{k=1}^2 \frac{\partial \sigma}{\partial x_k} E^k,\end{eqnarray*}
 i.e.,
   \begin{eqnarray} \label{19.7.21-1} & \sigma \big( \frac{\partial E^3}{\partial x_3} + \sum\limits_{l=1}^3 \Gamma_{3l}^l E^3 \big) + \frac{\partial \sigma}{\partial x_3} E^3 =- \sigma \big( \sum\limits_{\beta} \frac{\partial E^\beta}{\partial x_\beta} +\sum\limits_{\beta} \sum\limits_{\gamma} \Gamma_{\beta \gamma}^\gamma E^\beta\big) - \sum\limits_{\beta} \frac{\partial \sigma}{\partial x_\beta} E^\beta.\end{eqnarray}
It follows from Proposition 2.5  that
$$\bigg(\frac{\partial}{\partial x_3} \, I_3\bigg)E= \Phi E, $$
or equivalently,  \begin{eqnarray*} \begin{bmatrix} \frac{\partial}{\partial x_3} & 0 & 0 \\
0 &  \frac{\partial}{\partial x_3}& 0 \\
0 & 0 & \frac{\partial}{\partial x_3}\end{bmatrix} \begin{bmatrix} E^1 \\
E^2 \\ E^3 \end{bmatrix} = \begin{bmatrix} \Phi^{11} & \Phi^{12} & \Phi^{13}\\
\Phi^{21} & \Phi^{22} & \Phi^{23}\\
\Phi^{31} & \Phi^{32} & \Phi^{33}\end{bmatrix} \begin{bmatrix} E^1 \\
E^2 \\ E^3 \end{bmatrix} .\end{eqnarray*}
This can be written as  \begin{eqnarray}\label{19.7.21-2} \frac{\partial E^k}{\partial x_3} = \Phi^{k1} E^1 +\Phi^{k2} E^2 +\Phi^{k3} E^3 \quad \, \mbox{for}\;\; k=1,2,3.\end{eqnarray}
By taking $k=3$ in (\ref{19.7.21-2}) and  inserting the result into (\ref{19.7.21-1}) we get
\begin{eqnarray*}\sigma\big( \Phi^{31} E^1 +\Phi^{32} E^2 +\Phi^{33} E^3 + \sum\limits_{\beta} \Gamma_{3\beta}^\beta E^3\big) +\frac{\partial \sigma}{\partial x_3} E^3  =
-\sigma \bigg( \sum\limits_{\beta} \frac{\partial E^\beta}{\partial x_\beta} + \sum\limits_{\beta} \big( \sum\limits_{\gamma} \Gamma_{\beta}^\gamma\big) E^\beta\bigg)-\sum\limits_{\beta} \frac{\partial \sigma}{\partial x_\beta}E^\beta,\end{eqnarray*}
 i.e., \begin{eqnarray} && \;\;\,\bigg(\! \sigma \big( \!\Phi^{33}\! +\!\sum\limits_{\beta}\! \Gamma_{3\beta}^\beta\big) \!+\!\frac{\partial \sigma}{\partial x_3}\!\bigg) E^3\! \!=\! \!  - \sigma \!\bigg( \!\!\Phi^{31} E^1\! + \!\Phi^{32} E^2 \!+\!\sum_{\beta}\! \frac{\partial E^\beta}{\partial x_\beta} \!+\! \sum\limits_{\beta}\!
 \big(\!\sum\limits_{\gamma} \!\Gamma_{\beta \gamma}^\gamma\! \big) E^\beta \!\!\bigg)\!\! -\!\sum\limits_{\beta}\! \frac{\partial \sigma}{\partial x_\beta} E^\beta.\end{eqnarray}
 Let $Q(x', D_{x'})$ be a pseudodifferential operator of order $-1$ in $x'$ such that \begin{eqnarray}\label{19.8.4-12} Q\bigg( \sigma \big(\Phi^{33} +\sum\limits_{\beta} \Gamma_{3\beta}^\beta\big) + \frac{\partial \sigma}{\partial x_3} \bigg)=I\end{eqnarray}
  modulo a smoothing operator (We will determine $Q$ by calculating the full symbol of $Q$ in section 3).
Then \begin{eqnarray} \label{19.7.21-4} && E^3 =- Q\bigg[  \sigma \bigg( \!\Phi^{31} E^1\! + \!\Phi^{32} E^2 \!+\!\sum_{\beta} \frac{\partial E^\beta}{\partial x_\beta} \!+\! \sum\limits_{\beta}
 \big(\sum\limits_{\gamma} \Gamma_{\beta \gamma}^\gamma \big) E^\beta \bigg)\! +\!\sum\limits_{\beta} \frac{\partial \sigma}{\partial x_\beta} E^\beta\bigg]\\
&& \;\;\qquad  -Q\bigg[ \sigma \bigg( \sum\limits_{\beta}  \Phi^{3\beta} E^\beta + \sum\limits_{\beta}\frac{\partial E^\beta} {\partial x_\beta} +\sum\limits_{\beta} \sum\limits_{\gamma} \Gamma_{\beta\gamma}^\gamma E^\beta \bigg) +\sum\limits_\beta \frac{\partial \sigma}{\partial x_\beta}E^\beta\bigg] \nonumber\\
&&\quad \;\;  = -Q \bigg[ \sum\limits_{\beta} \bigg( \sigma \Big( \Phi^{3\beta} +\frac{\partial }{\partial x_\beta} +\sum\limits_{\gamma} \Gamma_{\beta\gamma}^\gamma \Big) +\frac{\partial \sigma}{\partial x_\beta} \bigg) E^\beta\bigg]. \nonumber\end{eqnarray}
Inserting this into the right-hand side of (\ref{19.7.21-5}), we find by (\ref{19.7.21-2}) and (\ref{19.7.21-4}) that  \begin{align*}  & \nu\times H  =     \frac{1}{i\omega \mu}
      \left[ \begin{matrix}\Phi^{11} E^1 +\Phi^{12} E^2 +\Phi^{13} E^3  + \sum\limits_{\alpha,\beta} g^{\beta 1} \,\frac{\partial g_{\beta \alpha }}{\partial x_3} \, E^\alpha  - \sum\limits_{\alpha} g^{\alpha 1} \frac{\partial E^3 }{\partial x_\alpha}
         \\
         \Phi^{21} E^1 +\Phi^{22} E^2 +\Phi^{23} E^3
         + \sum\limits_{\alpha,\beta} g^{\beta 2}\, \frac{\partial g_{\beta \alpha}}{\partial x_3} \, E^\alpha - \sum\limits_{\alpha} g^{\alpha 2} \frac{\partial E^3}{\partial x_\alpha}
       \\
   0 \end{matrix}\right]\\
  &                = \frac{1}{i\omega \mu}
      \left[ \begin{matrix}\sum\limits_{\beta} \Phi^{1\beta} E^\beta  +\big(\Phi^{13}  - \sum\limits_{\alpha} g^{\alpha 1} \frac{\partial  }{\partial x_\alpha}\big)E^3
        + \sum\limits_{\alpha,\beta} g^{\alpha 1} \,\frac{\partial g_{\alpha\beta  }}{\partial x_3} \, E^\beta   \\
         \sum\limits_{\beta} \Phi^{2\beta} E^\beta +\big(\Phi^{23} - \sum\limits_{\alpha} g^{\alpha 2} \frac{\partial }{\partial x_\alpha}\big)E^3
         + \sum\limits_{\alpha,\beta} g^{\alpha 2}\, \frac{\partial g_{ \alpha \beta}}{\partial x_3} \, E^\beta
       \\
   0 \end{matrix}\right]\\
   &     \!\!     =\!\frac{1}{i\omega \mu}\!
      \left[ \begin{matrix}\!\sum\limits_{\beta} \Phi^{1\beta} E^\beta \!+\! \big(\!-\Phi^{13} \! +\! \sum\limits_{\alpha} g^{\alpha 1} \frac{\partial }{\partial x_\alpha}\big)
        \Big\{ Q \Big[\sigma \Big(\sum\limits_{\beta} \big(  \frac{\partial }{\partial x_\beta} \!+\!\Phi^{3\beta} \!+\!\sum_{\gamma} \Gamma_{\beta \gamma}^\gamma\big) E^\beta  \Big)\! +\! \sum\limits_\beta  \frac{\partial \sigma}{\partial x_\beta} E^\beta \Big]\Big\}\!+\! \sum\limits_{\alpha,\beta} g^{\alpha 1} \frac{\partial g_{ \alpha \beta}}{\partial x_3}  E^\beta   \\
        \! \sum\limits_{\beta} \Phi^{2\beta} E^\beta \!+\! \big(\!-\Phi^{23}  \!+\! \sum\limits_{\alpha} g^{\alpha 2} \frac{\partial }{\partial x_\alpha} \big)
        \Big\{ Q \Big[\sigma \Big(\sum\limits_{\beta} \big(  \frac{\partial }{\partial x_\beta}\! +\!\Phi^{3\beta} \!+\!\sum_{\gamma} \Gamma_{\beta \gamma}^\gamma\big) E^\beta  \Big) \!+\! \sum\limits_\beta  \frac{\partial \sigma}{\partial x_\beta} E^\beta \Big] \Big\}  \!+\! \sum\limits_{\alpha,\beta} g^{\alpha 2} \frac{\partial g_{ \alpha\beta}}{\partial x_3}  E^\beta
       \\
  \! 0 \end{matrix}\right]\end{align*}
\begin{align*} \!\! \!\!\! \! \!\!\!\!\!\!\!\! \!\!&\! \! \! \!=\!\!\frac{1}{i\omega \mu}\!\!\!
      \left[ \begin{matrix}\!\sum\limits_{\beta}\! \Phi^{\!1\!\beta}\! E^\beta \!\!+\! \!\Big\{\!\!\big(\!\!-\!\!\Phi^{\!13} \!\! +\!\! \sum\limits_{\alpha}\! g^{\alpha\! 1}\! \frac{\partial}{\partial x_\alpha}\! \big) \!Q\! \sum\limits_{\beta}\!\Big[\!\sigma \! \big( \! \frac{\partial }{\partial x_\beta}\!\! +\!\!\Phi^{\!3\!\beta} \!\!+\!\!\sum_{\gamma} \!\Gamma_{\!\!\beta \gamma}^\gamma\!\big)  \!  \!+\! \!  \frac{\partial \sigma}{\partial x_\beta}\! \Big]\!\!\Big\}\! E^\beta\! \!+\! \! \sum\limits_{\alpha}\! g^{\!\alpha \!1}\! Q \!\sum\limits_{\beta} \!\Big[\!\sigma\!\big( \! \frac{\partial }{\partial x_\beta} \!\!+\!\!\Phi^{\!3\!\beta} \! \!+\!\!\sum\limits_{\gamma} \!\Gamma_{\!\!\beta \gamma}^\gamma\!\big)\!\!+\!\! \sum\limits_\beta\!\!
       \frac{\partial \sigma}{\partial x_\beta}\!\!\Big]\!\frac{\partial E^\beta}{\partial x_\alpha}\! \!+\!\! \sum\limits_{\alpha\!,\beta}\! g^{\alpha \!1} \!\frac{\partial g_{\alpha \!\beta }}{\partial x_3} \! E^\beta\\
     \!\!\sum\limits_{\beta}\! \Phi^{\!2\!\beta}\! E^\beta \!\!+\! \!\Big\{\!\!\big(\!\!-\!\!\Phi^{\!23} \!\! +\!\! \sum\limits_{\alpha}\! g^{\!\alpha \! 2}\! \frac{\partial }{\partial x_\alpha}\! \big)\! Q\! \sum\limits_{\beta}\!\Big[\!\sigma \! \big( \! \frac{\partial }{\partial x_\beta}\! \! +\!\!\Phi^{\!3\!\beta}\! \!+\!\!\sum_{\gamma} \!\Gamma_{\!\!\beta \gamma}^\gamma\big)   \! \!+\!\!   \frac{\partial \sigma}{\partial x_\beta}\! \Big]\!\Big\}\! E^\beta\! \!+\! \! \sum\limits_{\alpha}\! g^{\!\alpha\! 2}\! Q \!\sum\limits_{\beta} \!\Big[\!\sigma\!\big( \! \frac{\partial }{\partial x_\beta} \!\!+\!\!\Phi^{\!3\!\beta} \! \!+\!\!\sum\limits_{\gamma} \!\Gamma_{\!\!\beta \gamma}^\gamma\!\big)\!\!+\!\! \sum\limits_\beta\!\!
       \frac{\partial \sigma}{\partial x_\beta}\!\!\Big]\!\frac{\partial E^\beta}{\partial x_\alpha}\! \!+\!\! \sum\limits_{\alpha\!,\beta}\! g^{\!\alpha \! 2}\! \frac{\partial g_{\alpha \!\beta }}{\partial x_3} \! E^\beta
       \\
   0 \end{matrix}\!\right]\\
  \!\! \!  \!\!\! \!\!  \!\!  & \;\; =\frac{1}{i\omega \mu}  \left[ \begin{matrix} L^{11}  &  L^{12}\\
   L^{21}  & L^{22}  \\
    0  &0   \end{matrix} \right]\begin{bmatrix} E^1\\E^2 \end{bmatrix}, \end{align*}
   where \begin{align} \label{19.7.21-7} &  L^{jk} =\Phi^{jk}-\Phi^{j3} Q\! \Big[\sigma  \big( \frac{\partial }{\partial x_k} +\Phi^{3k} +\sum_{\gamma} \!\Gamma_{k \gamma}^\gamma\big)+\frac{\partial \sigma}{\partial x_k}\Big]
      +  \sum\limits_{\alpha}\bigg\{  g^{\alpha j}\frac{\partial }{\partial x_\alpha} \bigg( Q \Big[\sigma\big( \frac{\partial }{\partial x_k} +\Phi^{3k} +\sum\limits_{\gamma} \Gamma_{k \gamma}^\gamma\big)+ \frac{\partial \sigma}{\partial x_k}
      \Big]\bigg)
          \\
      &  +   g^{\alpha j}   Q\Big[\sigma  \big( \frac{\partial }{\partial x_k} +\Phi^{3k} +\sum_{\gamma} \Gamma_{k \gamma}^\gamma\!\big)+\frac{\partial \sigma}{\partial x_k}\Big] \frac{\partial}{\partial x_\alpha}
            +  g^{\alpha j} \!\frac{\partial g_{\alpha k }}{\partial x_3}
   \bigg\},   \qquad \,\;
  1\le j, k \le 2. \nonumber\end{align}
 We have used the fact that if $P$ is a pseudodifferential operator with full symbol $p(x,\xi)$, then for any function $u$,
 \begin{eqnarray} \label{19.9.22-1} \frac{\partial (Pu)}{\partial x_\alpha} = \bigg(\frac{\partial P}{\partial x_\alpha} +P \frac{\partial }{\partial x_\alpha}\bigg)u,\end{eqnarray}
where the operators $\frac{\partial P}{\partial x_\alpha}$ and $P\frac{\partial }{\partial x_\alpha}$ have the full symbols $\frac{\partial p(x,\xi)}{\partial x_\alpha}$ and $p(x,\xi) (i\xi_\alpha)$, respectively.
Indeed, according to the definition of the symbol in a Riemannian manifold (for example, \S 10 of Chapter 7 in \cite{Ta2}) in every local chart we have \begin{eqnarray*} Pu(x)= \int_{{\mathbb{R}}^n} e^{i\langle x,\xi\rangle} p(x,\xi) \hat{u}(\xi) \,d\xi,\end{eqnarray*}
 so that \begin{align*}  \frac{\partial (Pu(x))}{\partial x_\alpha} &=\frac{\partial }{\partial x_\alpha} \int_{{\mathbb{R}}^n}
e^{i\langle x,\xi\rangle} p(x,\xi)\hat{u}(\xi)\, d\xi \\
&= \int_{{\mathbb{ R}}^n}
e^{i\langle x,\xi\rangle} \bigg( \frac{\partial p(x,\xi)}{\partial x_\alpha} +p(x,\xi) (i\xi_\alpha)\bigg)\hat{u}(\xi)\, d\xi, \end{align*}
 which implies that the operator $u\mapsto \frac{\partial (Pu(x))}{\partial x_\alpha}$ has the  full symbol $\frac{\partial p(x,\xi)}{\partial x_\alpha}+ p(x, \xi) (i\xi_\alpha)$. Or equivalently, $\frac{\partial (Pu)}{\partial x_\alpha}$ can be written as (\ref{19.9.22-1}).

 Let $\wp: (C^\infty (M))^3 \to (C^\infty (M))^2$ be a linear isometric operator defined by  \begin{eqnarray*} \wp\left(\begin{bmatrix} A_1\\
 A_2\\ 0\end{bmatrix}\right)  \mapsto    \begin{bmatrix} A_1\\
 A_2\end{bmatrix}.\end{eqnarray*}
 Then  \begin{align} \label{19.8.4-2}\wp(\nu \times E)&=\wp \begin{bmatrix}\sqrt{|g|}(-g^{21} E^1 +g^{11} E^2)\\
 \sqrt{|g|}(-g^{22} E^1 +g^{12} E^2)\\ 0 \end{bmatrix}\\
  &  =\sqrt{|g|}\begin{bmatrix}-g^{21} E^1 +g^{11} E^2\\
   -g^{22} E^1 +g^{12} E^2 \end{bmatrix}= \sqrt{|g|} \begin{bmatrix} -g^{21}  & g^{11} \\
 - g^{22} & g^{12} \end{bmatrix} \begin{bmatrix} E^1\\ E^2\end{bmatrix} \nonumber   \end{align}
 and
  \begin{align} \label{19.8.4-3} \wp(\nu \times H)&=\frac{1}{i\omega \mu}  \left[ \begin{matrix} L^{11}  &  L^{12}\\
   L^{21}  & L^{22}   \end{matrix} \right]\begin{bmatrix} E^1\\E^2 \end{bmatrix}. \end{align}
From (\ref{19.8.4-2}) we have
\begin{align*} \begin{bmatrix} E^1\\ E^2 \end{bmatrix} = \frac{1}{\sqrt{|g|}} {\begin{bmatrix} -g^{21}  & g^{11} \\
 - g^{22} & g^{12} \end{bmatrix}}^{-1}\left( \wp(\nu \times E)\right)= \frac{1}{\sqrt{|g|}}  \begin{bmatrix} -g_{12} & - g_{22}\\
 g_{11} & g_{12}  \end{bmatrix} \wp(\nu\times E),\end{align*}
 so that
 \begin{eqnarray} \label{19.8.28-1} \wp(\nu \times H)&=\frac{1}{i\omega \mu\sqrt{|g|}}  \left[ \begin{matrix} L^{11}  &  L^{12}\\
   L^{21}  & L^{22}   \end{matrix} \right] \begin{bmatrix} -g_{12} & - g_{22}\\
 g_{11} & g_{12}  \end{bmatrix} \wp(\nu\times E).\end{eqnarray}

   Thus we have obtained the following proposition:

\vskip 0.28 true cm

\noindent{\bf Proposition 2.6.} \ {\it  In boundary normal coordinates, the electromagnetic Dirichlet-to-Neumann map $\Lambda_{g,\Gamma}$ is equivalent to the following operator which maps $\wp(\nu\times E)$ to $\wp(\nu \times H)$ defined by
\begin{align}  \label{19.8.4-1}       \wp(\nu \times H) = \left[ \begin{matrix} \Lambda^{11}  &  \Lambda^{12}\\
   \Lambda^{21}  & \Lambda^{22}  \end{matrix} \right]   \wp(\nu\times E) \quad \mbox{on}\;\; \Gamma,
\end{align}
where \begin{eqnarray} \label{19.8.6-1} \left. \begin{array}{ll} & \Lambda^{11}= \frac{1}{i\omega \mu\sqrt{|g|}} \big(-g_{12} L^{11} +g_{11} L^{12}\big), \;\quad \Lambda^{12} = \frac{1}{i\omega \mu\sqrt{|g|}} \big(-g_{22} L^{11} +g_{12} L^{12}\big), \\
& \Lambda^{21}= \frac{1}{i\omega \mu\sqrt{|g|}} \big(-g_{12} L^{21} +g_{11} L^{22}\big), \;\quad \Lambda^{22} = \frac{1}{i\omega \mu\sqrt{|g|}} \big(-g_{22} L^{21} +g_{12} L^{22}\big),\end{array} \right.\end{eqnarray} and $L^{11}, L^{12}, L^{21}, L^{22}$ are given by (\ref{19.7.21-7}).
}

\vskip 0.18 true cm

We will still denote the above equivalent operator by $\Lambda_{g,\Gamma}$.

\vskip 1.48 true cm

\section{Determining metric of manifold from the electromagnetic Dirichlet-to-Neumann map}

\vskip 0.48 true cm

We first calculate the full symbol of $Q$, which was introduced in section 2. Recall that (see (\ref{19.8.4-12}) and (\ref{19.9.1-1}))  \begin{eqnarray*} Q\Big\{ \sigma \big(\Phi^{33} +\frac{1}{2}\sum\limits_{\alpha,\beta} g^{\alpha \beta}\frac{\partial g_{\alpha\beta}}{\partial x_3}\big)+\frac{\partial \sigma}{\partial x_3} \Big\} =I.\end{eqnarray*} We wish to define $Q$ so that
\begin{eqnarray} \label{19.8.4-16} \iota\bigg(Q\Big(\sigma\big(\Phi^{33} +\frac{1}{2}\sum\limits_{\alpha,\beta} g^{\alpha \beta}\frac{\partial g_{\alpha\beta}}{\partial x_3}\big) +\frac{\partial \sigma}{\partial x_3}\Big) \bigg)  \sim 1,\end{eqnarray}
where $\iota(P)$ denotes the full symbol of pseudodifferential operator $P$.
Let $q(x,\xi')\sim \sum\limits_{l\le -1} q_{l}(x, \xi')$ and $\phi\sim \sum\limits_{l\le 1} \phi_{l}(x, \xi')$ be the full symbols of pseudodifferential operators $Q$ and $\Phi$, respectively.
Here, the definition of $\Phi$ is as in Proposition 2.4, and
 \begin{eqnarray*}\phi_l= \begin{bmatrix} \phi_l^{11} & \phi_l^{12}&\phi_l^{13}\\
\phi_l^{21}&\phi_l^{22} &\phi_l^{23}\\
\phi_l^{31}&\phi_l^{32} &\phi_l^{33}\end{bmatrix} \,\, \mbox{is the symbol matrix of $\Phi$ with homogenous of degree}\,\, l\,\, \mbox{in} \;\xi',\;\;l=1, 0, -1, -2,\cdots.\end{eqnarray*}
Then (\ref{19.8.4-16}) leads to the following full symbol equation:
\begin{eqnarray}\! \!\!\!\label{19.8.4-13} & \;\;\;\quad\quad  \sum\limits_{|\vartheta|\ge 0} \!\!\frac{(\!-i)^{|\vartheta|}}{\vartheta!}\!\big(\! \partial^\vartheta_{\xi'} (\!q_{-1} \!+\!
q_{-2}\! +\!q_{-3}\!+\!\cdots\!)\!\big)\!\Big\{\!\partial_{x'}^\vartheta \Big( \!\sigma (\phi^{33}_{1}\! +\! \phi^{33}_0\!+\!\frac{1}{2}\!\sum\limits_{\alpha,\beta} \!g^{\alpha\beta}\frac{\partial g_{\alpha\beta}}{\partial x_3}\!+\! \phi^{33}_{-1}\!+\! \phi^{33}_{-2} \!+\!\cdots\!) \!+\!\frac{\partial \sigma}{\partial x_3}\!\Big)\!\Big\}\!=\!1.\end{eqnarray}
Grouping the homogeneous terms of degree zero in (\ref{19.8.4-13}) we get
\begin{eqnarray*} q_{-1} \big(\sigma\phi^{33}_1) =1,\end{eqnarray*}
so that \begin{eqnarray} \label{19.8.4_14} q_{-1}= \frac{1}{\sigma \phi^{33}_1}=\frac{1}{\sigma\sqrt{\sum\limits_{\alpha,\beta} g^{\alpha\beta} \xi_\alpha\xi_\beta}} .\end{eqnarray}
The terms of degree $-1$ in  (\ref{19.8.4-13}) are
\begin{eqnarray*} q_{-1} \bigg( \sigma\Big( \phi^{33}_{0} +\frac{1}{2} \sum\limits_{\alpha,\beta}g^{\alpha \beta}\frac{\partial g_{\alpha\beta}}{\partial x_3}\Big) +\frac{\partial \sigma}{\partial x_3} \bigg)  +
q_{-2} (\sigma \phi^{33}_1) -i \sum\limits_{m=1}^2 \frac{\partial q_{-1}}{\partial \xi_m} \frac{\partial (\sigma \phi^{33}_{1})}{\partial x_m} =0,\end{eqnarray*}
which implies
\begin{eqnarray} \label{19.8.4-15} q_{-2} = -\frac{1}{\sigma \phi^{33}_{1}}\bigg\{ q_{-1} \Big( \sigma \big( \phi^{33}_{0}+\frac{1}{2}\sum\limits_{\alpha,\beta}
g^{\alpha\beta} \frac{g_{\alpha\beta}}{\partial x_3}\big) +\frac{\partial \sigma}{\partial x_3}\Big)   -i \sum\limits_{m=1}^2 \frac{\partial q_{-1} }{\partial \xi_m}\, \frac{\partial (\sigma \phi^{33}_{1})}{\partial x_m}\bigg\}.\end{eqnarray}
The terms of degree $-2$ in (\ref{19.8.4-13}) are
\begin{align*}  &\;\;\; q_{-1} (\sigma \phi_{-1}^{33}) \! +\!q_{-2} \Big(\!\sigma \big(\! \phi_{0}^{33} \!+\!\frac{1}{2}\! \sum\limits_{\alpha,\beta}\! g^{\alpha \beta} \frac{\partial g_{\alpha\beta}}{\partial x_3}\big)\!+\! \frac{\partial \sigma}{\partial x_3}\! \Big) \!+\!q_{-3} (\sigma\phi^{33}_1)  \!-i \sum\limits_{m=1}^2 \!\frac{\partial q_{-1}}{\partial \xi_m}\frac{\partial}{\partial x_m}
\!\Big(\sigma\big(\phi_0^{33} \!+\!\frac{1}{2} \sum\limits_{\alpha,\beta}g^{\alpha \beta}\frac{\partial g_{\alpha\beta}}{\partial x_3}\big)\!+\!\frac{\partial \sigma}{\partial x_3}\!\Big)\\
&\qquad\;\quad\;\, \;\;\; -i \sum\limits_{m=1}^2 \frac{\partial q_{-2}}{\partial \xi_m}\, \frac{\partial (\sigma\phi_1^{33})}{\partial x_m}-\frac{1}{2} \sum\limits_{m,k=1}^2 \frac{\partial^2q_{-1}}{\partial \xi_m\partial \xi_k} \,\frac{\partial^2 (\sigma\phi^{33}_1)}{\partial x_m \partial x_k}=0,\end{align*}
so \begin{align} \label{19.8.4-16,} &\;\, q_{-3} \!=\! -\! \frac{1}{\sigma\phi^{33}_1} \bigg\{\! q_{-1} (\sigma\phi_{-1}^{33}) +q_{-2}\Big(\sigma \big( \phi_{0}^{33} +\frac{1}{2} \sum\limits_{\alpha,\beta} g^{\alpha \beta} \frac{\partial g_{\alpha\beta}}{\partial x_3}\big) \!+\!\frac{\partial \sigma}{\partial x_3}\Big) \! -i \sum\limits_{m=1}^2 \frac{\partial q_{-1}}{\partial \xi_m} \,\frac{\partial}{\partial x_m}
\!\Big(\sigma\big(\phi_0^{33} \\
&\qquad\;\;\, +\frac{1}{2} \sum\limits_{\alpha,\beta}g^{\alpha \beta}\frac{\partial g_{\alpha\beta}}{\partial x_3}\big)+\frac{\partial \sigma}{\partial x_3}\Big) -i \sum\limits_{m=1}^2 \frac{\partial q_{-2}}{\partial \xi_m}\, \frac{\partial (\sigma\phi_1^{33})}{\partial x_m}-\frac{1}{2} \sum\limits_{m,k=1}^2 \frac{\partial^2q_{-1}}{\partial \xi_m\partial \xi_k} \,\frac{\partial^2 (\sigma\phi^{33}_1)}{\partial x_m \partial x_k} \bigg\}.\nonumber\end{align}
Proceeding recursively, the terms of degree $-m$, ($m>1$), are
\begin{eqnarray*} q_{-m-1} (\sigma\phi_1^{33})+ \sum_{\substack{k,l,\vartheta\\
-m=-k-|\vartheta|+l\\
-m\le -k\le -1, \,\; l\le 1}} \frac{(-i)^{|\vartheta|}}{\vartheta!} \big(\partial_{\xi'}^\vartheta q_{-k}\big)\big(\partial^\vartheta_{x'} {\tilde{\phi}}^{33}_l\big)=0,\end{eqnarray*}
where \begin{eqnarray*}  {\tilde{\phi}}^{33}_l=\left\{ \begin{array}{ll} \sigma\phi^{33}_l\quad &\mbox{if}\;\; l\ne 0,\\
 \sigma(\phi^{33}_0+\frac{1}{2} \sum\limits_{\alpha,\beta}g^{\alpha\beta}\frac{\partial g_{\alpha\beta}}{\partial x_3})+\frac{\partial \sigma}{\partial x_3} \quad &\mbox{if}\;\; l=0,\end{array} \right.\end{eqnarray*}
and hence \begin{eqnarray} \label{19.8.12-5}   q_{-m-1}  =-\frac{1}{\sigma \phi_1^{33}}\left( \sum_{\substack{k,l,\vartheta\\
-m=-k-|\vartheta|+l\\
 -m\le -k\le -1, \,\; l\le 1}} \frac{(-i)^{|\vartheta|}}{\vartheta!} \big(\partial_{\xi'}^\vartheta q_{-k}\big)\big(\partial^\vartheta_{x'} {\tilde{\phi}}^{33}_l\big)\right), \quad \,\; m>1. \end{eqnarray}

 The following calculation will be needed late. It follows from the symbol formula of product of two pseudodifferential operators (see p.$\,$37 of \cite{Tre}  or p.$\,$13 of \cite{Ta2}) that the full symbol of operator $Q\Big(\sigma(\frac{\partial}{\partial x_s} + \Phi^{3s}  +\frac{1}{2}\sum\limits_{\alpha,\beta} g^{\alpha\beta}\frac{\partial g_{\alpha\beta}}{\partial x_s})+\frac{\partial \sigma}{\partial x_3}\Big)$ for each $s=1,2$ is
\begin{eqnarray*} \sum\limits_{|\vartheta|\ge 0} \!\frac{(-i)^{|\vartheta|} }{\vartheta!} \! \left(\partial_{\xi'}^\vartheta (q_{-1} +q_{-2} \!+\!q_{-3}\! +\!\cdots) \right)\!\bigg(\partial_{x'}^\vartheta \Big(\sigma(i\xi_s+\phi_{1}^{3s} +\phi^{3s}_{0}+\frac{1}{2} \sum\limits_{\alpha,\beta}g^{\alpha \beta} \frac{\partial g_{\alpha\beta}}{\partial x_s} +\phi^{3s}_{-1}\! +\! \phi_{-2}^{3s} +\cdots)\!+\!\frac{\partial \sigma}{\partial x_3}\Big)\!\bigg). \end{eqnarray*}
  The corresponding terms of  degree zero, degree $-1$ and degree $-2$  are
  \begin{eqnarray} \label{19.8.8-1} \quad \quad\;\left. \begin{array}{ll} q_{-1} \Big(\sigma(i\xi_s+ \phi_1^{3s})\Big),\\
    q_{-1} \Big(\sigma\big( \phi_0^{3s} +\frac{1}{2} \!\sum\limits_{\alpha,\beta}\!g^{\alpha \beta} \frac{\partial g_{\alpha\beta}}{\partial x_s}\big) +\frac{\partial \sigma}{\partial x_3}\Big) +q_{-2} \Big(\sigma\big(i\xi_s +\phi^{3s}_1\big)\Big) -i\sum\limits_{m=1}^2 \frac{\partial q_{-1}}{\partial \xi_m} \, \frac{\partial (\sigma ( i\xi_s +\phi^{3s}_1))}{\partial x_m},\\
    q_{-1} (\sigma\phi_{-1}^{3s}) +q_{-2} \Big(\sigma\big(\phi_0^{3s}   \! +\!\frac{1}{2} \!\sum\limits_{\alpha,\beta}g^{\alpha \beta} \frac{\partial g_{\alpha\beta}}{\partial x_s}\big)+ \frac{\partial \sigma}{\partial x_3}\Big) + q_{-3} \big(\sigma(i\xi_s +\phi_1^{3s})\big)\! -\!i\!\sum\limits_{m=1}^2\! \frac{\partial q_{-1} }{\partial \xi_m} \frac{\partial }{\partial x_m}\!\big( \sigma(\phi_{0}^{3s}\! \\
    \qquad\quad \quad \,\;\;+ \frac{1}{2} \!\sum\limits_{\alpha,\beta}g^{\alpha \beta} \frac{\partial g_{\alpha\beta}}{\partial x_s}) +\frac{\partial \sigma}{\partial x_3}\big) -i \sum\limits_{m=1}^2  \frac{\partial q_{-2}}{\partial \xi_m} \, \frac{\partial \big(\sigma(i\xi_s+\phi^{3s}_1)\big)}{\partial x_m} - \frac{1}{2}
    \sum\limits_{m,l=1}^2 \frac{\partial^2 q_{-1}}{\partial \xi_m\partial \xi_l} \, \frac{\partial^2 \big(\sigma(i\xi_s+\phi_1^{3s})\big)}{\partial x_m\partial x_l} ,\end{array}\right.\end{eqnarray} respectively.
    Generally, the terms of degree $-m$ of the symbol of $Q\big(\sigma(\frac{\partial}{\partial x_s} + \Phi^{3s} +\frac{1}{2}\sum\limits_{\alpha,\beta} g^{\alpha\beta}\frac{\partial g_{\alpha\beta}}{\partial x_s})+\frac{\partial \sigma}{\partial x_3}\big)$ for each $s=1,2$ is
    \begin{eqnarray}\label{19.8.8-2} \sum_{\substack{k,l,\vartheta\\-m=-k -|\vartheta|+l\\-k\le -1, \;\, l\le 1}} \frac{(-i)^{|\vartheta|}}{\vartheta!} \big(\partial^\vartheta_{\xi'} q_{-k}\big)
    \big( \partial_{x'}^{\vartheta} {\tilde \phi}^{3s}_l \big), \end{eqnarray}
    where \begin{eqnarray*} {\tilde \phi}^{3s}_l =\left\{\begin{array} {ll} \sigma\phi^{3s}_l \quad &\mbox{if}\;\; l\le -1,\\
      \sigma( \phi^{3s}_0 +\frac{1}{2}\sum\limits_{\alpha,\beta} g^{\alpha\beta}\frac{\partial g_{\alpha\beta}}{\partial x_3}) +\frac{\partial \sigma}{\partial x_3} \quad &\mbox{if}\;\; l=0,
       \\
       \sigma( i\xi_s+\phi^{3s}_1)\quad &\mbox{if}\;\; l=1.\end{array} \right.\end{eqnarray*}

\vskip 0.39 true cm

 \noindent{\bf Proposition 3.1.} \ {\it Suppose $\mbox{dim}\, M=3$, and assume the real parts of electromagnetic parameters $\mu$ and $\sigma$ are positive functions in $\bar M$. Let $(x_1, x_{2})$ be any local coordinates for an open set $W\subset \Gamma\subset \partial M$,  and let $\psi\sim \sum_{j\le 1} \psi_j$ denote the full symbol of $\Lambda_{g,\Gamma}$ in these coordinates. Then, $g_{jk}$ and their partial derivatives up to order $m$ are determined by $\psi_1, \psi_0, \cdots, \psi_{-m+1}$ on $\Gamma$ for any $m\ge 0$. Furthermore, for any $x_0\in W$, the full Taylor series of $g$ at $x_0$ in boundary normal coordinates is given by explicit formula in terms of the matrix-valued functions $\{ \psi_j\}_{j\le 1}$
 and their tangential derivatives at $x_0$.}

\vskip 0.22 true cm

\noindent  {\it Proof.} \   Denote by $(x_1,x_2, x_3)$ the boundary normal coordinates associated with $x'=(x_1,x_2)$ as in section 2. According to the form of metric (\ref{18/a-1}) which we have chosen,  we immediately see that it suffices to show that the matrix-valued functions $\{\psi_j\}$ determine the metric $[g_{\alpha\beta}]_{2\times 2}$ and all its normal derivatives along $\Gamma\subset\partial M$. Noticing that $\frac{\partial g_{\alpha\beta}}{\partial x_3} =-\sum\limits_{\rho, \gamma} g_{\alpha \rho} \frac{\partial g^{\rho \gamma}}{\partial x_3} g_{\gamma \beta},$ it is also enough to determine the inverse matrix $[g^{\alpha\beta}]_{2\times 2}$ and all its normal derivatives.

First, according to (\ref{19.8.6-1}) and  (\ref{19.7.21-7}) we have
\begin{align*} \Lambda^{11}=& \frac{1}{i\omega \mu\sqrt{|g|}} \big(-g_{12} L^{11} +g_{11} L^{12}\big)\\
 =& \frac{1}{i\omega \mu \sqrt{|g|}} \bigg\{\!\! -\!g_{12} \bigg[\!\Phi^{11} -\Phi^{13} Q\Big(\sigma\big(\frac{\partial}{\partial x_1} \!+\!\Phi^{31}
 \!+\!\frac{1}{2}\!\sum\limits_{\alpha,\beta} g^{\alpha\beta}\frac{\partial g_{\alpha\beta}}{\partial x_1}\big)+\frac{\partial \sigma}{\partial x_1}\Big)
 +\!\sum\limits_{\alpha}\! g^{\alpha 1} \frac{\partial}{\partial x_\alpha}\!\bigg(Q \Big(\sigma \big(\frac{\partial}{\partial x_1} \!+\!\Phi^{31}
\\
&  +\!\frac{1}{2}\!\sum\limits_{\beta,\gamma}\! g^{\beta\gamma}\frac{\partial g_{\beta\gamma}}{\partial x_1})+\frac{\partial \sigma}{\partial x_1}\Big)\bigg)
 + \sum\limits_{\alpha}\! g^{\alpha 1} Q \Big(\sigma\big(\frac{\partial}{\partial x_1}\! +\!\Phi^{31}
 +\!\frac{1}{2}\!\sum\limits_{\beta,\gamma}\! g^{\beta\gamma}\frac{\partial g_{\beta\gamma}}{\partial x_1}\big)+\frac{\partial \sigma}{\partial x_1}\Big)\!\frac{\partial }{\partial x_\alpha}
 +\sum\limits_{\alpha}
 g^{\alpha 1} \frac{\partial g_{\alpha 1}}{\partial x_3}\bigg]\\
 &  +g_{11} \bigg[\!\Phi^{12} -\Phi^{13} Q\Big(\sigma\big(\frac{\partial}{\partial x_2} \!+\!\Phi^{32}
 \!+\!\frac{1}{2}\!\sum\limits_{\alpha,\beta} g^{\alpha\beta}\frac{\partial g_{\alpha\beta}}{\partial x_2}\big)+\frac{\partial \sigma}{\partial x_2}\Big)
 +\!\sum\limits_{\alpha}\! g^{\alpha 1} \frac{\partial}{\partial x_\alpha}\!\bigg(Q \Big(\sigma \big(\frac{\partial}{\partial x_2} \!+\!\Phi^{32}
\\
&  +\!\frac{1}{2}\!\sum\limits_{\beta,\gamma}\! g^{\beta\gamma}\frac{\partial g_{\beta\gamma}}{\partial x_2})+\frac{\partial \sigma}{\partial x_2}\Big)\bigg)
 + \sum\limits_{\alpha}\! g^{\alpha 1} Q \Big(\sigma\big(\frac{\partial}{\partial x_2}\! +\!\Phi^{32}
 +\!\frac{1}{2}\!\sum\limits_{\beta,\gamma}\! g^{\beta\gamma}\frac{\partial g_{\beta\gamma}}{\partial x_2}\big)+\frac{\partial \sigma}{\partial x_2}\Big)\!\frac{\partial }{\partial x_\alpha}
 +\sum\limits_{\alpha}
 g^{\alpha 1} \frac{\partial g_{\alpha 2}}{\partial x_3}\bigg]\bigg\}.\end{align*}
Let $\delta_{jk}$ be the Kronecker symbol defined by  \begin{eqnarray*}\delta_{jk} =\left\{ \begin{array}{ll} 1 \quad &\mbox{for}\,\; j=k,\\
0 \quad & \mbox{for}\;\, j\ne k.\end{array}\right. \end{eqnarray*}
Since the principal symbols of $Q$ and $\Phi^{jk}$, ($j,k=1,2,3$), are $$q_{-1}=\frac{1}{\sigma}\; \sqrt{\sum_{\alpha,\beta} g^{\alpha\beta}\xi_\alpha\xi_\beta}$$ and $$\phi_1^{jk}=\sqrt{\sum_{\alpha,\beta} g^{\alpha\beta}\xi_\alpha\xi_\beta}\,\delta_{jk},$$ respectively, we see that the symbol of the operator $\Lambda^{11}$  with homogeneous of degree one is
\begin{align} \label{2019.12.4-1}\!\!\!& \! \!  \psi_1^{11}= \frac{1}{i\omega \mu \sqrt{|g|}} \left\{ -g_{12} \bigg( \phi_1^{11}+  \sum\limits_{\alpha}  g^{\alpha 1} q_{-1} \big( \sigma\,i\xi_1\big) (i\xi_\alpha)  \bigg) + g_{11}\bigg( \sum\limits_\alpha g^{\alpha 1} q_{-1} \big(\sigma\, i\xi_2\big)(i\xi_\alpha)\bigg)\right\} \\
& \; \;\;\; \,= \frac{1}{i\omega \mu \sqrt{|g|}}\!\bigg(\!\!-g_{12}\sqrt{\sum\limits_{\alpha,\beta} g^{\alpha\beta} \xi_\alpha \xi_\beta}+
 g_{12} \!\sum\limits_{\alpha} \!g^{\alpha 1}\frac{\sigma\xi_1\xi_\alpha}{\sigma\sqrt{\sum\limits_{\alpha,\beta}\! g^{\alpha\beta} \xi_\alpha \xi_\beta}} - g_{11}  \sum\limits_{\alpha} g^{\alpha 1}\frac{\sigma\xi_2\xi_\alpha}{\sigma\sqrt{\sum\limits_{\alpha,\beta} g^{\alpha\beta} \xi_\alpha \xi_\beta}} \bigg)\nonumber\\
 &     \; \;  \; \;\,= \frac{1}{i\omega \mu \sqrt{|g|}} \!\bigg(\!\!-g_{12} \sqrt{\sum\limits_{\alpha,\beta} \!g^{\alpha\beta} \xi_\alpha \xi_\beta} + g_{12}\! \sum\limits_{\alpha}\! g^{\alpha 1}\frac{\xi_1\xi_\alpha}{\sqrt{\sum\limits_{\alpha,\beta} g^{\alpha\beta} \xi_\alpha \xi_\beta}} + g_{12}  \sum\limits_{\alpha} g^{\alpha 2}\frac{\xi_2\xi_\alpha}{\sqrt{\sum\limits_{\alpha,\beta} g^{\alpha\beta} \xi_\alpha \xi_\beta}}-\frac{\xi_1\xi_2}{ \sqrt{\sum\limits_{\alpha,\beta} g^{\alpha\beta} \xi_\alpha \xi_\beta}}\bigg) \nonumber\\
      &     \; \; \;\;\, = \frac{1}{i\omega \mu \sqrt{|g|}}\bigg( \!\!-g_{12} \sqrt{\sum\limits_{\alpha,\beta} g^{\alpha\beta} \xi_\alpha \xi_\beta}+ g_{12} \frac{\sum_{\alpha,\beta}g^{\alpha\beta} \xi_\alpha\xi_\beta}{\sqrt{\sum\limits_{\alpha,\beta} g^{\alpha\beta} \xi_\alpha \xi_\beta}}
     - \frac{\xi_1\xi_2}{ \sqrt{\sum\limits_{\alpha,\beta} g^{\alpha\beta} \xi_\alpha \xi_\beta}}\bigg)\nonumber\\
   &   \;\;\; \;\,= -\,\frac{\xi_1\xi_2}{i\omega \mu \,\sqrt{ \sum\limits_{\alpha,\beta}|g| g^{\alpha\beta}\xi_\alpha\xi_\beta}}. \nonumber\end{align}
Here the third equality follows from the relation  $g_{11} g^{11} +g_{12}g^{12}=1$ and $g_{11}g^{21}+g_{12} g^{22}=0$ (see (\ref{19.8.3;1})).
Since $\Lambda^{11}$ is uniquely determined by $\Lambda_{g,\Gamma}$, so the principal symbol $-\frac{\xi_1\xi_2}{i\omega \mu\sqrt{ \sum\limits_{\alpha,\beta}|g|g^{\alpha\beta} \xi_\alpha\xi_\beta}}$ of $\Lambda^{11}$ is uniquely determined by $\Lambda_{g,\Gamma}$ at each boundary point $x_0$. This implies that $|g| g^{\alpha\beta}$ are determined by $\Lambda_{g,\Gamma}$ at each boundary point $x_0$ for all $1\le \alpha,\beta\le 2$. However, $\mbox{det} (|g|g^{\alpha\beta})= |g|$, so the boundary values of $|g|g^{\alpha\beta}$  determine those of $|g|$, and we can thus recover the values of $g^{\alpha\beta}$ itself along $\Gamma$.

Next, recall that \begin{eqnarray*}  \wp(\nu \times H)=\frac{1}{i\omega \mu\sqrt{|g|}}  \left[ \begin{matrix} L^{11}  &  L^{12}\\
   L^{21}  & L^{22}   \end{matrix} \right] \begin{bmatrix} -g_{12} & - g_{22}\\
 g_{11} & g_{12}  \end{bmatrix} \wp(\nu\times E).\end{eqnarray*}
Since $g_{\alpha\beta}$ have been determined by the principal symbol $\psi^{11}_1$ of component $\Lambda^{11}$, it follows that $\Lambda$  determines
the operator $$\left[ \begin{matrix} L^{11}  &  L^{12}\\
   L^{21}  & L^{22}   \end{matrix} \right].$$
     From (\ref{19.7.21-7}) and (\ref{19.8.8-1}), we see that
     the terms of degree zero in $\xi'$ of the symbol of $L^{jj}$, ($j=1,2$), are
\begin{align} \label{19.9.6-6}  l_0^{jj}  =& \phi_0^{jj} -\bigg\{ \phi_1^{j3} \bigg( q_{-1}\Big(\sigma  \big(\phi_0^{3j}+\frac{1}{2} \sum\limits_{\alpha,\beta}g^{\alpha\beta}\frac{\partial g_{\alpha\beta}}{\partial x_j} \big)+\frac{\partial \sigma}{\partial x_j}\Big) +q_{-2} \Big(\sigma (i\xi_j +\phi_1^{3j} )\Big) -i\sum\limits_{m=1}^2 \frac{\partial q_{-1}}{\partial \xi_m} \frac{\partial \big(\sigma(i\xi_j+\phi_1^{3j})\big)}{\partial x_m} \bigg)
\\  & + \phi_0^{j3} q_{-1} \Big(\sigma(i\xi_j+\phi_1^{3j} )\Big) -i \sum\limits_{m=1}^2 \frac{\partial \phi_1^{j3} }{\partial \xi_m }\,\frac{\partial }{\partial x_m} \!\Big(q_{-1}\big(\sigma (i\xi_j+\phi_1^{3j} )\big)\Big) \bigg\} \nonumber\\
& + \sum\limits_{\alpha} g^{\alpha j} \frac{\partial}{\partial x_\alpha}\Big(q_{-1} \big(\sigma(i\xi_j +\phi_1^{3j} ) \big)\Big)
+\sum\limits_{\alpha}g^{\alpha j}\bigg( q_{-1} \Big(\sigma\big(\phi_0^{3j} +\frac{1}{2} \sum\limits_{\beta,\gamma}g^{\beta\gamma}\frac{\partial g_{\beta\gamma}}{\partial x_j}\big)+\frac{\partial \sigma}{\partial x_j}\Big) \nonumber\\
& +q_{-2} \Big(\sigma(i\xi_j +\phi_1^{3j} )\Big)\! -  i\sum\limits_{m=1}^2 \! \frac{\partial q_{-1} }{\partial \xi_m}\frac{\partial \big(\sigma(i\xi_j+\phi_1^{3j})\big) }{\partial x_m}\! \bigg)  i\xi_\alpha  + \sum\limits_{\alpha}
g^{\alpha j} \frac{\partial g_{\alpha j}}{\partial x_3}  \nonumber\\
     =& \phi_0^{jj} -\phi_0^{j3} q_{-1} \big(\sigma i\xi_j \big)+\sum\limits_{\alpha} g^{\alpha j}\frac{\partial }{\partial x_\alpha} \Big( q_{-1} (\sigma i\xi_j)\Big)  +\sum\limits_{\alpha} g^{\alpha j} \bigg(q_{-1} \Big(\sigma \big( \phi_0^{3j}
     +\frac{1}{2} \sum\limits_{\alpha,\beta} g^{\alpha\beta} \frac{\partial g_{\alpha\beta}}{\partial x_j} \big) +\frac{\partial \sigma}{\partial x_j}\Big) \nonumber\\
     &\, +q_{-2} \big(\sigma i\xi_j\big) -i\sum\limits_{m=1}^2 \frac{\partial q_{-1}}{\partial \xi_m} \frac{\partial }{\partial x_m} \big(\sigma i \xi_j\big)\bigg)
     \, i\xi_\alpha +\sum\limits_{\alpha} g^{\alpha j} \frac{\partial g_{\alpha j}}{\partial x_3} \nonumber\\
                     =&\phi_0^{jj}+\frac{
                     i \sum_{\alpha,\beta}g^{\alpha\beta} \Gamma_{3\alpha}^j \xi_\beta   }{\sum\limits_{\alpha,\beta} g^{\alpha\beta} \xi_\alpha\xi_\beta} i\xi_j -\frac{i\sum\limits_{\alpha} g^{\alpha j} \sum\limits_{\eta, \beta} g^{\eta \beta} \Gamma_{j\eta}^3 \xi_\beta}{\sum\limits_{\alpha,\beta} g^{\alpha\beta} \xi_\alpha\xi_\beta}i\xi_\alpha  +\sum\limits_{\alpha}  g^{\alpha j} q_{-2} \sigma (i\xi_j)(i\xi_\alpha)  +  \sum\limits_{\alpha}  g^{\alpha j} \frac{\partial g_{\alpha j}}{\partial x_3} \nonumber\\
                     &+T_0^{(1)} (g_{\alpha\beta}),\nonumber\end{align}
                     where and throughout the proof, each $T_0^{(s)}(g_{\alpha\beta})$ is an expressions involving only the boundary values of $g_{\alpha\beta}$, $g^{\alpha\beta}$, and their tangential derivatives.
The last two equalities follow from the fact that, for $1\le \gamma\le 2$, \begin{eqnarray}\label{19.8.29,3} & \;\;\quad\;\;\;\phi_1^{\gamma3}\!=\!\phi_1^{3\gamma}\!=\!0,\; \, \;\phi_0^{\gamma3}\!=\!-\frac{i \sum_{\alpha,\beta}g^{\alpha\beta} \Gamma_{3\alpha}^\gamma \xi_\beta}{ \sqrt{\sum_{\alpha,\beta} g^{\alpha\beta}\xi_\alpha\xi_\beta}}+T_0^{(2)}(g_{\alpha\beta}),\quad
   \phi_0^{3\gamma }\! =\! -\frac{i \sum_{\eta,\beta}g^{\eta\beta} \Gamma_{\gamma \eta}^3 \xi_\beta}{ \sqrt{\sum_{\alpha,\beta} \! g^{\alpha\beta}\xi_\alpha\xi_\beta}}+T_0^{(3)}(g_{\alpha\beta}).\end{eqnarray}
It follows from (\ref{19.8.4-15}) and (\ref{19.8.2-9}) that $q_{-2}$  and $\phi_0^{jk}$ can be written as  \begin{eqnarray} \label{19.8.12,7} q_{-2} = -\frac{1}{\sigma \sum\limits_{\alpha,\beta} g^{\alpha\beta} \xi_\alpha\xi_\beta} \big( \phi_0^{33} +\frac{1}{2} \sum\limits_{\alpha,\beta}g^{\alpha\beta}\frac{\partial g_{\alpha\beta}}{\partial x_3}\big)+ T_0^{(4)}(g_{\alpha\beta}),\end{eqnarray}
and \begin{eqnarray} \label{19.8.10-1}   \,\qquad\;\;\; \quad \phi_0^{jk}\! \!= \!\frac{1}{2\sqrt{\sum\limits_{\alpha\!,\beta}\!\! g^{\alpha \beta}\xi_\alpha \xi_\beta}} \!\bigg\{\!\!\bigg(\!\frac{1}{2} \!\sum\limits_{\alpha, \beta}\! g^{\alpha\beta}  \frac{\partial g_{\alpha\beta}}{\partial x_3}\!\bigg)\phi_1^{jk}\!
\!+ \!2\!\sum\limits_{l}\! \Gamma_{l3}^j \phi_1^{lk}\! +\!\frac{\partial \phi_1^{jk}}{\partial x_3} \!-\! 2i\!\sum\limits_{\alpha,\beta} \! g^{\alpha\beta}\Gamma_{\!k\alpha}^j \xi_\beta\! \bigg\}\!\!+\! T^{(5)}_0\!(g_{\alpha\beta}),\end{eqnarray}
respectively. Noting that $\phi_1^{jk}=0$ if $j\ne k$, we find by (\ref{19.9.6-6}),  (\ref{19.8.12,7}) and (\ref{19.8.10-1}) that  \begin{align*} & l_0^{11}\! +\!l_0^{22}
   = \phi_0^{11}
 +\phi_0^{22}\!-\! \frac{\sum\limits_{\alpha,\beta,\eta}
 g^{\alpha\beta} \Gamma_{3\alpha}^\eta \xi_\beta \xi_\eta}{\sum_{\alpha,\beta} g^{\alpha\beta} \xi_\alpha\xi_\beta} \!+\!\frac{\sum\limits_{\alpha,\beta,\gamma,\eta} g^{\alpha\gamma} g^{\eta \beta} \Gamma_{\gamma\eta}^3 \xi_\alpha \xi_\beta}
 {\sum_{\alpha,\beta} g^{\alpha\beta}\xi_\alpha\xi_\beta} \!-\!q_{-\!2} \sigma \sum\limits_{\alpha,\beta} g^{\alpha\beta} \xi_\alpha \xi_\beta \!+\!
 \sum\limits_{\alpha, \beta} g^{\alpha\beta}\frac{\partial g_{\alpha\beta}}{\partial x_3} \! +\!T_{0}^{(6)} (g_{\alpha\beta}).\end{align*}
  Since  \begin{eqnarray*}
 && \Gamma_{3\alpha}^\eta = \frac{1}{2} \sum_{\rho} g^{\eta \rho} \big(\frac{\partial g_{\alpha \rho}}{\partial x_3} +\frac{\partial g_{3 \rho}}{\partial x_\alpha} -\frac{\partial g_{3\alpha}}{\partial x_\rho}\big)=\frac{1}{2}\sum_{\rho} g^{\eta \rho}
\frac{\partial g_{\alpha\rho}}{\partial x_3} +T_0^{(7)}(g_{\alpha\beta}),\\
&&\Gamma^3_{\gamma\eta}= \frac{1}{2} \sum_{s=1}^3 g^{3s} \big(\frac{\partial g_{\eta s}}{\partial x_\gamma} +\frac{\partial g_{\gamma s}}{\partial x_\eta} -\frac{\partial g_{\gamma\eta}}{\partial x_s}\big)=-\frac{1}{2}
\frac{\partial g_{\gamma\eta}}{\partial x_3} +T_0^{(8)}(g_{\alpha\beta}),
\end{eqnarray*}
 and \begin{eqnarray*} \sum\limits_{\gamma,\eta} g^{\alpha\gamma} \frac{\partial g_{\gamma \eta}}{\partial x_3} g^{\eta\beta} =- \frac{\partial g^{\alpha\beta}}{\partial x_3},\end{eqnarray*}
 we see by applying (\ref{19.8.12,7}) that
 \begin{align*} \! \!\! l_0^{11}+l_0^{22}&= \phi_0^{11} +\phi_0^{22}
  -\frac{\sum_{\alpha,\beta,\eta,\rho}\frac{1}{2} g^{\alpha \beta} g^{\eta \rho} \frac{\partial g_{\alpha\rho}}{\partial x_3} \xi_\beta \xi_\eta}
  {\sum_{\alpha,\beta} g^{\alpha\beta} \xi_\alpha\xi_\beta} -\frac{\sum_{\alpha,\beta,\gamma,\eta} \frac{1}{2} g^{\alpha\gamma} \frac{\partial g_{\gamma\eta }}{\partial x_3} g^{\eta\beta} \xi_\alpha \xi_\beta}{\sum_{\alpha,\beta} g^{\alpha\beta} \xi_\alpha\xi_\beta} \\
 \!\!\!  &   \;\quad  -\bigg(\sum\limits_{\alpha,\beta} g^{\alpha\beta}\xi_\alpha\xi_\beta \bigg)\bigg(\frac{-1}{\sum_{\alpha,\beta} g^{\alpha\beta}\xi_\alpha\xi_\beta }\bigg) \big(   \phi_0^{33} +\frac{1}{2} \sum\limits_{\alpha,\beta} g^{\alpha\beta}\frac{\partial g_{\alpha\beta}}{\partial x_3}\big)+ \sum\limits_{\alpha,\beta} g^{\alpha\beta}\frac{\partial g_{\alpha\beta}}{\partial x_3}
  +T_0^{(9)}(g_{\alpha\beta})\\
 \!\!\!& = \phi_0^{11}+\phi_0^{22} +\phi_0^{33} -\frac{\sum_{\beta,\eta} \big(\!\!-\!\frac{1}{2} \frac{\partial g^{\beta\eta}}{\partial x_3} \xi_\beta \xi_\eta\big)}{\sum_{\alpha,\beta} g^{\alpha\beta} \xi_\alpha\xi_\beta}-\frac{\sum_{\alpha,\beta} \big(\!\!-\!\frac{1}{2} \frac{\partial g^{\alpha\beta}}{\partial x_3} \xi_\alpha \xi_\beta\big)}{\sum_{\alpha,\beta} g^{\alpha\beta} \xi_\alpha\xi_\beta}
  + \frac{3}{2} \sum\limits_{\alpha,\beta} g^{\alpha\beta}\frac{\partial g_{\alpha\beta}}{\partial x_3}
 +T_0^{(9)}(g_{\alpha\beta})\\
  \!\!\! &    =\phi_0^{11}+\phi_0^{22} +\phi_0^{33} +\frac{\sum_{\alpha,\beta}  \frac{\partial g^{\alpha\beta}}{\partial x_3} \xi_\alpha \xi_\beta}{\sum_{\alpha,\beta} g^{\alpha\beta} \xi_\alpha\xi_\beta} +  \frac{3}{2} \sum\limits_{\alpha,\beta} g^{\alpha\beta}\frac{\partial g_{\alpha\beta}}{\partial x_3} +T_0^{(9)}(g_{\alpha\beta}).\end{align*}
But, it follows from (\ref{19.8.10-1}) and $\phi_1^{kk}=\sqrt{\sum_{\alpha,\beta} g^{\alpha\beta} \xi_\alpha\xi_\beta}$ for $k=1,2,3$ that \begin{eqnarray*}\!\!\!&&  \phi_0^{11}+\phi_0^{22}+\phi_0^{33} =\frac{1}{2\sqrt{\sum\limits_{\alpha,\beta} g^{\alpha\beta}\xi_\alpha\xi_\beta}}\, \sum\limits_{k=1}^3 \bigg(
\frac{1}{2} \sum\limits_{\alpha,\beta} g^{\alpha\beta}  \frac{\partial g_{\alpha\beta}}{\partial x_3} \phi_1^{kk} + 2 \sum\limits_{l=1}^3
\Gamma_{l3}^k \phi_{1}^{lk} +\frac{\partial \phi_1^{kk}}{\partial x_3}\bigg) +T_0^{(10)} (g_{\alpha\beta})\\
 \!\!\!&& = \frac{1}{2\sqrt{\sum\limits_{\alpha,\beta} g^{\alpha\beta}\xi_\alpha\xi_\beta}}\bigg\{ \bigg( \frac{3}{2} \sum\limits_{\alpha,\beta} g^{\alpha\beta}  \frac{\partial g_{\alpha\beta}}{\partial x_3}  +  \sum\limits_{\alpha,\beta} g^{\alpha\beta}  \frac{\partial g_{\alpha\beta}}{\partial x_3}\bigg) \sqrt{\sum_{\alpha,\beta} g^{\alpha\beta} \xi_\alpha\xi_\beta} + \frac{3 \sum_{\alpha,\beta} \frac{\partial g^{\alpha\beta}}{\partial x_3} \xi_\alpha\xi_\beta}{2\sqrt{\sum_{\alpha,\beta} g^{\alpha\beta}\xi_\alpha\xi_\beta} }\bigg\} +T_0^{(10)} (g_{\alpha\beta})\\
\!\!\! &&= \frac{5}{4} \sum\limits_{\alpha,\beta} g^{\alpha\beta} \frac{\partial g_{\alpha\beta}}{\partial x_3} + \frac{3}{4}\,\frac{\sum_{\alpha,\beta} \frac{\partial g^{\alpha\beta} }{\partial x_3} \xi_\alpha\xi_\beta}{\sum_{\alpha,\beta} g^{\alpha\beta}\xi_\alpha\xi_\beta}
  +T_0^{(10)} (g_{\alpha\beta}),\end{eqnarray*} where, in the second equality we have used  $\sum_{k=1}^3 \Gamma_{k3}^k =\frac{1}{2} \sum_{{}_{\alpha,\beta}} g^{\alpha\beta} \frac{\partial g_{\alpha\beta}}{\partial x_3}$.
Hence we have
\begin{align*} l_0^{11} +l_0^{22}
 = &  \frac{11}{4} \sum\limits_{\alpha,\beta} g^{\alpha\beta}\frac{\partial g_{\alpha\beta}}{\partial x_3}+\frac{7}{4}\,\frac{\sum_{\alpha,\beta} \frac{\partial g^{\alpha\beta} }{\partial x_3} \xi_\alpha\xi_\beta}{\sum_{\alpha,\beta} g^{\alpha\beta}\xi_\alpha\xi_\beta} +T^{(11)}_0(g_{\alpha\beta}).\end{align*}
  In view of  \begin{eqnarray}\label{19.8.6;7}\sum\limits_{\alpha,\beta} g^{\alpha\beta} g_{\alpha\beta}=2,\end{eqnarray} we have \begin{eqnarray*} \sum\limits_{\alpha,\beta} \bigg( \frac{\partial g^{\alpha\beta}}{\partial x_3} g_{\alpha\beta} +g^{\alpha\beta} \frac{\partial g_{\alpha\beta}}{\partial x_3}\bigg)
 =0,\end{eqnarray*}  i.e., \begin{eqnarray} \label{19.8.6,1} \sum\limits_{\alpha,\beta} g^{\alpha\beta} \frac{\partial g_{\alpha\beta}}{\partial x_3} =- \sum\limits_{\alpha,\beta} g_{\alpha\beta} \frac{\partial g^{\alpha\beta}}{\partial x_3}.\end{eqnarray}  Therefore, \begin{eqnarray*} l_{0}^{11} +l_{0}^{22} = -\frac{11}{4} \sum\limits_{\alpha,\beta}
 g_{\alpha\beta} \frac{\partial g^{\alpha\beta}}{\partial x_3} +\frac{7}{4} \frac{\sum_{\alpha,\beta}\frac{\partial g^{\alpha\beta}}{\partial x_3}\xi_\alpha\xi_\beta }{\sum_{\alpha,\beta} g^{\alpha \beta} \xi_\alpha\xi_\beta} +T^{(11)}_0(g_{\alpha\beta}).\end{eqnarray*}
 If we set $h_1^{\alpha\beta} =\frac{\partial g^{\alpha\beta}}{\partial x_3}$ and $h_1=\sum\limits_{\alpha,\beta} g_{\alpha\beta}h_1^{\alpha\beta}$, then
 \begin{align}  & l_0^{11}+l_0^{22} = \frac{1}{4\sum\limits_{\alpha,\beta} g^{\alpha\beta} \xi_\alpha\xi_\beta} \, \sum\limits_{\alpha,\beta} \big( 7h^{\alpha\beta}_1 -11h_1 g^{\alpha\beta}\big) \xi_\alpha\xi_\beta +T^{(11)}_0(g_{\alpha\beta}).\end{align}
  Evaluating this on unit vectors $\xi\in T^*(\Gamma)$ implies that
   $l_0^{11}+l_0^{22}$ and the values of $g_{\alpha\beta}$ along $\partial M$ completely determine the quadratic form
        \begin{eqnarray}\label{19.9.20,2}  k^{\alpha\beta}_1 =7h^{\alpha\beta}_1 -11h_1 g^{\alpha\beta} \end{eqnarray}
        along $\Gamma$.
          From (\ref{19.9.20,2}), we have  \begin{eqnarray*} \sum\limits_{\alpha,\beta} g_{\alpha\beta} k^{\alpha\beta}_1 = 7\sum\limits_{\alpha,\beta} g_{\alpha\beta} h^{\alpha\beta}_1 -11h_1 \sum\limits_{\alpha,\beta} g_{\alpha\beta} g^{\alpha\beta} = 7h_1- 22h_1=-15h_1,\end{eqnarray*}
so \begin{eqnarray} h_1= - \frac{1}{15} \sum\limits_{\alpha,\beta} g_{\alpha\beta} k^{\alpha\beta}_1,\end{eqnarray}
  and hence  \begin{eqnarray} \frac{\partial g^{\alpha\beta}}{\partial x_3}= h^{\alpha\beta}_1 =  \frac{1}{7} \bigg( k^{\alpha\beta}_1 -\big(\frac{11}{15}
    \sum\limits_{\gamma,\rho} g_{\gamma\rho} k^{\gamma\rho}_1\big)g^{\alpha\beta}\bigg).\end{eqnarray}
This shows that $\frac{\partial g^{\alpha\beta}}{\partial x_3}$ on $\partial M$ are uniquely determined by the symbol of degree zero of $L^{11}+L^{22}$ as well as  the symbol of order $1$ of $\Lambda^{11}$ along $\Gamma$, and hence $\frac{\partial g^{\alpha\beta}}{\partial x_3}$ on $\Gamma$ are uniquely determined by $\Lambda_{g,\Gamma}$.

Now, we calculate the terms of degree $-1$ in $\xi'$ of the symbol of $L^{11}+L^{22}$. It follows from (\ref{19.7.21-7}), (\ref{19.8.8-1}) and the symbol formula for the product of two pseudodifferential operators that
\begin{align}\label{2019.12.6-1} &     l^{jj}_{-1} =    \phi^{jj}_{-1} - \bigg\{ \phi_1^{j3} \bigg( q_{-1} \big(\sigma \phi_{-1}^{3j}\big) +q_{-2} \Big( \sigma\big(\phi_0^{3j} +\frac{1}{2} \sum\limits_{\alpha,\beta} g^{\alpha\beta} \frac{\partial g_{\alpha\beta}}{\partial x_j}\big) +\frac{\partial \sigma}{\partial x_j} \Big)  +q_{-3} \sigma (i\xi_j +\phi_1^{3j} )\\
 &       \quad \quad\;  \, -i \sum\limits_{m=1}^2
\frac{\partial q_{-1}}{\partial \xi_m} \, \frac{\partial }{\partial x_m}\! \Big(\sigma\big( \phi_0^{3j} +\frac{1}{2} \sum\limits_{\alpha,\beta} g^{\alpha\beta} \frac{\partial g_{\alpha\beta}}{\partial x_j}\big)+\frac{\partial \sigma}{\partial x_j}\Big) -i\sum\limits_{m=1}^2 \frac{\partial q_{-2}}{\partial \xi_m} \, \frac{\partial} {\partial x_m}\big( \sigma( i\xi_j+\phi_1^{3j})\big)\nonumber\\
&  \quad \quad\;  \, -
 \frac{1}{2} \sum\limits_{m,l} \frac{\partial^2 q_{-1}}{\partial \xi_m \partial \xi_l}\, \frac{\partial^2 }{\partial x_m\partial x_l}\big(\sigma(i\xi_j+\phi_1^{3j})\big)  \bigg)  + \phi_0^{j3} \bigg( q_{-1} \Big(\sigma\big( \phi_0^{3j} + \frac{1}{2} \sum\limits_{\alpha,\beta} g^{\alpha\beta} \frac{\partial g_{\alpha\beta}}{\partial x_j}\big)+\frac{\partial \sigma}{\partial x_j}\Big) \nonumber\\
 &  \quad \quad\;\, +q_{-2} \big(\sigma(i\xi_j+\phi_1^{3j} ) \big)  -i  \sum\limits_{m=1}^2 \frac{\partial q_{-1}}{\partial \xi_m} \, \frac{\partial  }{\partial x_m}\big(\sigma (i\xi_j+\phi_1^{3j})\big)
\bigg)  +\phi_{-1}^{j3} q_{-1} \big(\sigma(i\xi_j+ \phi_1^{3j})\big)\nonumber\\
&    \quad \quad \,\; - \!i\! \sum\limits_{m=1}^2 \!\frac{\partial \phi_1^{j3} }{\partial \xi_m}  \frac{\partial}{\partial x_m}\!\bigg(\!
q_{-1} \!\Big(\!\sigma\big( \phi_0^{3j} \!+\! \frac{1}{2}\! \sum\limits_{\alpha,\beta}\! g^{\alpha\beta} \frac{\partial g_{\alpha\beta}}{\partial x_j}\big)\!+\!\frac{\partial \sigma}{\partial x_j}\!\Big) \!+\!q_{-2} \!\big(\sigma(i\xi_j\!+\!\phi_1^{3j} )\big) \!-\!i \! \sum\limits_{m=1}^2\! \frac{\partial q_{-1}}{\partial \xi_m}  \frac{\partial }{\partial x_m}\!\big(\sigma (i\xi_j +\phi_1^{3j})\big)\!\bigg) \nonumber\\
&     \quad \quad\, \; -\frac{1}{2} \sum\limits_{m,l=1}^2
\frac{\partial^2 \phi_1^{j3} }{\partial \xi_m\partial \xi_l}\, \frac{\partial^2 }{\partial x_m\partial x_j} \! \big(q_{-1}\sigma(i\xi_j+\phi_1^{3j})\big)-i \sum\limits_{m=1}^2 \frac{\partial \phi_0^{j3}}{\partial \xi_m} \frac{\partial}{\partial x_m} \!\big( q_{-1}\sigma(i\xi_j+\phi_1^{3j})\big)\bigg\}\nonumber\\
& \quad \quad \;\,  + \! \sum\limits_{\alpha}\! g^{\alpha j}\! \frac{\partial }{\partial x_\alpha}\!\bigg\{\!q_{-1} \Big(\sigma\big(\phi_{0}^{3j}\!+\!\frac{1}{2}\!\sum\limits_{\alpha,\beta}\! g^{\alpha\beta}\frac{\partial g_{\alpha\beta}}{\partial x_j} \big)+\frac{\partial\sigma}{\partial  x_j}\Big) \!+\!q_{-2} \big(\sigma(i\xi_j\!+\! \phi_1^{3j})\big)  \! -\! i \!\sum\limits_{m=1}^2\! \frac{\partial q_{-1}}{\partial \xi_m} \, \frac{\partial  }{\partial x_m}\big(\sigma(i\xi_j\!+\!\phi_1^{3j})\big)\!\bigg\}\nonumber\\
 &  \quad \quad\; \; +  \!\sum\limits_{\alpha}\! g^{\alpha j} \bigg\{ q_{-1} \sigma\phi_{-1}^{3j} +q_{-2} \Big(\sigma\big(\phi_0^{3j} +\frac{1}{2} \!\sum\limits_{\alpha,\beta} \!g^{\alpha\beta}
\frac{\partial g_{\alpha}}{\partial x_j}\big)+\frac{\partial \sigma}{\partial x_j}\Big) +q_{-3} \sigma(i\xi_j +\phi_1^{3j})\nonumber\\
& \quad \quad \;\, -i\! \sum\limits_{m=1}^2 \!\frac{\partial q_{-1}}{\partial \xi_m} \, \frac{\partial }{\partial x_m}\!\Big(\sigma\big(\phi_0^{3j}\! +\!\frac{1}{2}\! \sum\limits_{\alpha,\beta}\! g^{\alpha\beta} \frac{\partial g_{\alpha\beta}}{\partial x_j}\big)+\frac{\partial \sigma}{\partial x_j}\Big)   -\! i \!\sum\limits_{m=1}^2 \!\frac{\partial q_{-2}}{\partial \xi_m} \, \frac{\partial }{\partial x_m}\big(\sigma(i\xi_j +\phi_1^{3j})\big) \nonumber\\
&    \quad  \quad\;\, -\frac{1}{2}\! \sum\limits_{m,l}^2 \!\frac{\partial^2 q_{-1}}{\partial \xi_m \partial \xi_l}\, \frac{\partial^2 }{\partial x_m\partial x_j}\big( \sigma(i\xi_j+\phi_1^{3j})\big) \! \bigg\}i\xi_\alpha     \quad \,\mbox{for}\,\, j=1,2.\nonumber\end{align}
 Clearly, there is not other terms containing $\frac{\partial^2 g_{\alpha\beta}}{\partial x_3^2}$ except for $\phi_{-1}^{jj}$, $\phi_{-1}^{3j}$, $\phi_{-1}^{j3}$ and $q_{-3}$ on the right-hand side of the above equality. Thus
 \begin{eqnarray*} l^{jj}_{-1}=   \phi_{-1}^{jj} - \phi_{-1}^{j3}  q_{-1} \sigma\,i\xi_j+ \sum\limits_{\alpha} g^{\alpha j}\! \bigg(\! q_{-1} \sigma \phi_{-1}^{3j}+ q_{-3}( \sigma i\xi_j)
          \!\bigg)i\xi_\alpha  +T^{(1)}_{-1}(g_{\alpha\beta}), \quad\, j=1,2, \end{eqnarray*}
        where  each $T_{-1}^{(s)}(g_{\alpha\beta})$ only involves the boundary values of $g^{\alpha\beta}$, $g_{\alpha\beta}$, and their normal derivatives of order at most one. By (\ref{19.8.4-16,}) we have \begin{eqnarray*} q_{-3} = -\frac{1}{\sigma\phi_1^{33}} q_{-1} \sigma \phi_{-1}^{33} +T^{(2)}_{-1} (g_{\alpha\beta}),\end{eqnarray*}
   so that \begin{eqnarray*}   l_{-1}^{jj} = \phi_{-1}^{jj} -\phi_{-1}^{j3} q_{-1}\sigma i\xi_j+\sum\limits_{\alpha} g^{\alpha j}  q_{-1} \sigma \phi_{-1}^{3j} i\xi_\alpha + \sum\limits_{\alpha} g^{\alpha j} \bigg(\frac{\xi_j\xi_\alpha}{\phi_1^{33} } \bigg)q_{-1}\sigma \phi_{-1}^{33}   + T^{(3)}_{-1}(g_{\alpha\beta}), \;\; \, j=1,2.\end{eqnarray*}
It follows that
\begin{align} \label{19.8.12;1} & \!\!\;\, l^{11}_{-1}+ l^{22}_{-1} =  \phi_{-1}^{11} +\phi_{-1}^{22}
- \!\sum\limits_{\beta} \phi_{-1}^{\beta 3} q_{-1}\sigma i\xi_\beta+\sum\limits_{\alpha, \beta} g^{\alpha\beta} q_{-1}\sigma \phi_{-1}^{3\beta} \, i\xi_\alpha +  \sum\limits_{\alpha,\beta} g^{\alpha \beta} \bigg(\frac{\xi_\beta\xi_\alpha}{\phi_1^{33} } \bigg)q_{-1} \sigma\phi_{-1}^{33}  \!+\! T^{(4)}_{-1}(g_{\alpha\beta})\\
 &\qquad \qquad  \, =   \! \phi_{-1}^{11}\! +\!\phi_{-1}^{22}\!+\! \phi_{-1}^{33}\!-\! \frac{\sum_{\beta}  \phi_{-1}^{\beta 3}  i\xi_\beta}{\sqrt{\sum_{\alpha,\beta} g^{\alpha\beta} \xi_\alpha \xi_\beta }}  \!+\!\frac{\sum_{\alpha, \beta} g^{\alpha \beta} \phi_{-1}^{3\beta}  i\xi_\alpha}{\sqrt{\sum_{\alpha,\beta} g^{\alpha\beta} \xi_\alpha \xi_\beta }}\! + \!  T^{(4)}_{-1}(g_{\alpha\beta}).\nonumber \end{align}
From (\ref{19.8.2-10}) and (\ref{19.8.10-1})  we have
   \begin{align*}    \phi_{-1}^{kk} =& \frac{1}{2\sqrt{\sum\limits_{\alpha,\beta} g^{\alpha\beta} \xi_\alpha\xi_\beta}}\left\{
 -(\phi_0^2)^{kk} +i \sum\limits_{l=1}^3\sum\limits_{m=1}^2 \big( \frac{\partial \phi_1^{kl}} {\partial \xi_m} \, \frac{\partial \phi_0^{lk}}{\partial x_m} +\frac{\partial \phi_0^{kl}}{\partial \xi_m}\, \frac{\partial \phi_{1}^{lk}}{\partial x_m} \big) +\frac{1}{2} \sum\limits_{m,l=1}^2 \frac{\partial^2 \phi_1^{kk}}{\partial \xi_m\partial \xi_l} \, \frac{\partial^2 \phi_1^{kk}}{\partial x_m\partial x_l} \right.\\
 &       \left. + \big(\frac{1}{2} \sum\limits_{\alpha,\beta} g^{\alpha\beta} \frac{\partial g_{\alpha\beta}}{\partial x_3}\big) \phi^{kk}_0 + 2\sum\limits_{l=1}^3 \Gamma_{l3}^k  \phi_0^{lk} + \frac{\partial \phi_0^{kk}}{\partial x_3}-a_{kk} +  R_k^k -\omega^2 \mu \sigma \right\}+ T^{(5)}_{-1}(g_{\alpha\beta})\\
 = & \frac{1}{2\sqrt{\sum\limits_{\alpha,\beta} g^{\alpha\beta} \xi_\alpha\xi_\beta}}\bigg(\frac{\partial \phi_0^{kk}}{\partial x_3}-a_{kk}+R^k_k\bigg) +T^{(6)}_{-1}(g_{\alpha\beta}), \quad  \; k=1,2,3.\end{align*}
According to the definitions of $a_{kk}$ and $R^k_k$, we have that for $k=1,2,3$ (see (\ref{19.9.13-1}) and (\ref{19.8.29-2})),  \begin{eqnarray} \label{19.9.1-4} && a_{kk}
= \sum\limits_{l,m=1}^3 g^{ml} \bigg(\frac{\partial \Gamma_{km}^k}{\partial x_l} + \sum\limits_{h}\Gamma_{hl}^k \Gamma_{km}^h -\sum\limits_{h} \Gamma_{hk}^k\Gamma_{ml}^h \bigg)\\
&& \quad \;\;\; =\sum\limits_{l,m=1}^3g^{ml} \frac{\partial
\Gamma_{km}^k}{\partial x_l} +T^{(7)}_{-1}(g_{\alpha\beta})\nonumber\\
&&  \quad \;\;\; =
 \frac{\partial \Gamma_{k3}^k }{\partial x_3} +  T^{(7)}_{-1}(g_{\alpha\beta})\nonumber\\
 && \quad \;\;\; =
 \frac{\partial }{\partial x_3} \bigg( \frac{1}{2} \sum\limits_{s=1}^3  g^{ks}\big( \frac{\partial g_{ks}}{\partial x_3} +\frac{\partial g_{3s}}{\partial x_k} -\frac{\partial g_{k3}}{\partial x_s}\big)\bigg) + T^{(7)}_{-1}(g_{\alpha\beta})\nonumber\\
 &&   \quad \;\;\; =\left\{ \begin{array} {ll}     \frac{1}{2} \sum\limits_{\beta} g^{k\beta}\frac{\partial^2 g_{k\beta}}{\partial x_3^2}+T^{(8)}_{-1}(g_{\alpha\beta}) \quad &\mbox{when}\;\,k=1,2,\\
 0 \quad &\mbox{when}\;\, k=3,\end{array}\right.  \nonumber \\
   \label{19.9.1-5} &&  R_k^k  = \sum\limits_{m=1}^3 g^{km} R_{mk}\\
 &&     \;\; \quad =  \sum\limits_{l,m=1}^3 g^{km} \bigg(\frac{\partial \Gamma_{mk}^l}{\partial x_l} -\frac{\partial \Gamma_{ml}^l}{\partial x_k}\bigg)+T^{(9)}_{-1}(g_{\alpha\beta})\nonumber\\
                &&   \; \;\quad =\left\{ \!\begin{array}{ll}
             \! -\frac{1}{2} \sum_{\beta} g^{k\beta} \frac{\partial^2 g_{k\beta}}{\partial x_3^2}+ T^{(10)}_{-1}(g_{\alpha\beta})& \quad \mbox{when}\,\, k=1,2,\\
             -\frac{1}{2} g^{\alpha\beta}\frac{\partial^2 g_{\alpha\beta}}{\partial x_3^2}   +T^{(11)}_{-1}(g_{\alpha\beta})  & \quad \mbox{when}\,\, k=3.  \end{array}\right.\nonumber\end{eqnarray}
            Also, \begin{eqnarray} \label{19.9.3-1} \left.\begin{array}{ll}   a_{3\beta}= 0+ T_{-1}^{(12)} (g_{\alpha\beta}), \quad a_{\beta 3} =  0+ T_{-1}^{(13)} (g_{\alpha\beta}),\\
             R^{3}_\beta = 0 +T_{-1}^{(14)} (g_{\alpha\beta}), \quad R^\beta_3 = 0 +T_{-1}^{(15)} (g_{\alpha\beta}).\,\,\end{array}\right.\end{eqnarray}
              (\ref{19.9.1-4}) and (\ref{19.9.1-5})  imply that \begin{align*}\sum_{k=1}^3 (R_k^k-a_{kk})&=
       -\frac{3}{2} \sum_{\alpha,\beta} g^{\alpha \beta}\frac{\partial^2 g_{\alpha\beta}}{\partial x_3^2}+T^{(16)}_{-1}(g_{\alpha\beta}).\end{align*}
 It follows from (\ref{19.8.2-9} ) that  \begin{align*}  \frac{\partial \phi^{kk}_0}{\partial x_3} & = \frac{\partial }{\partial x_3}
 \bigg\{\frac{1}{2\sqrt{\sum\limits_{\alpha,\beta} g^{\alpha\beta} \xi_\alpha\xi_\beta}} \bigg[\bigg( \frac{1}{2} \sum\limits_{\alpha,\beta}
 g^{\alpha\beta} \frac{\partial g_{\alpha\beta}}{\partial x_3}  + 2 \Gamma_{k3}^k \bigg)\sqrt{\sum\limits_{\alpha,\beta} g^{\alpha\beta}\xi_\alpha\xi_\beta}
 +\frac{\partial}{\partial x_3} \sqrt{\sum\limits_{\alpha,\beta} g^{\alpha\beta} \xi_\alpha\xi_\beta} \bigg]\bigg\} +T^{(17)}_{-1}(g_{\alpha\beta})\\
    & = \frac{\partial }{\partial x_3}
 \bigg\{ \bigg( \frac{1}{4} \sum\limits_{\alpha,\beta}
 g^{\alpha\beta} \frac{\partial g_{\alpha\beta}}{\partial x_3}  +  \Gamma_{k3}^k \bigg)
 +\frac{\sum_{\alpha,\beta}\frac{\partial g^{\alpha\beta}}{\partial x_3}\xi_\alpha\xi_\beta}{4\sum_{\alpha,\beta} g^{\alpha\beta} \xi_\alpha\xi_\beta}\bigg\} +T^{(17)}_{-1}(g_{\alpha\beta}) \quad \, \mbox{for}\;\, k=1,2,3.\end{align*}
       Therefore  \begin{align} \label{19.9.3,3} & \; \phi_{-1}^{11}+\phi_{-1}^{22} +\phi_{-1}^{33} =  \frac{1}{2\sqrt{\sum\limits_{\alpha,\beta} g^{\alpha\beta} \xi_\alpha\xi_\beta}}\bigg\{  \frac{\partial }{\partial x_3}
 \bigg[\bigg( \frac{3}{4} \sum\limits_{\alpha,\beta}
 g^{\alpha\beta} \frac{\partial g_{\alpha\beta}}{\partial x_3}  +  \sum\limits_{k=1}^3\Gamma_{k3}^k \bigg)
 +\frac{3}{4}\,\frac{\sum_{\alpha,\beta}\frac{\partial g^{\alpha\beta}}{\partial x_3}\xi_\alpha\xi_\beta}{\sum_{\alpha,\beta} g^{\alpha\beta} \xi_\alpha\xi_\beta}\;\bigg]\\
 &\quad  \quad \, -\frac{3}{2} \sum_{\alpha,\beta} g^{\alpha \beta}\frac{\partial^2 g_{\alpha\beta}}{\partial x_3^2}
  \bigg\}
 +\!T^{(18)}_{-1}(g_{\alpha\beta})\nonumber\\
 &\quad =\frac{1}{2\sqrt{\sum\limits_{\alpha,\beta} g^{\alpha\beta} \xi_\alpha\xi_\beta}}\!\bigg\{\!  \frac{\partial }{\partial x_3}
 \bigg[ \frac{5}{4} \sum\limits_{\alpha,\beta}
 g^{\alpha\beta} \frac{\partial g_{\alpha\beta}}{\partial x_3}
 \!+\!\frac{3}{4} \frac{\sum_{\alpha,\beta} \frac{\partial g^{\alpha\beta}}{\partial x_3}\xi_\alpha\xi_\beta}{\sum_{\alpha,\beta} g^{\alpha\beta}\xi_\alpha\xi_\beta}  \;\bigg]\!-\frac{3}{2} \sum_{\alpha,\beta} g^{\alpha \beta}\frac{\partial^2 g_{\alpha\beta}}{\partial x_3^2}
  \!\bigg\}
\! +\!T^{(18)}_{-1}(g_{\alpha\beta})\nonumber\\
   & \quad  = \frac{1}{2\sqrt{\sum\limits_{\alpha,\beta} g^{\alpha\beta} \xi_\alpha\xi_\beta}}\bigg\{
 \bigg[- \frac{1}{4} \sum\limits_{\alpha,\beta}
 g^{\alpha\beta} \frac{\partial^2 g_{\alpha\beta}}{\partial x_3^2}
 +\frac{3}{4} \frac{\sum_{\alpha,\beta} \frac{\partial^2 g^{\alpha\beta}}{\partial x_3^2}\xi_\alpha\xi_\beta}{\sum_{\alpha,\beta} g^{\alpha\beta}\xi_\alpha\xi_\beta}  \;\bigg]
  \bigg\}
 +\!T^{(19)}_{-1}(g_{\alpha\beta})\nonumber
 \\  &  \quad  = \frac{1}{2\sqrt{\sum\limits_{\alpha,\beta} g^{\alpha\beta} \xi_\alpha\xi_\beta}}\bigg\{
 \bigg[ \frac{1}{4} \sum\limits_{\alpha,\beta}
 g_{\alpha\beta} \frac{\partial^2 g^{\alpha\beta}}{\partial x_3^2}
 +\frac{3}{4} \frac{\sum_{\alpha,\beta} \frac{\partial^2 g^{\alpha\beta}}{\partial x_3^2}\xi_\alpha\xi_\beta}{\sum_{\alpha,\beta} g^{\alpha\beta}\xi_\alpha\xi_\beta}  \;\bigg]
  \bigg\}
 +\!T^{(20)}_{-1}(g_{\alpha\beta}).\nonumber
 \end{align}
The last equality follows from the relation $\sum_{{}_{\alpha, \beta}} g^{\alpha\beta} \frac{\partial^2 g_{\alpha\beta}}{\partial x_3^2} =- \sum_{{}_{\alpha,\beta}}
  g_{\alpha\beta} \frac{\partial^2 g^{\alpha\beta}}{\partial x_3^2} +T^{(21)}_{-1} (g_{\alpha\beta})$ which can be derived by the equality $\frac{\partial}{\partial x_3} \big(\sum_{{}_{\alpha,\beta}} g^{\alpha\beta}  \frac{\partial g_{\alpha\beta}}{\partial x_3}\big) =-\frac{\partial }{\partial x_3} \big(\sum_{\alpha,\beta} g_{\alpha\beta}  \frac{\partial g^{\alpha\beta}}{\partial x_3}\big)$.
  In addition, from (\ref{19.8.2-9}), (\ref{19.8.2-10}) and (\ref{19.9.3-1}) we have
 \begin{align*}
\phi_{-1}^{\beta 3} &= \frac{1}{2\sqrt{\sum_{\alpha,\beta} g^{\alpha\beta}\xi_\alpha\xi_\beta}} \bigg(\frac{\partial \phi_0^{\beta 3}}{\partial x_3}-a_{\beta 3}+  R^{\beta}_3\bigg)  +T_{-1}^{(22)}(g_{\alpha\beta}) \\
  &= \frac{1 }{2\sqrt{\sum_{\alpha,\beta} g^{\alpha\beta}\xi_\alpha\xi_\beta}} \bigg(\frac{\partial \phi_0^{\beta 3}}{\partial x_3}\bigg)  +T_{-1}^{(23)}(g_{\alpha\beta}) \\
&=\frac{1}{2\sqrt{\sum_{\alpha,\beta} g^{\alpha\beta} \xi_\alpha\xi_\beta}} \bigg\{\frac{ \partial }{\partial x_3}  \bigg[\frac{- 2i\sum_{\gamma,\eta} g^{\eta \gamma} \Gamma_{3\eta}^\beta\,  \xi_\gamma}{2 \sqrt{\sum_{\alpha,\beta} g^{\alpha\beta}\xi_\alpha\xi_\beta}}\bigg] \bigg\}+T_{-1}^{(24)}(g_{\alpha\beta}) \\
&= \frac{1}{2\sqrt{\sum_{\alpha,\beta} g^{\alpha\beta} \xi_\alpha\xi_\beta}} \bigg\{\frac{\partial }{\partial x_3} \bigg[\frac{-i  \sum_{\eta,\gamma,\rho}
g^{\eta\gamma}  g^{\beta \rho} \frac{\partial g_{\eta\rho}}{\partial x_3}  \xi_\gamma} {2\sqrt{\sum_{\alpha,\beta}
g^{\alpha\beta} \xi_\alpha\xi_\beta}}\bigg]\bigg\}+ T^{(25)}_{-1}(g_{\alpha\beta})\\
 & = \frac{1}{2\sqrt{\sum_{\alpha,\beta} g^{\alpha\beta} \xi_\alpha\xi_\beta}} \bigg\{\frac{\partial }{\partial x_3} \bigg[ \frac{\sum_{\gamma} \frac{\partial g^{\gamma\beta}}{\partial x_3}\, i \xi_\gamma }{2 \sqrt{\sum_{\alpha,\beta} g^{\alpha\beta}\xi_\alpha\xi_\beta}}\bigg]\bigg\} + T^{(25)}_{-1}(g_{\alpha\beta})\\
& =  \frac{\sum_{\gamma}\frac{\partial^2 g^{\gamma\beta}}{\partial x_3^2} \,i \xi_\gamma }{4\sum_{\alpha,\beta} g^{\alpha\beta} \xi_\alpha\xi_\beta}  + T^{(26)}_{-1}(g_{\alpha\beta})\end{align*}
 and
 \begin{align*} \label{19.8.30-1'}
\phi_{-1}^{3\beta} &= \frac{1}{2\sqrt{\sum_{\alpha,\beta} g^{\alpha\beta}\xi_\alpha\xi_\beta}}   \bigg(
 \frac{\partial \phi_0^{3\beta}}{\partial x_3}- a_{3\beta}+R^{3}_\beta\bigg)  +T_{-1}^{(27)}(g_{\alpha\beta})\\
 &= \frac{1}{2\sqrt{\sum_{\alpha,\beta} g^{\alpha\beta}\xi_\alpha\xi_\beta}}   \bigg(
 \frac{\partial \phi_0^{3\beta}}{\partial x_3}\bigg)  +T_{-1}^{(28)}(g_{\alpha\beta})\\
 & =  \frac{1}{2\sqrt{\sum_{\alpha,\beta} g^{\alpha\beta} \xi_\alpha\xi_\beta}} \bigg\{\frac{\partial }{\partial x_3}\bigg[ \frac{- 2i \sum_{\eta,\gamma} g^{\eta\gamma} \Gamma^3_{\beta \eta} \xi_\gamma}{2\sqrt{\sum_{\alpha,\beta}
g^{\alpha\beta} \xi_\alpha\xi_\beta}} \bigg]\bigg\} + T^{(29)}_{-1}(g_{\alpha\beta})\\
  &=   \frac{1}{2\sqrt{\sum_{\alpha,\beta} g^{\alpha\beta} \xi_\alpha\xi_\beta}} \bigg\{ \frac{\partial }{\partial x_3} \bigg[\frac{- 2i \sum_{\eta,\gamma} g^{\eta\gamma} \big(-\frac{1}{2} \frac{\partial g_{\beta\eta}}{\partial x_3} \big) \xi_\gamma }{2\sqrt{\sum_{\alpha,\beta}
g^{\alpha\beta} \xi_\alpha\xi_\beta}} \bigg]\bigg\}  + T^{(29)}_{-1}(g_{\alpha\beta})\\
& = \frac{1}{2\sqrt{\sum_{\alpha,\beta} g^{\alpha\beta} \xi_\alpha\xi_\beta}}\bigg\{\frac{\partial}{\partial x_3} \bigg( \frac{\sum_{\eta,\gamma} g^{\eta \gamma} \frac{\partial g_{\beta \eta}}{\partial x_3} \,i \xi_\gamma}{2 \sqrt{\sum_{\alpha,\beta} g^{\alpha\beta}\xi_\alpha\xi_\beta}}\bigg)\bigg\} +  T^{(29)}_{-1}(g_{\alpha\beta}),\end{align*}
so that
\begin{eqnarray}  \label{19.9.3;4'}
  - \frac{\sum_{\beta}  \phi_{-1}^{\beta 3} \, i\xi_\beta}{\sqrt{\sum_{\alpha,\beta} g^{\alpha\beta} \xi_\alpha \xi_\beta }}\!\!\!&=\!& \!\!\! \frac{-1}{4\big(\sum_{\alpha,\beta}g^{\alpha\beta}\xi_\alpha\xi_\beta\big)^{3/2}} \bigg(\!
\sum_{\beta,\gamma} \frac{\partial^2 g^{\gamma\beta}}{\partial x_3^2} (i\xi_\gamma)(i\xi_\beta) \bigg)  +   T^{(30)}_{-1}(g_{\alpha\beta})\\
\!\!\!&=\!& \!\!\! \frac{\sum_{\alpha,\beta}
\frac{\partial^2 g^{\alpha\beta}}{\partial x_3^2}\xi_\alpha \xi_\beta}{4\big(\sum_{\alpha,\beta} g^{\alpha\beta}\xi_\alpha\xi_\beta\big)^{3/2}}+ T_{-1}^{(30)} (g_{\alpha\beta})\nonumber\end{eqnarray}
and
\begin{eqnarray}  \label{19.9.3;4-1}\quad\;\;\;\;\;\;\;
 \frac{\sum_{\alpha\!, \beta} g^{\alpha \beta} \phi_{-1}^{3\beta}  i\xi_\alpha}{\sqrt{\sum_{\alpha,\beta} g^{\alpha\beta} \xi_\alpha \xi_\beta }} \!\!\!\!&=\!\!& \!\!\! \frac{1}{2\sum_{\alpha,\beta}g^{\alpha\beta}\xi_\alpha\xi_\beta}\!\left\{ \!\sum\limits_{\alpha,\beta} \!g^{\alpha\beta}\! \bigg[\!\frac{\partial }{\partial x_3}\! \bigg( \!\frac{\sum_{\eta,\gamma} g^{\eta\gamma}\frac{\partial g_{\beta\eta}}{\partial x_3} \,i\xi_\gamma}{2\sqrt{\sum_{\alpha,\beta} g^{\alpha\beta} \xi_\alpha\xi_\beta}}\!\bigg)\! \bigg] i\xi_\alpha \!\right\} \! + \!  T^{(31)}_{\!-1}(g_{\alpha\beta})\\
\!\!\!&=\!& \!\!\! \frac{1}{2\sum_{\alpha,\beta}g^{\alpha\beta}\xi_\alpha\xi_\beta}\!\left\{  \bigg[\frac{\partial }{\partial x_3} \bigg( \frac{\sum_{\alpha,\beta,\eta,\gamma} g^{\alpha\beta} g^{\eta\gamma}\frac{\partial g_{\beta\eta}}{\partial x_3} \,i\xi_\gamma}{2\sqrt{\sum_{\alpha,\beta} g^{\alpha\beta} \xi_\alpha\xi_\beta}}\bigg) \bigg] i\xi_\alpha  \!\right\} \! +  \! T^{(32)}_{-1}(g_{\alpha\beta})\nonumber\\
\!\!\!&=\!& \!\!\! \frac{1}{2\sum_{\alpha,\beta} g^{\alpha\beta}\xi_\alpha\xi_\beta} \bigg[ \frac{\partial}{\partial x_3} \bigg(\frac{-\sum_{\alpha,\gamma}
\frac{\partial g^{\alpha\gamma}}{\partial x_3}\,(i\xi_\alpha)(i \xi_\gamma)}{2\sqrt{\sum_{\alpha,\beta} g^{\alpha\beta} \xi_\alpha\xi_\beta}}\bigg) \bigg] +T_{-1}^{(32)} (g_{\alpha\beta})\nonumber\\
\!\!\!&=\!& \!\!\!\frac{ \sum_{\alpha,\beta}
\frac{\partial^2 g^{\alpha\beta}}{\partial x_3^2}\xi_\alpha \xi_\beta}{4\big(\sum_{\alpha,\beta} g^{\alpha\beta}\xi_\alpha\xi_\beta\big)^{3/2}}+ T_{-1}^{(33)} (g_{\alpha\beta}).\nonumber\end{eqnarray}
 Combining (\ref{19.8.12;1}), (\ref{19.9.3,3}), (\ref{19.9.3;4'}) and (\ref{19.9.3;4-1}) we have
\begin{eqnarray} \label{19.9.6-10.} l_{-1}^{11} +l_{-1}^{22}
 \!\!\! && \!\!\!\!\!\!\!= \frac{1}{2\sqrt{\sum\limits_{\alpha,\beta} g^{\alpha\beta} \xi_\alpha\xi_\beta}}\bigg\{
 \bigg[ \frac{1}{4} \sum\limits_{\alpha,\beta}
 g_{\alpha\beta} \frac{\partial^2 g^{\alpha\beta}}{\partial x_3^2}
 +\frac{7}{4} \frac{\sum_{\alpha,\beta} \frac{\partial^2 g^{\alpha\beta}}{\partial x_3^2}\xi_\alpha\xi_\beta}{\sum_{\alpha,\beta} g^{\alpha\beta}\xi_\alpha\xi_\beta}  \;\bigg]
  \bigg\}
 +\!T^{(34)}_{-1}(g_{\alpha\beta})\nonumber\\
      \!\!\!&&\!\!\!\!\!\!\!= \frac{1}{8 \,\sqrt{\sum\limits_{\alpha,\beta} \!g^{\alpha\beta} \!\xi_\alpha\xi_\beta}} \bigg\{\! \sum\limits_{\alpha,\beta}\! g_{\alpha\beta} \frac{\partial^2 g^{\alpha\!\beta}}{\partial x_3^2}
+7\,\frac{\sum_{\alpha,\beta} \frac{\partial^2 g^{\alpha\beta}}{\partial x_3^2} \xi_\alpha \xi_\beta }{ \sum_{\alpha,\beta} g^{\alpha\beta}\xi_\alpha\xi_\beta} \!\bigg\}
+ T^{(34)}_{-1}(g_{\alpha\!\beta})\nonumber
 .\end{eqnarray}
 Setting  \begin{eqnarray} \label{19.8.6;4} h_2:= \sum\limits_{\alpha,\beta} g_{\alpha \beta} \frac{\partial^2 g^{\alpha\beta}}{\partial x_3^2},\end{eqnarray}
 we thus have  \begin{eqnarray*}  \label{19.8.6;2}  l_{-1}^{11}+l_{-1}^{22} = -\frac{1} {8 \big(\sum\limits_{\alpha,\beta} g^{\alpha\beta} \xi_\alpha\xi_\beta\big)^{3/2}} \bigg[ \sum\limits_{\alpha, \beta} \bigg( h_2g^{\alpha\beta} +7\,\frac{\partial^2 g^{\alpha\beta}}{\partial x_3^2}\bigg) \xi_\alpha\xi_\beta\bigg]+T_{-1}^{(34)}(g_{\alpha\beta}).\end{eqnarray*}
 Note that $g^{\alpha\beta}$ have been known by $\Lambda^{11}$ according to the discussion before. Evaluating the above result for all $\xi'\in T^*(\Gamma)$ shows that
the symbol of degree $-1$ of $L^{11}+L^{22}$ completely determines the quadratic form
  \begin{eqnarray}\label{19.8.6;3}  k_2^{\alpha\beta}= h_2g^{\alpha\beta} +7 \frac{\partial^2 g^{\alpha\beta}}{\partial x_3^2}\end{eqnarray}  along $\partial M$.  From this quadratic form, (\ref{19.8.6;7}) and (\ref{19.8.6;4}) we find that
  \begin{eqnarray*} && \sum\limits_{\alpha,\beta} g_{\alpha\beta} k^{\alpha\beta}_2 = h_2\sum\limits_{\alpha,\beta} g_{\alpha\beta}  g^{\alpha\beta} +7\sum\limits_{\alpha,\beta} g_{\alpha\beta}\, \frac{\partial^2 g^{\alpha\beta}}{\partial x_3^2} \\
  && \qquad \qquad \quad   = 2h_2 +7h_2=9h_2.\end{eqnarray*}
  Therefore (\ref{19.8.6;3}) imply that \begin{eqnarray*} \frac{\partial^2 g^{\alpha\beta} }{\partial x_3^2} =\frac{1}{7}\bigg\{   k^{\alpha\beta}_2-\frac{1}{9} \bigg(\sum\limits_{\gamma,\rho} g_{\gamma\rho} k^{\gamma\rho}_2\bigg)g^{\alpha\beta}\bigg\}\end{eqnarray*}
   is uniquely determined by $l_{-1}^{11}+ l_{-1}^{22}$ and $\psi^{11}_1$, and hence $\frac{\partial^2 g^{\alpha\beta} }{\partial x_3^2}$ are uniquely determined on $\Gamma$ by $\Lambda_{g,\Gamma}$.

Finally, we will consider the symbol $l^{11}_{-m-1}+l^{22}_{-m-1}$ of degree $-m-1$ for $L^{11}+L^{22}$ for the general integer $m\ge 1$.
From (\ref{19.7.21-7}), we see that    \begin{eqnarray*} l^{jj}_{-m-1} = \phi_{-m-1}^{jj} - \phi_{-m-1}^{j3} q_{-1}\sigma \,i\xi_j  + \sum\limits_{\alpha} g^{\alpha j} \bigg(q_{-1}\sigma \phi^{3j}_{-m-1}   +  q_{-m-3}\sigma \, i\xi_j \bigg) i\xi_\alpha+T_{-m-1}^{(1)}(g_{\alpha\beta}), \quad \mbox{for}\,\, j=1,2,\end{eqnarray*}
 where $T_{-m-1}^{(s)}(g_{\alpha\beta})$, $(s=1,2,\cdots)$, are expressions involving only the boundary values of $g^{\alpha\beta}$, $g_{\alpha\beta}$, and their normal derivatives of order at most $m+1$.
From (\ref{19.8.12-5}) we have
\begin{eqnarray*} q_{-m-3} = -\frac{ 1}{\sigma\phi_1^{33}} \,q_{-1} \,\sigma \phi^{33}_{-m-1} + T^{(2)}_{-m-1}(g_{\alpha\beta}).\end{eqnarray*}
  It follows that \begin{eqnarray} \label{19.9.4-21} && \;\;\,\;  l^{11}_{-m-1} +l^{22}_{-m-1}= \phi_{-m-1}^{11} +\phi_{-m-1}^{22} +\phi_{-m-1}^{33} -q_{-1} \sum\limits_{j=1}^2 \phi_{-m-1}^{j3} \,\sigma\,i\xi_j \\
   &&  \qquad \quad +  \sum\limits_{j,\alpha=1}^2  g^{\alpha j} q_{-1}\sigma \phi_{-m-1}^{3j} \,i\xi_\alpha
   +T_{-m-1}^{(3)} (g_{\alpha\beta})\nonumber\\
   &&\quad\;\,\;  \; =\phi_{-m-1}^{11}  +\phi_{-m-1}^{22} +\phi_{-m-1}^{33} - \frac{\sum_{\beta} \phi_{-m-1}^{\beta 3} i\xi_\beta}{\sqrt{\sum_{\alpha,\beta} g^{\alpha\beta} \xi_\alpha\xi_\beta}}
    +\frac{\sum_{\alpha, \beta} g^{\alpha\beta} \phi_{-m-1}^{3\beta} i\xi_\alpha}{ \sqrt{ \sum_{\alpha,\beta} g^{\alpha\beta}\xi_\alpha\xi_\beta}}
    + T_{-m-1}^{(3)} (g_{\alpha\beta}).\nonumber\end{eqnarray}
 From (\ref{19.8.3-1}) we have
 \begin{eqnarray} \label{19.9.4-1,} \phi_{-m-1}^{jk} = \frac{1}{2\sqrt{\sum_{\alpha,\beta} g^{\alpha\beta} \xi_\alpha\xi_\beta}} \frac{\partial \phi_{-m}^{jk}}{\partial x_3} +T_{-m-1}^{(4)} (g_{\alpha\beta})\quad \,\mbox{for}\;\, j,k=1,2,3.\end{eqnarray}
 We will give the corresponding estimates for $\sum_{k=1}^3 \phi_{-m-1}^{kk}$, $ \frac{\sum_{\alpha, \beta} g^{\alpha\beta} \phi_{-m-1}^{3\beta} i\xi_\alpha}{ \sqrt{ \sum_{\alpha,\beta} g^{\alpha\beta}\xi_\alpha\xi_\beta}}$ and $ - \frac{\sum_{\beta} \phi_{-m-1}^{\beta 3} i\xi_\beta}{\sqrt{\sum_{\alpha,\beta} g^{\alpha\beta} \xi_\alpha\xi_\beta}}$ by induction.
 We first show that for all $m\ge 1$,
 \begin{eqnarray} \label{2020.3.29-5} \quad\;\,\; \;\quad\;\;\sum_{k=1}^3\! \phi_{-m}^{kk}\!=\! \bigg(\!2\sqrt{\sum\limits_{\alpha,\beta}\! g^{\alpha\beta} \xi_\alpha\xi_\beta}\!\bigg)^{\!\!-m}\! \bigg\{\!
\frac{1}{4}\! \sum\limits_{\alpha,\beta}\! g^{\alpha\beta} \frac{\partial^{\!m+1} g^{\alpha\beta}}{\partial x_3^{m+1}}\!+\!\frac{3}{4}\frac{\sum_{\alpha,\beta} \frac{\partial^{m+1} g^{\alpha\beta} }{\partial x_3^{m+1}}\xi_\alpha\xi_\beta}{ \sum_{\alpha,\beta}
g^{\alpha\beta} \xi_\alpha\xi_\beta}\!\bigg\} \!\!+\!\!  T_{\!-\!m\!-\!1}^{(5)} \!(g_{\alpha\beta}).\end{eqnarray}
  Suppose we have shown that for $1\le r\le m$,
\begin{eqnarray} \label{19.9.4-6} \quad \quad\;\;\;\,\;\;\sum_{k=1}^3\! \phi_{-r}^{kk}\!=\! \bigg(\!2\sqrt{\sum\limits_{\alpha,\beta}\! g^{\alpha\beta} \xi_\alpha\xi_\beta}\,\bigg)^{\!\!-r} \!\bigg\{\!
\frac{1}{4}\! \sum\limits_{\alpha,\beta}\! g^{\alpha\beta} \frac{\partial^{r+1} g^{\alpha\beta}}{\partial x_3^{r+1}}\!+\!\frac{3}{4}\frac{\sum_{\alpha,\beta} \frac{\partial^{r+1} g^{\alpha\beta} }{\partial x_3^{r+1}}\xi_\alpha\xi_\beta}{ \sum_{\alpha,\beta}
g^{\alpha\beta} \xi_\alpha\xi_\beta}\!\bigg\} \!+\!  T_{\!-\!r\!-\!1}^{(6)} (g_{\alpha\beta}).\end{eqnarray}
Clearly, this estimates actually holds when $r=1$ by (\ref{19.9.3,3}).
 Then, from  (\ref{19.9.4-1,}) and  (\ref{19.9.4-6}), we have  \begin{eqnarray} \label{19.9.4-9} &&\; \;\quad \;\,\;\;\;\sum_{k=1}^3 \phi_{-m-1}^{kk}
=\frac{1}{2\sqrt{\sum_{\alpha,\beta} g^{\alpha\beta} \xi_\alpha\xi_\beta}} \frac{\partial}{\partial x_3} \big(\sum_{k=1}^3 \phi_{-m}^{kk} \big) +T_{-m-1}^{(7)} (g_{\alpha\beta})\\
&& \!\!=\!\! \frac{1}{2\sqrt{\sum_{\alpha,\beta} \!g^{\alpha\beta} \xi_\alpha\xi_\beta}} \frac{\partial}{\partial x_3}\!\bigg[
\!\bigg(\!2\sqrt{\sum\limits_{\alpha,\beta}\! g^{\alpha\beta} \xi_\alpha\xi_\beta}\!\bigg)^{\!\!-m} \!\bigg(
\!\frac{1}{4}\! \sum\limits_{\alpha,\beta}\! g^{\alpha\beta} \frac{\partial^{m\!+\!1} g^{\alpha\beta}}{\partial x_3^{m+1}}\!+\!\frac{3}{4}\!\frac{\sum_{\alpha,\beta} \!\frac{\partial^{m\!+\!1} g^{\alpha\beta} }{\partial x_3^{m+1}}\xi_\alpha\xi_\beta}{ \sum_{\alpha,\beta}
\!g^{\alpha\beta} \xi_\alpha\xi_\beta}\!\bigg)\!\bigg]\!\!+\! \! T_{\!-\!m\!-\!1}^{(8)} (g_{\alpha\beta})\nonumber\\
&&= \bigg(2\sqrt{\sum\limits_{\alpha,\beta}\! g^{\alpha\beta} \xi_\alpha\xi_\beta}\,\bigg)^{-m-1} \frac{\partial}{\partial x_3}
 \bigg(
\!\frac{1}{4}\! \sum\limits_{\alpha,\beta}\! g^{\alpha\beta} \frac{\partial^{m+1} g^{\alpha\beta}}{\partial x_3^{m+1}}+\frac{3}{4}\frac{\sum_{\alpha,\beta} \frac{\partial^{m+1} g^{\alpha\beta} }{\partial x_3^{m+1}}\xi_\alpha\xi_\beta}{ \sum_{\alpha,\beta}
g^{\alpha\beta} \xi_\alpha\xi_\beta}\bigg)+ T_{\!-\!m\!-\!1}^{(9)} (g_{\alpha\beta})\nonumber\\
&& =  \bigg(2\sqrt{\sum\limits_{\alpha,\beta}\! g^{\alpha\beta} \xi_\alpha\xi_\beta}\,\bigg)^{-m-1}
 \bigg(
\!\frac{1}{4}\! \sum\limits_{\alpha,\beta}\! g^{\alpha\beta} \frac{\partial^{m+2} g^{\alpha\beta}}{\partial x_3^{m+2}}+\frac{3}{4}\frac{\sum_{\alpha,\beta} \frac{\partial^{m+2} g^{\alpha\beta} }{\partial x_3^{m+2}}\xi_\alpha\xi_\beta}{ \sum_{\alpha,\beta}
g^{\alpha\beta} \xi_\alpha\xi_\beta}\bigg)+ T_{\!-\!m\!-\!1}^{(10)} (g_{\alpha\beta}). \nonumber\end{eqnarray}
Thus, by induction we get that (\ref{2020.3.29-5}) holds for all $m\ge 1$.
 Next, we will prove by induction that
 \begin{eqnarray} \label{19.9.4-10} \;\;\;\;\frac{\sum\limits_{\alpha,\beta} g^{\alpha\beta}\phi_{-m-1}^{3\beta} i\xi_\alpha}{\sqrt{\sum_{\alpha,\beta}g^{\alpha\beta} \xi_\alpha\xi_\beta}} \!=\! \frac{1}{2} \bigg(2\sqrt{\sum\limits_{\alpha,\beta}\! g^{\alpha\beta} \xi_\alpha\xi_\beta}\,\bigg)^{\!\!-\!m-1}\bigg(\frac{\sum_{\alpha,\beta} \frac{\partial^{m+2} g^{\alpha\beta}}{\partial x_3^{m+2}}\xi_\alpha\xi_\beta}{\sum_{\alpha,\beta}g^{\alpha\beta} \xi_\alpha\xi_\beta}\bigg)\! +\!  T_{\!-m-1}^{(11)} (g_{\alpha\beta}).\end{eqnarray}
 Suppose we have also shown that, when $1\le r\le m$,
 \begin{eqnarray}  \label{19.9.4-11}  \;\; \frac{\sum\limits_{\alpha,\beta} g^{\alpha\beta}\phi_{-r}^{3\beta} i\xi_\alpha}{\sqrt{\sum_{\alpha,\beta}g^{\alpha\beta} \xi_\alpha\xi_\beta}} \!=\! \frac{1}{2} \bigg(2\sqrt{\sum\limits_{\alpha,\beta}\! g^{\alpha\beta} \xi_\alpha\xi_\beta}\,\bigg)^{\!\!-r}\bigg( \frac{\sum_{\alpha,\beta} \frac{\partial^{r+1} g^{\alpha\beta}}{\partial x_3^{r+1}}\xi_\alpha\xi_\beta}{\sum_{\alpha,\beta} g^{\alpha\beta} \xi_\alpha\xi_\beta} \bigg)\!+\!  T_{\!-m-1}^{(12)} (g_{\alpha\beta}).\end{eqnarray}
Again, this actually is true when $r=1$ by (\ref{19.9.3;4-1}).
 By (\ref{19.9.4-1,})  we find that
\begin{eqnarray*} \label{19.9.4-1:,} \phi_{-m-1}^{3\beta} = \frac{1}{2\sqrt{\sum_{\alpha,\beta} g^{\alpha\beta} \xi_\alpha\xi_\beta}} \frac{\partial \phi_{-m}^{3\beta}}{\partial x_3} +T_{-m-1}^{(13)} (g_{\alpha\beta})\quad \,\mbox{for}\;\, j,k=1,2,3.\end{eqnarray*}
It follows from this and (\ref{19.9.4-11}) that
\begin{eqnarray} \label{19.9.4-17}  &&  \;\;\;\;\;\frac{\sum\limits_{\alpha,\beta}g^{\alpha\beta} \phi_{-m-1}^{3\beta} i\xi_\alpha}{\sqrt{\sum_{\alpha,\beta}g^{\alpha\beta} \xi_\alpha\xi_\beta}} =
  \frac{1}{2\sum_{\alpha,\beta} g^{\alpha\beta} \xi_\alpha\xi_\beta} \sum\limits_{\alpha,\beta}g^{\alpha\beta}\frac{\partial \phi_{-m}^{3\beta}}{\partial x_3}\,i\xi_\alpha +T_{-m-1}^{(14)} (g_{\alpha\beta})\\
  && \quad  =  \frac{1}{2\sqrt{\sum_{\alpha,\beta} g^{\alpha\beta} \xi_\alpha\xi_\beta}}\; \frac{\partial }{\partial x_3} \bigg( \frac{\sum\limits_{\alpha,\beta}
  g^{\alpha\beta} \phi_{-m}^{3\beta}i \xi_\alpha}{\sqrt{\sum_{\alpha,\beta}g^{\alpha\beta}\xi_\alpha\xi_\beta}}\bigg) + T_{-m-1}^{(15)} (g_{\alpha\beta})\nonumber\\
  && \quad =  \frac{1}{2\sqrt{\sum_{\alpha,\beta} g^{\alpha\beta} \xi_\alpha\xi_\beta}}\;\frac{\partial }{\partial x_3} \bigg[\frac{1}{2}\bigg( 2 \sqrt{ \sum\limits_{\alpha,\beta}
  g^{\alpha\beta}\xi_\alpha\xi_\beta}\,\bigg)^{\!-m}\bigg( \frac{\sum_{\alpha,\beta} \frac{\partial^{m+1} g^{\alpha\beta}}{\partial x_3^{m+1}}\xi_\alpha\xi_\beta}{\sum_{\alpha,\beta} g^{\alpha\beta} \xi_\alpha\xi_\beta}\bigg)\bigg]  + T_{-m-1}^{(16)} (g_{\alpha\beta})\nonumber\\
&& \quad = \frac{1}{2}\bigg( 2 \sqrt{ \sum\limits_{\alpha,\beta}
  g^{\alpha\beta}\xi_\alpha\xi_\beta}\,\bigg)^{\!-m-1}\bigg( \frac{\sum_{\alpha,\beta} \frac{\partial^{m+2} g^{\alpha\beta}}{\partial x_3^{m+2}}\xi_\alpha\xi_\beta}{\sum_{\alpha,\beta} g^{\alpha\beta} \xi_\alpha\xi_\beta}\bigg)  + T_{-m-1}^{(17)} (g_{\alpha\beta}).\nonumber\end{eqnarray}
  Thus (\ref{19.9.4-10}) holds for all $m\ge 1$ by induction.
  Similarly, by applying induction we can easily prove that for all $m\ge 1$,
  \begin{eqnarray} \label{19.9.4-18}\;\; \;\;\;-\frac{\sum_{\beta} \phi_{-m-1}^{\beta 3}i\xi_\beta}{\sum_{\alpha,\beta} g^{\alpha\beta} \xi_\alpha\xi_\beta} = \frac{1}{2}\bigg( \!2 \sqrt{ \sum\limits_{\alpha,\beta}
  g^{\alpha\beta}\xi_\alpha\xi_\beta}\bigg)^{\!-m-1} \bigg(\frac{\sum_{\alpha,\beta} \frac{\partial^{m+2} g^{\alpha\beta}}{\partial x_3^{m+2}}\xi_\alpha\xi_\beta}{\sum_{\alpha,\beta} g^{\alpha\beta} \xi_\alpha\xi_\beta}\bigg) \! +\! T_{\!-m-1}^{(18)} (g_{\alpha\beta}).\end{eqnarray}
Hence \begin{eqnarray} \label{19.9.4-19}&&-  \frac{ \sum\limits_{\beta}  \phi_{-m-1}^{\beta 3} i\xi_\beta}{\sqrt{\sum_{\alpha,\beta} g^{\alpha\beta} \xi_\alpha\xi_\beta}} + \frac{ \sum\limits_{\alpha,\beta}g^{\alpha\beta}  \phi_{-m-1}^{3\beta} i\xi_\alpha}{\sqrt{\sum_{\alpha,\beta} g^{\alpha\beta} \xi_\alpha\xi_\beta}}+ T_{-m-1}^{(19)} (g_{\alpha\beta})\\
&& \quad \;\quad \;= \bigg(2\sqrt{\sum\limits_{\alpha,\beta}\! g^{\alpha\beta} \xi_\alpha\xi_\beta}\,\bigg)^{-m-1}
 \bigg(\frac{\sum_{\alpha,\beta} \frac{\partial^{m+2} g^{\alpha\beta} }{\partial x_3^{m+2}}\xi_\alpha\xi_\beta}{\sum_{\alpha,\beta} g^{\alpha\beta} \xi_\alpha\xi_\beta}\bigg)+ T_{\!-\!m\!-\!1}^{(20)} (g_{\alpha\beta}).
\nonumber\end{eqnarray}
Combining  (\ref{19.9.4-21}), (\ref{19.9.4-9}) and (\ref{19.9.4-19}), we get
\begin{eqnarray} \label{19.9.4-23} \;\quad \quad &\;\,\;\;l_{-m-1}^{11}\! +\! l_{-m-1}^{22} \!=\!\bigg(\!2\sqrt{\sum\limits_{\alpha,\beta}\! g^{\alpha\beta} \xi_\alpha\xi_\beta}\!\bigg)^{\!\!-\!m\!-\!1}
 \!\bigg(\!
\frac{1}{4}\! \sum\limits_{\alpha,\beta}g^{\alpha\beta}\frac{\partial^{m+2} g^{\alpha\beta}}{\partial x_3^{m+2}} \!+\!\frac{7}{4}  \frac{\sum_{\alpha,\beta} \! \frac{\partial^{m+2} g^{\alpha\beta} }{\partial x_3^{m+2}}\xi_\alpha\xi_\beta}{\sum_{\alpha,\beta} g^{\alpha\beta} \xi_\alpha\xi_\beta}\! \bigg)\! \!+\! T_{\!-\!m\!-\!1}^{(21)} (g_{\alpha\beta})\quad\\
&\qquad \quad\;\; \,\qquad =\!\frac{1}{4}\bigg(\!2\sqrt{\sum\limits_{\alpha,\beta}\! g^{\alpha\beta} \xi_\alpha\xi_\beta}\,\bigg)^{\!\!-\!m\!-\!1}
 \!\bigg(\!
 \sum\limits_{\alpha,\beta}g^{\alpha\beta}\frac{\partial^{m+2} g^{\alpha\beta}}{\partial x_3^{m+2}} \!+\!7  \,\frac{\sum_{\alpha,\beta} \! \frac{\partial^{m+2} g^{\alpha\beta} }{\partial x_3^{m+2}}\xi_\alpha\xi_\beta}{\sum_{\alpha,\beta} g^{\alpha\beta} \xi_\alpha\xi_\beta} \! \bigg)\! \!+\! T_{\!-\!m\!-\!1}^{(21)} (g_{\alpha\beta}).\nonumber\end{eqnarray}
Putting \begin{eqnarray*}
 h_{m+2} =\sum\limits_{\alpha,\beta} g^{\alpha\beta} \frac{\partial^{m+2} g^{\alpha\beta}}{\partial x_3^{m+2}},\quad
\;  \kappa_{m+2}^{\alpha\beta} = h_{{}_{m+2}} \,g^{\alpha\beta} +7\, \frac{\partial^{m+2}g^{\alpha\beta}}{\partial x_3^{m+2}}.\end{eqnarray*}
 we see  that \begin{eqnarray} l_{-m-1}^{11}+l_{-m-1}^{22}= \frac{1}{\big(2\sqrt{\sum\limits_{\alpha,\beta}\! g^{\alpha\beta} \xi_\alpha\xi_\beta}\;\big)^{m+3}}
\sum\limits_{\alpha,\beta} \bigg( h_{{}_{m+2}}\,g^{\alpha\beta} +7\,\frac{\partial^{m+2} g^{\alpha\beta}}{\partial x_3^{m+2}}\bigg) \xi_\alpha\xi_\beta,\end{eqnarray}
i.e., the quadratic form $\kappa_{m+2}^{\alpha\beta}$ is completely determined by $l^{11}_{-m-1}+l_{-m-1}^{22}$ and  $\psi_1^{11}$ (i.e., $g^{\alpha\beta}$).
 Since \begin{eqnarray*} \sum\limits_{\alpha,\beta} g_{{}_{\alpha\beta}} \kappa_{m+2}^{\alpha\beta} =h_{m+2}\; \sum\limits_{\alpha,\beta}
 g_{{}_{\alpha\beta}} g^{\alpha\beta} +7 \sum\limits_{\alpha,\beta} g_{\alpha\beta} \frac{\partial^{m+2} g^{\alpha\beta}}{\partial x_3^{m+2}}=9h_{m+2},\end{eqnarray*}
 we have \begin{eqnarray*} h_{m+2}= \frac{1}{9} \sum\limits_{\alpha,\beta} g_{{}_{\alpha\beta}} \kappa^{\alpha\beta}_{m+2},\end{eqnarray*}
 so that \begin{eqnarray*} \frac{\partial^{m+2}g^{\alpha\beta}}{\partial x_3^{m+2}}=\frac{1}{7} \bigg( \kappa_{m+2}^{\alpha\beta} -\frac{1}{9}
 \big(\sum\limits_{\gamma,\rho} g_{\gamma\rho}\kappa_{m+2}^{\gamma\rho}\big)g^{\alpha\beta}\bigg).\end{eqnarray*}
 This implies that $\frac{\partial^{m+2} g^{\alpha\beta}}{\partial x_3^{m+2}}$ is uniquely determined by $\psi_1^{11}$, $l_0^{11}+l_0^{22}$, $\cdots$, $l^{11}_{-m-1}+l_{-m-1}^{22}$ for any fixed $m\ge 1$.
 Consequently, the operator $\Lambda_{g,\Gamma}$ uniquely determines $g^{\alpha \beta}$ and all order normal derivatives of $g^{\alpha\beta}$ along on $\Gamma$, and the proof of Proposition 3.1 is completed.  $\qquad \qquad \qquad \square$

\vskip 0.34 true cm

 \noindent{\bf Lemma 3.2.} \ {\it Let $\bar M$ be a three-dimensional Riemannian manifold having compact closure and $C^1$-smooth boundary $\partial M$,
  and let $\Gamma$ be a real analytic piece of boundary $\partial M$. Assume that $g$ and $\tilde{g}$ be real-analytic metrics in $M$ up to $\Gamma$. Assume also that the electromagnetic parameters $\mu$ and $\sigma$ are real analytic in $M$ up to $\Gamma$, and that $\mu$ and} $\,\mbox{Re}(\sigma)$ {\it are strictly positive on $\bar M$. If $\Lambda_{g,\Gamma}=\Lambda_{\tilde g,\Gamma}$, then there exists a neighborhood $\mathfrak{U}$ of $\,\Gamma$ in $\bar M$ and a real-analytic map $\varrho_0: \mathfrak{U}\to \bar M$ such that $\varrho_0\big|_{\Gamma}=\mbox{identity}$ and $g=\varrho_0^* \tilde {g}$.}

 \vskip 0.26 true cm

    \noindent  {\it Proof.} \  Let $\mathfrak{V}$ be some connected open set in the half-space $\{ x_3\ge 0\}\subset {\mathbb{R}}^3$ containing the origin. For any $x_0\in \Gamma$ we define a real-analytic diffeomorphism $\varsigma_{x_0}: \mathfrak{V}\to {\mathfrak{U}}_{x_0}$ (respectively, ${\tilde\varsigma}_{x_0}: \mathfrak{V}\to {\tilde{\mathfrak{U}}}_{x_0}$), where  $( x_1, x_2, x_{3})$ (respectively, $( {\tilde x}_1, {\tilde x}_{2}, {\tilde x}_{3})$) denote the corresponding boundary normal coordinates for $g$ (respectively, $\tilde g$) defined on a connected neighborhood ${\mathfrak{U}}_{x_0}$ (respectively, ${\tilde{\mathfrak{U}}}_{x_0}$) of $x_0$ (see \cite{LU}).     It follows from Proposition 3.1 that the metric $\varsigma_{x_0}^* g$ and ${\tilde{\varsigma}}_{x_0}^* \tilde{g}$ are real-analytic metric on $\mathfrak{V}$, where Taylor series at origin by explicit formulas involving the symbol of $\Xi_g$ in $\{x_\alpha\}$ coordinates. Clearly, these two metric must be identical on $\mathfrak{V}$. Set $ \varrho_{x_0} ={\tilde{\varsigma}}_{x_0} \circ \varsigma_{x_0}^{-1} :{\mathfrak{U}}_{x_0}\to {\tilde{\mathfrak{U}}}_{x_0}$.  Then, $\varrho_{x_0}$
is a real-analytic diffeomorphism which fixes the portion of $\Gamma$ lying in ${\mathfrak{U}}_{x_0}$ and satisfies ${\varrho}^*_{x_0}\tilde {g}= g$.
Therefore, repeating this construction for each $x_0\in \Gamma$, we obtain a real analytic map $\varrho_0: \mathfrak{U}\to \bar M$ such that $\varrho_0\big|_{\Gamma}=\mbox{identity}$ and $\varrho_0^*\tilde{g} =g$. \qed

\vskip 1.48 true cm

\section{Determining metric tensor from electromagnetic Dirichlet-to-Neumann map on a piece of analytic surface}

\vskip 0.48 true cm

Recall that $\Gamma$ is a real analytic piece of boundary, and that the metric tensor of
$M$ is real analytic up to $\Gamma$. Let $f \in C(\partial M)\cap TH^{\frac{1}{2}}(\partial M)$ with supp $f \subset \Gamma$, and let $(E,H)\in (L^2(M))^3\times \big(L^2(M))^3$ be the unique solution of the Maxwell's equations:
 \begin{eqnarray}\label{2020.1.31-2}  \left\{\!\! \begin{array} {ll}
  \mbox{curl}\, E= i \omega \mu H \quad &\mbox{in}\,\; M,\\
 \mbox{curl}\, H=-i \omega \sigma E \quad &\mbox{in}\,\; M,\\
     \nu \times E = f  \quad &\mbox{on}\,\; \partial M.\end{array}\right. \end{eqnarray}
We assume that we know the all boundary data $(\mu\times E, \nu\times H)$ on $\Gamma$ of all possible solutions of (\ref{2020.1.31-2}), for $f$ supported on $\Gamma$. Equivalently, we assume that we know the electromagnetic Dirichlet-to-Neumann map
\begin{eqnarray} \Lambda_{g,\Gamma} :f \rightarrow \nu\times H\big|_{\Gamma} \end{eqnarray}

\vskip 0.2 true cm

\vskip 0.12 true cm

Let $\mathcal{O}$ be a boundary normal coordinate collar neighborhood of $\Gamma$ (i.e., $\mathcal{O}=\Gamma \times (-r,r)$ in the boundary normal coordinates). Set $U=\mathcal{O}\setminus \bar M$. Then, we define a manifold $\tilde{M}$ by gluing
to $M\cup \Gamma$ the boundary normal coordinate collar $U$ (i.e., U=$\partial M \times (-r,0)$) with metric described as follows.
In section 3, we have shown that for a connected Riemannian manifold with compact closure, the $\Lambda_{g, \Gamma}$ determines $g_{jk}$ and all order normal derivatives $\frac{\partial^k g_{jk}}{\partial x_\nu^k}(x)$, $k\ge 0$ of the metric tensor $g$ on
  $\Gamma$ (This result is based on the local fact that when
$\Lambda_{g, \Gamma}$ is considered as a pseudodifferential operator, its full symbol determines
the all order derivatives of the metric in boundary normal coordinates). By Taylor's series we can continue the metric so that the new metric is real-analytic in $\tilde{M}$ when $r$ is small enough (In other words, the real analytic metric $g_{jk}$ is uniquely determined by $\Lambda_{g,\Gamma}$ in $\tilde{M}$, see the proof of Lemma 3.2). We denote the new metric of $\tilde{M}$ also by $g_{jk}$.

\vskip 0.38 true cm

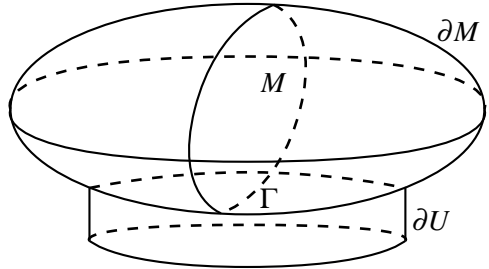
\begin{figure}[h]
\centering
\begin{tikzpicture}[scale=0.7,line width=1]

\draw [thick] (0,0) ellipse (4.5 and 2);
\draw [thick] plot[smooth, tension=2] coordinates {(-4.5,0) (0,-1) (4.5,0)};
\draw [dash pattern=on 4pt off 4pt] plot[smooth, tension=2] coordinates {(-4.5,0) (0,1) (4.5,0)};
\draw [thick] plot[smooth, tension=1.5] coordinates {(0.5,1.98) (-1,0) (-0.5,-1.98)};
\draw [dash pattern=on 4pt off 4pt] plot[smooth, tension=1.5] coordinates {(0.5,1.98) (1,0) (-0.5,-1.98)};
\draw [dash pattern=on 4pt off 4pt] plot[smooth, tension=1] coordinates {(-3,-1.5) (0,-1.2) (3,-1.5)};
\draw [thick] (-3,-1.5) -- (-3,-2.5);
\draw [thick] (3,-1.5) -- (3,-2.5);
\draw [dash pattern=on 4pt off 4pt] plot[smooth, tension=2] coordinates {(-3,-2.5) (0,-2.2) (3,-2.5)};
\draw [thick] plot[smooth, tension=1.5] coordinates {(-3,-2.5) (0,-3) (3,-2.5)};
\node at (0.5,0.5) {$M$};
\node at (4,1.5) {$\partial M$};
\node at (3.5,-2.1) {$\partial U$};
\node at (0.4,-1.7) {$\Gamma$};
\end{tikzpicture}

\caption{\label{figure1}The manifold $M$, the set $\Gamma \subset \partial M$ and the extension of the manifold over $\Gamma$.}

\end{figure}

\vskip 0.46 true cm

From (\ref{2020.4.11-4}) and (\ref{19.9.6-1}), we see that if $(E,H)$ is a solution of Maxwell's equations (\ref{2020.1.31-2}) then $E$ satisfies
\begin{eqnarray} \label{20.4.4-2}  \left\{\!\!\begin{array}{ll} \mbox{curl}\; \mbox{curl}\, E
-  \big(\mbox{grad}\, (\log \mu)\big) \times \mbox{curl}\; E -\omega^2 \mu \sigma E=0  \;\;&\mbox{in}\;\; M,\\
\mbox{div}\, (\sigma E)=0 \;\; &\mbox{in}\;\; M,\\
\nu\times E=f \;\; &\mbox{on} \;\; \partial M.\end{array}\right.\end{eqnarray}
 Conversely, let $E$ be a solution of (\ref{20.4.4-2}) then $E$ and $H:=\frac{1}{i\omega \mu} \mbox{curl}\, E$ satisfies Maxwell's equations (\ref{2020.1.31-2}).
  In the local (or local boundary normal) coordinates, (\ref{20.4.4-2}) can be written as
\begin{eqnarray}\label{20.4.11-8}  \left\{\!\begin{array}{ll}  {\mathcal{M}}_g E=0, \;\;&\mbox{in}\;\; M,\\
\mbox{div}\, (\sigma E)=0 \;\; &\mbox{in}\;\; M,\\
\nu\times E=f \;\; &\mbox{on} \;\; \partial M.\end{array}\right.\end{eqnarray}
where  ${\mathcal{M}}_gE:=\mbox{curl}\; \mbox{curl}\, E
-  \big(\mbox{grad}\, (\log \mu)\big) \times \mbox{curl}\; E -\omega^2 \mu \sigma E$ is as in (\ref{19.11.26-1}) (or (\ref{19.7.7-2,1})).
Note that ${\mathcal{M}}_g$  is a system of second-order linear elliptic equations in $M$. We define (electric) dyadic Green's function on ${\tilde M}$ which satisfies the equation \begin{eqnarray} \label{2020.2.1-5} \left\{\!\begin{array}{ll}  {\mathcal{M}}_g {G}^e (\cdot,y) = \delta_y \quad
\mbox{in} \,\; \tilde{M},\\
\mbox{div}\, G^e(\cdot,y) =0 \quad \;\; \mbox{in} \,\; \tilde{M}\setminus\{y\},\\
G^e (\cdot, y)\big|_{\partial \tilde{M}}=0. \end{array}\right.\end{eqnarray}
 It is well-known that the (electric) dyadic Green function $G^e(x, y)$ is a real-analytic function of $x$
when $x \notin \{y\} \cup \partial {\tilde{M} } $.
Moreover, when $x$ is near to a given $y$
it has the asymptotics (see \cite{Chew})
\begin{eqnarray} G^e(x, y)\to \infty  \quad \mbox{as}\;\, x\to y.\end{eqnarray}

 The (electric) dyadic Green function $G^e(x,y)$ is a generalization of the following classical (electric) dyadic Green function for Maxwell's equations:

 \vskip 0.28 true cm

\noindent{\bf Example 4.1} (see p.$\,31$ of \cite{Chew}) \
{\it For the classical Maxwell's equations (where $g_{jk}=\delta_{jk}$, and $\mu$ and  $\sigma$ are constants),  the classical dyadic Green's function (for electric field) satisfies the equation \begin{eqnarray}\label{2020.2.3-11}  \mbox{curl}\; \mbox{curl}\; G^e (x, y) - \omega^2\mu\sigma  G^e(x, y) =\delta(x - y)\,I. \end{eqnarray}
Then, after post-multiplying (\ref{2020.2.3-11}) by $G^e(x, y)$, pre-multiplying (\ref{2020.2.3-11}) by $E(x)$, subtracting the resultant equations and integrating the difference over $M$, we have \begin{eqnarray*} E(y) =\int_{M}  \big[ E(x) \cdot \big(\mbox{curl}\;\mbox{curl}\; G^e(x,y)\big) +\big(\mbox{curl}\;\mbox{curl}\; E(x)\big) \cdot G^e(x,y)\big] dV_x.\end{eqnarray*}
Next, using the vector identity that \begin{eqnarray*}  -\nabla\cdot \!\left[E(x) \cdot \big(\mbox{curl}\;\mbox{curl}\; G^e(x,y)\big) \!+\!\big(\mbox{curl}\;\mbox{curl}\; E(x)\big) \cdot G^e(x,y)\right]\!=\! E(x) \cdot  \mbox{curl}\; \mbox{curl}\; G^e(x,y) \!-\! \mbox{curl}\; \mbox{curl}\;E (x) \cdot  G^e(x,y),\end{eqnarray*} and applying Gauss' divergence theorem, one has
 \begin{eqnarray} \label{2020.2.3-13} E(y) = -\int_{\partial M}  \big[\nu\times E(x) \cdot \mbox{curl}\;  G^e(x, y) +i\omega \mu \nu\times H (x) \cdot  G^e(x,y)\big]dS.\end{eqnarray}
 Again, notice that (\ref{2020.2.3-13}) is derived via the use of (\ref{2020.2.3-11}), but no boundary condition has yet been imposed on $G^e(x, y)$ on $\partial M$. Now, if we require that $\nu\times  G^e(x, y) = 0$ for $x\in \partial M$, then (\ref{2020.2.3-13}) becomes \begin{eqnarray} \label{2020.2.3-14} E(y) = -\int_{\partial M}  \big[\nu\times E(x) \cdot \mbox{curl}\;  G^e(x, y)\big] dS.\end{eqnarray}
Note that if no boundary condition is imposed on $G^e(x,y)$, then \begin{eqnarray*} \left[G^e (x,y)\right]^T=G^e (y,x) = \left[I+\frac{\nabla_y\nabla_y}{\omega^2\mu\sigma}\right] g(y-x),  \quad g(y-x)=\frac{e^{i\omega \sqrt{\mu \sigma}|y-x|}}{4\pi |y-x|}. \end{eqnarray*}}

  \vskip 0.22 true cm

Let us next consider two manifolds $M_1$ and $M_2$  for which we have identified
$\Gamma_1 = \Gamma_2 = \Gamma$ and  $\Lambda_{g_1, \Gamma}
=\Lambda_{g_2, \Gamma} $.  Using the previous construction of the
set $U$ and the metric tensor on $U$, which is the same for both manifolds, we
can attach this set and the metric on it to both manifolds, i.e.,
\begin{eqnarray*}  {\tilde{M}}_1= M_1 \cup U, \;  {\tilde{M}}_2= M_2 \cup U.\end{eqnarray*}
Then the  corresponding (electric) dyadic Green functions of ${\tilde{M}}_j$, ($j=1,2$), satisfying
\begin{eqnarray}  \label{2020.2.1-9} \left\{ \!\begin{array} {ll}  {\mathcal{M}}_g {G}^e_j (\cdot,y) = \delta_y \quad
\mbox{in} \,\; {\tilde{M}}_j,\\
\mbox{div}\, G^e_j(\cdot, y)=0 \quad \;\;\mbox{in}\;\, {\tilde{M}}_j\setminus \{y\},\\
G^e_j (\cdot, y)\big|_{\partial {\tilde{M}}_j}=0.\end{array} \right.\end{eqnarray}

\vskip 0.46 true cm

\noindent  {\it Proof of Theorem 1.2.} \  First, we will show that the (electric) dyadic Green's function $G_j^e(x,y)$ satisfy
\begin{eqnarray} \label{2020.2.11-1}  G^e_1(x,y) =G^e_2(x,y),  \quad (x,y)\in U\times U.\end{eqnarray}
  Pick $y\in U$, and define $V_0 \in C(\partial M_2)$ by
 \begin{eqnarray} \label{2020.2.11-2} V_0( x)= G_1^e (x, y), \,\, x\in \Gamma \cap \mathcal{O}; \;\, V_0(x)=0, \, x\in \partial M_2\setminus \mathcal{O}.\end{eqnarray}
Let $V$ be the solution on $M_2$ to
 \begin{eqnarray*} \left\{\! \begin{array}{ll} {\mathcal{M}}_g V=0 \;\, &\mbox{in}\;\, M_2,\\
 \mbox{div}\; (\sigma E)=0\;\, &\mbox{in}\;\, M_2,\\
  V=V_0 \;\, &\mbox{on}\;\, \partial M_2.\end{array}\right. \end{eqnarray*}
 From the hypothesis $\Lambda_{\Gamma, g_1} =\Lambda_{\Gamma,g_2}$, we have
 \begin{eqnarray*} \nu\times V(x) =\nu\times  G_1^e(x, y) \;\,\mbox{and}\;\,  \nu\times \big(\frac{1}{i\omega \mu} \mbox{curl}\; E(x)\big)= \nu\times \big(\frac{1}{i\omega \mu} \mbox{curl}\; G^e_1 (x,y)\big) \quad \mbox{for}\;\, x\in \Gamma,\end{eqnarray*}
so that, by Holmgren's theorem for Maxwell's equations (see Theorem 6.5 on p.$194$ of \cite{CK2})
and unique continuation of real analytic function (see, for example, p.$\,$65 in \cite{John}), $V$ continues analytically to $\tilde{V}\in C^\infty( {\tilde{M}}_2\setminus \{y\})$, with
$\tilde{V}(x) =G_1^e(x,y)$ for $x\in U\setminus \{y\}$. This satisfies
\begin{eqnarray*}\left\{\!\begin{array}{ll}   {\mathcal{M}}_g\,{\tilde{V}}=  \delta_y\;\,  \mbox{in}\;\,  {\tilde{M}}_2, \\
 \mbox{div}\, \tilde{V}(x)=0 \;\,\mbox{in}\;\, x\in \tilde{M} \setminus \{y\},\\
   {\tilde{V}}\big|_{\partial {\tilde{M}}_2} =0.\end{array}\right.\end{eqnarray*}
    Clearly£¬ $\tilde{V}(x) = G_2^e(x,y)$ for $x\in {\tilde{M}}_2\setminus \{y\}$ since ${\tilde{M}}_2$ has compact closure. Thus, (\ref{2020.2.11-1}) holds.

 \vskip 0.19 true cm
  Secondly, as in \cite{LTU}, we introduce the maps
\begin{eqnarray*} \label{2020.2.2-1}  {\mathcal{G}}_j : {\tilde{M}}_j\to (H^s (U))^3, \quad (\mbox{any} \,\, s<-\frac{1}{2}),\end{eqnarray*}
defined by
\begin{eqnarray*} \label{2020.2.2-2}  {\mathcal{G}}_j (x) (y) =G_j (x,y) , \quad x\in {\tilde{M}}_j , \, y\in U.\end{eqnarray*}
Since $\delta_x\in (H^{s-2} (U))^3$ depends continuously on $x$, we see that ${\mathcal{G}}_j(x) \in (H^s (U))^3$ depends continuously on $x$  for $s<-\frac{1}{2}$, and the maps ${\mathcal{G}}_j, j=1,2$ are $C^1$. Because $G_j(x,y)$ are real-analytic functions of $x$ in ${\tilde{M}}_j \setminus \{y\}$, we conclude that maps ${\mathcal{G}}_j$ are real analytic on ${M}_j$. Noting that the derivative of ${\mathcal{G}}_j$
 \begin{eqnarray*} D{\mathcal{G}}_j(x): T_x {\tilde{M}}_j\to (H^s(U))^3\end{eqnarray*}
 is defined by
 \begin{eqnarray*} D{\mathcal{G}}_j(x) v = v G_j(x,\cdot) =\sum_{k=1}^3 v^k \frac{\partial G_j(x,\cdot)}{\partial x^k} \bigg|_{x},\end{eqnarray*}
 where $v= v^k \frac{\partial }{\partial x_k} \in T_x {\tilde{M}}$, we immediately see that the map $D{\mathcal{G}}_j (x)$ is injective for each $x\in {\tilde{M}}_j$. Furthermore,  it can be  shown that the map ${\mathcal{G}}_j: {\tilde{M}}_j\to (H^s (U))^3$ is an embedding.
  In fact, it remains that to show that $x_1\ne x_2$ in ${\tilde{M}}_j$ implies that ${\mathcal{G}}_j(x_1)\ne {\mathcal{G}}_j(x_2)$. Suppose  this is not the case, then
   \begin{eqnarray} \label{20.4.4-4} G_j^e(x_1,y)=G_j^e(x_2,y)\end{eqnarray}
   for all $y\in U$, therefore, by analyticity, (\ref{20.4.4-4}) holds for all $y\in {\tilde{M}}_j\setminus \{x_1, x_2\}$. However, $G^e(x_1,\cdot)$ is singular only at $y=x_1$ and $G_j(x_2, \cdot)$ is singular only at $y=x_2$, which implies that $x_1=x_2$.

  Finally, it is clear that ${\mathcal{G}}_1$ and ${\mathcal{G}}_2$ coincide in the set $U$ according to (\ref{2020.2.11-1}). Completely similar to proof of Theorem  3.3 of \cite{LTU} we can show that the sets ${\mathcal{G}}_1 ({\tilde{M}}_1)$ and  ${\mathcal{G}}_2 ({\tilde{M}}_2)$ are identical subsets of $(H^s(U))^3$, and the map $J:={\mathcal{G}}_2\circ {\mathcal{G}}_1^{-1} : {\tilde{M}}_1\to  {\tilde{M}}_2$ is an isometry. Combining the arguments in section 3 (see also the proof of Lemma 3.2) and the definition of $J$, we immediately see that $J$ is an identity on $\Gamma$. Therefore, the desired result is immediately obtained. \qed

\vskip 0.33 true cm

\noindent{\bf Remark 4.2.} \  {\it In \cite{Liu}, among other things, the author of this paper also proved that the elastic Dirichlet-to-Neumann map determines the metric $g$ uniquely up to isometry for a strong convex or extendable real-analytic manifold $M$ with boundary. Actually, we can show the isometric result for the elastic problem without topology assumption other than compactness and connectedness by completely similar to the proof of Theorem 1.2.}

\vskip 1.48 true cm

\section{Determining electromagnetic parameters from the electromagnetic Dirichlet-to-Neumann map}

\vskip 0.48 true cm

\noindent  {\it Proof of Theorem 1.4.} \   Let $(x_1,x_2, x_3)$ be the local boundary normal coordinates associated with $(x_1,x_2)$ for $\Gamma \subset \partial \Omega$ as in section 2. According to
(\ref{2019.12.4-1}), the principal symbol $\psi_1^{11}$ of the operator $\Lambda^{11}$ is
\begin{eqnarray*} \label{2019.12.4-3} \psi_1^{11}\!\!&\!\!=\!\!&\! \!\frac{1}{i\omega \mu \sqrt{|g|}} \!\left\{\!\! -\!g_{12}\! \Big( \!\phi_1^{11}\!+\!  \sum\limits_{\alpha} \! g^{\alpha 1} q_{-1} \!\big( \!\sigma i\xi_1\big)\! (i\xi_\alpha)  \!\Big) \!+ \!g_{11}\!\Big(\! \sum\limits_\alpha \! g^{\alpha 1} \!q_{-1} \!\big(\sigma i\xi_2\big)\!(i\xi_\alpha)\!\Big)\!\right\}\\
\!&\!\!=\!\! &\!\!-\,\frac{\xi_1\xi_2}{i\omega \mu \sqrt{ \sum\limits_{\alpha,\beta}\!|g| g^{\alpha\beta}\xi_\alpha\xi_\beta}}. \end{eqnarray*}
Since $g=(g_{jk})$ are the known Riemannian metric on $M$, we see that the parameter $\mu$ is uniquely determined by $\psi_1^{11} (x',\xi')$ on $\Gamma$. Furthermore, all their tangential derivatives along $\Gamma$ are determined by $\psi_1^{11} (x',\xi')$.
By a completely similar discussion of impedance map $\Lambda_{g,\Gamma}^{-1}$, we can show that the highest terms of homogeneity in the expansion of the symbol for $\Lambda_{g,\Gamma}^{-1}$ uniquely determines the $\sigma$ on $\Gamma$.

Next, by (\ref{19.8.28-1}) we see that $\Lambda_{g,\Gamma}$ uniquely determines $L^{jk}$, $\, (1\le j,k\le 2)$. It follows from (\ref{19.9.6-6}) that
\begin{eqnarray}  \label{19.12.10-6} &&l_0^{jj}  = \phi_0^{jj} -\phi_0^{j3} q_{-1} \big(\sigma i\xi_j \big)+\sum\limits_{\alpha} g^{\alpha j}\frac{\partial }{\partial x_\alpha} \Big( q_{-1} (\sigma i\xi_j)\Big)  +\sum\limits_{\alpha} g^{\alpha j} \bigg\{q_{-1} \Big(\sigma \big( \phi_0^{3j}
     \\
     &&\qquad \;\;   +\frac{1}{2} \sum\limits_{\alpha,\beta} g^{\alpha\beta} \frac{\partial g_{\alpha\beta}}{\partial x_j} \big) +\frac{\partial \sigma}{\partial x_j}\Big)+q_{-2} \big(\sigma i\xi_j\big) -i\sum\limits_{m=1}^2 \frac{\partial q_{-1}}{\partial \xi_m} \frac{\partial }{\partial x_m} \big(\sigma i \xi_j\big)\bigg\}
     \, i\xi_\alpha +\sum\limits_{\alpha} g^{\alpha j} \frac{\partial g_{\alpha j}}{\partial x_3} \nonumber.
                    \end{eqnarray}
 Observe that there is not other terms containing $\frac{\partial \mu}{\partial x_3}$ or $\frac{\partial \sigma}{\partial x_3}$ except for $\phi_{0}^{jj}$, $\phi_{0}^{3j}$, $\phi_{0}^{j3}$ and $q_{-2}$ on the right-hand side of the above equality. In view of (\ref{19.8.4-15}) we have
\begin{align*} q_{-2}\! =&\!-\!\frac{1}{\sigma \phi_1^{33}}\Big\{\! q_{-1}\big(\sigma \phi_0^{33} +\frac{\partial \sigma}{\partial x_3}\big)\!\Big\}+S_0^{(1)}(\mu,\sigma),\end{align*}
here and throughout the proof,  each $S_0^{(l)}(\mu,\sigma)$ is an expression involving only the $\mu$, $\sigma$ and their tangential derivatives.
Noting that $\phi_0^{3j} =0 +S_0^{(3)} (\mu,\sigma)$ for $j=1,2$ (see (\ref{19.8.29,3})), we get
\begin{eqnarray} \label{19.12.7-2}\end{eqnarray}
\begin{align*}  l_0^{jj} =& \phi_0^{jj} -\phi_0^{j3} q_{-1} \sigma i\xi_j + \sum\limits_{\alpha} g^{\alpha j} \big(q_{-1} \sigma \phi_0^{3j} +q_{-2} \sigma i \xi_j\big) +S_0^{(2)}(\mu,\sigma)\\
=& \frac{1}{2\sqrt{\sum_{\alpha,\beta} g^{\alpha\beta} \xi_\alpha\xi_\beta}} \bigg(\!\!-\! \frac{1}{\mu} \frac{\partial \mu}{\partial x_3} \sqrt{\sum_{\alpha,\beta} g^{\alpha\beta} \xi_\alpha\xi_\beta}\,\bigg)-\frac{1}{2\sqrt{\sum_{\alpha,\beta} g^{\alpha\beta} \xi_\alpha\xi_\beta}} \bigg( \!-i \!\sum\limits_{\beta} g^{j\beta} \big(\frac{1}{\sigma} \frac{\partial \sigma}{\partial x_3}\big)\xi_\beta\\
& - \!\frac{(-1)^{\!j\!-\!1}i}{\mu|g|} \frac{\partial \mu}{\partial x_3} (g_{3\!-\!j, 2} \xi_1 \!-\!g_{3\!-\!j, 1} \xi_2) \!\bigg) \! \frac{i\xi_j}{\!\sqrt{\sum\limits_{\alpha,\beta} g^{\alpha\beta} \xi_\alpha\xi_\beta\!}\!}+\!\sum\limits_\alpha g^{\alpha j} \!\bigg\{ \! \!\! -\!\frac{1}{\sigma \phi_1^{\!33} }\Big(\!q_{-1} \big(\sigma \phi_0^{33}\! +\!\frac{\partial \sigma}{\partial x_3} \big)\!\Big) \sigma  i\xi_j\!\!\bigg\}i\xi_\alpha \!+\!S_0^{(3)}\!(\mu,\sigma)\nonumber\\
=& -\frac{1}{2\mu} \frac{\partial \mu}{\partial x_3} -\Big(\frac{1}{2\sigma} \frac{\partial \sigma}{\partial x_3} \Big)\frac{\sum_{\beta} g^{j\beta}\xi_\beta\xi_j}{\sum_{\alpha,\beta} g^{\alpha\beta}\xi_\alpha\xi_\beta} -\Big(\frac{1}{ 2\mu |g|} \frac{\partial \mu}{\partial x_3}\Big)  \frac{(-1)^{j-1}}{\sum_{\alpha,\beta} g^{\alpha\beta} \xi_\alpha\xi_\beta} (g_{3-j, 2} \xi_1\xi_j - g_{3-j, 1} \xi_2\xi_j) \nonumber\\
& +\sum\limits_\alpha  g^{\alpha j} \frac{1}{\sigma\sum_{\alpha,\beta} g^{\alpha\beta} \xi_\alpha\xi_\beta} \bigg\{ \sigma\bigg(  \frac{1}{2\sqrt{\sum_{\alpha,\beta} g^{\alpha\beta} \xi_\alpha\xi_\beta} } \Big(\frac{1}{\sigma} \frac{\partial \sigma}{\partial x_3}\Big)\sqrt{\sum_{\alpha,\beta} g^{\alpha\beta} \xi_\alpha\xi_\beta} \bigg) \!+\!\frac{\partial \sigma}{\partial x_3}\bigg\} \xi_j\xi_\alpha +S_0^{\!(3)}\!(\mu,\sigma)\nonumber\\
=& -\frac{1}{2\mu} \frac{\partial \mu}{\partial x_3} -\bigg(\frac{1}{2\sigma} \frac{\partial \sigma}{\partial x_3}\bigg) \frac{\sum_{\beta} g^{j \beta}\xi_\beta\xi_j}{\sum_{\alpha,\beta} g^{\alpha\beta}\xi_\alpha\xi_\beta} -\Big(\frac{1}{ 2\mu |g|} \frac{\partial \mu}{\partial x_3}\bigg)  \frac{(-1)^{j-1}}{\sum_{\alpha,\beta} g^{\alpha\beta} \xi_\alpha\xi_\beta} (g_{3-j, 2} \xi_1\xi_j - g_{3-j, 1} \xi_2\xi_j)\nonumber \\
&+ \sum\limits_\alpha g^{\alpha j} \frac{1}{\sum_{\alpha,\beta} g^{\alpha\beta} \xi_\alpha\xi_\beta} \Big( \frac{3}{2\sigma} \frac{\partial \sigma}{\partial x_3} \Big) \xi_j\xi_\alpha +S_0^{(3)}(\mu,\sigma)\nonumber\\
=& \frac{\partial \mu}{\partial x_3}\! \Big(\!\! \!-\!\frac{1}{2\mu}\! -\! \frac{ (-1)^{j-1}}{2\mu |g|\!\sum_{\alpha,\beta} g^{\alpha\beta} \xi_\alpha\xi_\beta} (g_{3-\!j,2} \xi_1\xi_j\!-\! g_{3-\!j, 1} \xi_2\xi_j)\!\Big)\! +\!\frac{\partial \sigma}{\partial x_3} \Big( \frac{\sum_{\alpha} g^{\alpha j}\xi_\alpha \xi_j}{\sigma\! \sum_{\alpha,\beta} g^{\alpha\beta} \xi_\alpha\xi_\beta}\Big)\!+\! S_0^{(3)}(\mu,\sigma) \,\, \mbox{for} \,\, j\!=\!1,2.  \nonumber\end{align*}
For each $x\in \Gamma$, since the determinant
\begin{align*} \left|\begin{matrix} -\!\frac{1}{2\mu} \!-\! \frac{ g_{22} \xi_1\xi_1- g_{21} \xi_2\xi_1}{2\mu |g|\sum_{\alpha,\beta} g^{\alpha\beta} \xi_\alpha\xi_\beta}  & \frac{\sum_{\alpha} g^{\alpha 1}\xi_\alpha \xi_1}{\sigma \sum_{\alpha,\beta} g^{\alpha\beta} \xi_\alpha\xi_\beta}\\
-\!\frac{1}{2\mu}\! + \!\frac{g_{12} \xi_1\xi_2- g_{11} \xi_2\xi_2}{2\mu |g|\sum_{\alpha,\beta} g^{\alpha\beta} \xi_\alpha\xi_\beta}  & \frac{\sum_{\alpha} g^{\alpha 2}\xi_\alpha \xi_2}{\sigma \sum_{\alpha,\beta} g^{\alpha\beta} \xi_\alpha\xi_\beta}\end{matrix}\right|
= \frac{g^{11}\xi_1^2 -g^{22}\xi_2^2}{2\mu\sigma \sum_{\alpha,\beta} g^{\alpha,\beta} \xi_\alpha\xi_\beta},\end{align*}
we get the above determinant does not vanish when $(\xi_1, \xi_2)\in {\mathbb{R}}^2 \setminus \mathfrak{B}$
because the set $\mathfrak{B}:=\{(\xi_1,\xi_2)\big| g^{11}(x) \xi_1^2 -g^{22}(x) \xi_2^2=0\}$ is at most two one-dimensional curves in ${\mathbb{R}}^2$. For each $x\in \Gamma$, we may choose $(\xi_1,\xi_2)\in {\mathbb{ R}}^2\setminus \mathfrak{B}$. By solving the linear equation system
\begin{align*} \left\{ \begin{array}{ll}\! \frac{\partial \mu}{\partial x_3} \Big(\!\! -\!\frac{1}{2\mu} - \frac{ 1}{2\mu |g|\sum_{\alpha,\beta} g^{\alpha\beta} \xi_\alpha\xi_\beta} (g_{22} \xi_1\xi_j- g_{21} \xi_2\xi_1)\Big) +\frac{\partial \sigma}{\partial x_3} \Big(\! \frac{\sum_{\alpha} g^{\alpha 1}\xi_\alpha \xi_1}{\sigma \sum_{\alpha,\xi_\beta} g^{\alpha\beta} \xi_\alpha\xi_\beta}\Big) =l_0^{11},\\
\! \frac{\partial \mu}{\partial x_3} \Big( \!\!-\!\frac{1}{2\mu} + \frac{ 1}{2\mu |g|\sum_{\alpha,\beta} g^{\alpha\beta} \xi_\alpha\xi_\beta} (g_{12} \xi_1\xi_2- g_{11} \xi_2\xi_2)\Big) +\frac{\partial \sigma}{\partial x_3} \Big( \!\frac{\sum_{\alpha} g^{\alpha 2}\xi_\alpha \xi_2}{\sigma \sum_{\alpha,\xi_\beta} g^{\alpha\beta} \xi_\alpha\xi_\beta}\Big)=l_0^{22},\end{array}\right.\end{align*}
 we see that  $\frac{\partial \mu}{\partial x_3}$ and $\frac{\partial \sigma}{\partial x_3}$ are uniquely determined on $\Gamma$ by $l_1^{11}$ and $l_0^{22}$ (and hence by the operator $\Lambda_{g,\Gamma}$).

 Because there are not $\frac{\partial^2 \mu}{\partial x_3^2}$ and $\frac{\partial^2 \sigma}{\partial x_3^2}$ in the expression of $l_0^{jj}$  except for $\phi^{jj}_{-1}$, $\phi_{-1}^{j3}$, $\phi^{3j}_{-1}$ and $q_{-3}$,
it follows from (\ref{2019.12.6-1}) that \begin{align*} j_{-1}^{jj}= \phi_{-1}^{jj} -\phi_{-1}^{j3}q_{-1}(\sigma \,i\xi_j) +\sum\limits_{\alpha} g^{\alpha j} \Big( q_{-1} \sigma \phi_{-1}^{3j} +q_{-3} \sigma i \xi_j\Big) i\xi_\alpha+S_{-1}^{(1)}(\mu,\sigma), \quad \, j=1,2,\end{align*}
where each $S_{-1}^{(m)}(\mu,\sigma)$ denotes an expression involving only $\mu$, $\sigma$, $\frac{\partial \mu}{\partial x_3}$ and $\frac{\partial \sigma}{\partial x_3}$.
According to (\ref{19.8.4-16,}), we see that
\begin{align*} q_{-3}& =-\frac{1}{\sigma \phi_{1}^{33}} \Big(q_{-1} \sigma \phi_{-1}^{33}\Big)+ S_{-1}^{(2)} (\mu,\sigma)\\
=& -\frac{1}{\sigma \sum_{\alpha,\beta}g^{\alpha\beta}\xi_\alpha\xi_\beta} \phi_{-1}^{33}+S_{-1}^{(2)} (\mu,\sigma)\\
=&  - \frac{1}{\sigma \sum_{\alpha,\beta}g^{\alpha\beta}\xi_\alpha\xi_\beta}\bigg\{ \frac{1}{2\sqrt{\sum_{\alpha,\beta} g^{\alpha\beta} \xi_\alpha \xi_\beta}}\Big(\frac{\partial \phi_0^{33}}{\partial x_3} -\frac{\partial }{\partial x_3} \big(\frac{1}{\sigma} \frac{\partial \sigma}{\partial x_3}\big)\Big) \bigg\}+ S_{-1}^{(3)} (\mu,\sigma)\\
=& -\frac{1}{2 \sigma \big( \sum_{\alpha,\beta} g^{\alpha\beta} \xi_\alpha \xi_\beta\big)^{3/2}}\bigg\{ \frac{\partial}{\partial x_3} \Big(\frac{1}{2\sqrt{\sum_{\alpha,\beta} g^{\alpha\beta} \xi_\alpha\xi_\beta}} \big( \frac{1}{\sigma} \frac{\partial \sigma}{\partial x_3}
\sqrt{\sum_{\alpha,\beta} g^{\alpha\beta} \xi_\alpha\xi_\beta}\big) \Big)
-\frac{\partial }{\partial x_3} \big(\frac{1}{\sigma} \frac{\partial \sigma}{\partial x_3}\big)\bigg\} + S_{-1}^{(4)} (\mu,\sigma)\\
=&   \frac{1}{4\sigma \big( \sum\limits_{\alpha,\beta} g^{\alpha\beta} \xi_\alpha\xi_\beta\big)^{\frac{3}{2}}}\, \frac{1}{\sigma} \frac{\partial^2 \sigma}{\partial x_3^2} +S_{-1}^{(5)}(\mu,\sigma),\end{align*}
so that for $j=1,2$,
\begin{eqnarray} \label{20.3.30-7} \end{eqnarray}
\begin{align*} l_{-1}^{jj} \!= & \phi_{-1}^{jj}\! -\!\phi_{-1}^{j3} \frac{ i\xi_j}{ \sqrt{\sum\limits_{\alpha,\beta} g^{\alpha\beta}\xi_\alpha\xi_\beta}}+\! \sum\limits_{\alpha} g^{\alpha j}\Bigg\{\! \frac{1}{\sqrt{\sum\limits_{\alpha,\beta} g^{\alpha\beta}\xi_\alpha\xi_\beta}} \phi_{-1}^{3j} +\Big(\!\frac{1}{4\sigma \big(\sum\limits_{\alpha,\beta} g^{\alpha\beta} \xi_\alpha\xi_\beta\big)^{\!\frac{3}{2}}}\frac{1}{\sigma} \frac{\partial^2\sigma}{\partial x_3^2} \Big) \sigma i\xi_j\! \Bigg\} i\xi_\alpha \!+\! S_{-1}^{(6)}(\mu,\sigma) \\
\!=    & \phi_{-1}^{jj} \!- \!\phi_{-1}^{j3} \frac{i\xi_j}{\sqrt{\sum_{\alpha,\beta} g^{\alpha\beta}\xi_\alpha\xi_\beta }} \!+\!\sum\limits_{\alpha}\!  \!\frac{g^{\alpha j}}{\sqrt{\sum_{\alpha,\beta} g^{\alpha\beta}\xi_\alpha\xi_\beta}} \phi_{-1}^{3j} i\xi_\alpha \!-\!
\sum\limits_{\alpha} \!\bigg( \!\frac{g^{\alpha j}}{ 4\sigma \big(\sum_{\alpha,\beta}\! g^{\alpha\beta} \xi_\alpha\xi_\beta\big)^{\!\frac{3}{2}} }\, \frac{\partial^2 \sigma}{\partial x_3^2} \xi_j\xi_\alpha\!\bigg)\!+\!S_{-1}^{(7)} (\mu,\sigma)\\
\!=      & \frac{1}{2\sqrt{\sum\limits_{\alpha,\beta} g^{\alpha\beta}\xi_\alpha\xi_\beta}}  \frac{\partial \phi_0^{jj}}{\partial x_3} -  \frac{1}{2\sqrt{\sum\limits_{\alpha,\beta} g^{\alpha\beta}\xi_\alpha\xi_\beta}} \Big(\frac{\partial \phi_{0}^{j3}}{\partial x_3}\Big)
\frac{i\xi_j}{\sqrt{\sum_{\alpha,\beta} g^{\alpha\beta} \xi_\alpha\xi_\beta}}  + \sum\limits_{\alpha} g^{\alpha j}
 \frac{1}{2 \sum_{\alpha,\beta} g^{\alpha\beta}
\xi_\alpha\xi_\beta} \Big(\frac{\partial \phi_0^{3j}}{\partial x_3} \Big) i\xi_\alpha \\
\!&        - \sum\limits_{\alpha}
g^{\alpha  j}  \Big( \frac{1}{4\sigma \big(\sum_{\alpha, \beta} g^{\alpha\beta}\xi_\alpha\xi_\beta\big)^{\frac{3}{2}}}\Big) \frac{\partial^2 \sigma}{\partial x_3^2} \xi_j\xi_\alpha\!+\! S_{-1}^{(8)} (\mu,\sigma) \\
\!=&  \frac{1}{2\sqrt{\sum_{\alpha,\beta} g^{\alpha\beta} \xi_\alpha\xi_\beta}} \,\frac{\partial }{\partial x_3} \Big\{
\frac{1}{2\sqrt{\sum_{\alpha,\beta} g^{\alpha\beta} \xi_\alpha\xi_\beta}}\Big(\! -\frac{1}{\mu} \frac{\partial \mu}{\partial x_3} \sqrt{\sum_{\alpha,\beta} g^{\alpha\beta} \xi_\alpha\xi_\beta} \Big)\Big\} \\
\!&  +
\frac{i\xi_j}{4(\sum_{\alpha,\beta} g^{\alpha\beta} \xi_\alpha\xi_\beta)^{\frac{3}{2}}} \, \frac{\partial }{\partial x_3} \Big(\! -i \sum\limits_{\beta}
g^{j\beta} \big(\frac{1}{\sigma} \frac{\partial \sigma}{\partial x_3} \big) \xi_\beta -\frac{(-1)^{j-1} i}{\mu|g|} \frac{\partial \mu}{\partial x_3}(g_{3-j,2}\xi_1- g_{3-j, 1} \xi_2) \Big)+0\\
\!&  - \sum\limits_{\alpha} g^{\alpha j} \Big(\frac{1}{4\sigma\big( \sum_{\alpha,\beta} g^{\alpha\beta} \xi_\alpha\xi_\beta\big)^{\!\frac{3}{2}}}\Big)
\frac{\partial^2 \sigma}{\partial x_3^2} \xi_j\xi_\alpha  +\!S_{-1}^{(9)} (\mu,\sigma)\\
\!=& \frac{\partial^2 \mu}{\partial x_3^2} \Big(\!- \frac{1}{4\mu\sqrt{\sum_{\alpha,\beta}g^{\alpha\beta}\xi_\alpha\xi_\beta}}  +
\frac{(-1)^{j-1}(g_{3-j,2}\xi_1-g_{3-j,1}\xi_2) \xi_j}{4\mu|g|\big(\sum_{\alpha,\beta}g^{\alpha\beta}\xi_\alpha\xi_\beta\big)^{3/2}}\Big)+ \frac{\sum_{\alpha}g^{j\beta} \xi_j\xi_\beta}{4\sigma \big(\sum_{\alpha} g^{\alpha\beta} \xi_\alpha\xi_\beta\big)^{\!\frac{3}{2}} } \, \frac{\partial^2\sigma}{\partial x_3^2}+ \!S_{-1}^{(10)} (\mu,\sigma).\end{align*}
By taking $j=1,2$ and solving the above linear equation system, we immediately get that the $l^{11}_{-1}$ and $l_{-1}^{22}$ uniquely determines  $\frac{\partial^2 \mu}{\partial x_3^2}$ and $\frac{\partial^2 \sigma}{\partial x_3^2}$ on $\Gamma$.

Finally, we will consider the symbol $l^{11}_{-m-1}+l^{22}_{-m-1}$ with homogeneous of degree $-m-1$, $\,(m\ge 1)$, for $L^{11}+L^{22}$.
Obviously, there is not other terms containing $\frac{\partial^{m+2}\mu}{\partial x_3^{m+2} }$ and $\frac{\partial^{m+2} \sigma}{\partial x_3^{m+2}}$ except for  $\phi_{-1}^{jj}$, $\phi_{-1}^{3j}$, $\phi_{-1}^{j3}$ and $q_{-3}$ on the right-hand side of (\ref{19.7.21-7}), so that
        \begin{eqnarray*} l^{jj}_{-m-1} = \phi_{-m-1}^{jj} - \phi_{-m-1}^{j3} q_{-1}\sigma \,i\xi_j  + \sum\limits_{\alpha} g^{\alpha j} \bigg(q_{-1}\sigma \phi^{3j}_{-m-1}   +  q_{-m-3}\sigma \, i\xi_j \bigg) i\xi_\alpha+S_{-m-1}^{(1)}(\mu,\sigma), \quad \mbox{for}\,\, j=1,2,\end{eqnarray*}
 where each $S_{-m-1}^{(s)}$ is an expression involving only the boundary values of $\mu$, $\sigma$, and their normal derivatives of order at most $m+1$.
From (\ref{19.8.12-5}) we have
\begin{eqnarray*} q_{-m-3} = -\frac{ 1}{\sigma\phi_1^{33}} \,q_{-1} \,\sigma \phi^{33}_{-m-1} + S^{(2)}_{-m-1}(\mu,\sigma).\end{eqnarray*}
  It follows that \begin{eqnarray} \label{19.12.7-4} && l_{-m-1}^{jj} = \phi_{-m-1}^{jj} -\frac{i\xi_j}{ \sqrt{ \sum_{\alpha,\beta} g^{\alpha\beta}\xi_\alpha\xi_\beta}}\phi_{-m-1}^{j3}  + \sum\limits_{\alpha} g^{\alpha j} \bigg( \frac{1}{ \sqrt{ \sum_{\alpha, \beta} g^{\alpha\beta} \xi_\alpha\xi_\beta}}  \phi_{-m-1}^{3j} \\
  &&\quad  -\frac{1}{ \phi_1^{33}} \frac{1}{\sqrt{\sum_{\alpha,\beta} g^{\alpha\beta} \xi_\alpha\xi_\beta}} \phi_{-m-1}^{33}  i\xi_j\bigg) i\xi_\alpha
   +S_{-m-1}^{(3)} (\mu,\sigma)\nonumber\\
   && =\phi_{-m-1}^{jj} -\frac{i\xi_j}{\sqrt{\sum_{\alpha,\beta} g^{\alpha\beta} \xi_\alpha\xi_\beta}}\phi_{-m-1}^{j3}
      +\frac{1}{\sqrt{\sum_{\alpha,\beta}g^{\alpha\beta} \xi_\alpha\xi_\beta }} \sum\limits_{\alpha} g^{\alpha j} i\xi_\alpha \phi_{-m-1}^{3j}\nonumber\\
      && \quad  +
      \frac{1}{\sum_{\alpha,\beta}g^{\alpha\beta}\xi_\alpha\xi_\beta } \big(\sum\limits_{\alpha} g^{\alpha j} \xi_j\xi_\alpha\big)\phi_{-m-1}^{33}+ S_{-m-1}^{(3)}(\mu,\sigma).\nonumber \end{eqnarray}
 From (\ref{19.8.3-1}) we have
 \begin{eqnarray} \label{19.9.4-1:} \phi_{-m-1}^{jk} = \frac{1}{2\sqrt{\sum_{\alpha,\beta} g^{\alpha\beta} \xi_\alpha\xi_\beta}} \frac{\partial \phi_{-m}^{jk}}{\partial x_3} +S_{-m-1}^{(4)}(\mu,\sigma) \quad \,\mbox{for}\;\, j,k=1,2,3.\end{eqnarray}
   We will end our proof by induction. Suppose we have shown that for $1\le r\le m$,
\begin{eqnarray} \label{19.19.4-6}\;\; l_{-r}^{jj}=&&\!\!\!\!\!\!\!\!\!\! \bigg(2\sqrt{\sum\limits_{\alpha,\beta}\! g^{\alpha\beta} \xi_\alpha\xi_\beta}\,\bigg)^{-r} \bigg\{
\frac{\partial^{r+1} \mu}{\partial x_3^{r+1}} \Big(\!- \frac{1}{2\mu}  +
\frac{(-1)^{j-1}(g_{3-j,2}\xi_1-g_{3-j,1}\xi_2) \xi_j}{2\mu|g|\big(\sum_{\alpha,\beta}g^{\alpha\beta}\xi_\alpha\xi_\beta\big)}\Big)\\
&& \!\!\!\!\!\!\!\!\!+ \frac{\sum_{\alpha}g^{j\beta} \xi_j\xi_\beta}{2\sigma \sum_{\alpha} g^{\alpha\beta} \xi_\alpha\xi_\beta}\,  \frac{\partial^{r+1}\sigma}{\partial x_3^{r+1}}\!\bigg\} +  S_{-r}^{(5)}(\mu,\sigma).\nonumber\end{eqnarray}
Clearly, this estimates actually holds when $r=1$ by (\ref{20.3.30-7}). It follows from (\ref{19.12.7-4}) and (\ref{19.8.3-1}) that
\begin{eqnarray} \label{19.12.7-5}  l_{-m-1}^{jj}&&\!\!\!\!\!\!\!\!\! = \frac{1}{2\sqrt{\sum_{\alpha,\beta} g^{\alpha\beta} \xi_\alpha\xi_\beta}}
    \frac{\partial }{\partial x_3}\left\{
      \phi_{-m}^{jj} -\frac{i\xi_j}{\sqrt{\sum_{\alpha,\beta} g^{\alpha\beta} \xi_\alpha\xi_\beta}}\phi_{-m}^{j3}
      +\frac{1}{\sqrt{\sum_{\alpha,\beta}g^{\alpha\beta} \xi_\alpha\xi_\beta }} \sum\limits_{\alpha} g^{\alpha j} i\xi_\alpha \phi_{-m}^{3j}\right.\nonumber\\
      &&\!\!\left.   +
      \frac{1}{\sum_{\alpha,\beta}g^{\alpha\beta}\xi_\alpha\xi_\beta } \big(\sum\limits_{\alpha} g^{\alpha j} \xi_j\xi_\alpha\big)\phi_{-m}^{33}+ S_{-m}^{(6)}(\mu,\sigma)\right\} +S_{-m-1}^{(7)}(\mu,\sigma)\nonumber\\
        &&\!\!\!\!\!\!\!\!\!=\frac{1}{2\sqrt{\sum_{\alpha,\beta} g^{\alpha\beta} \xi_\alpha\xi_\beta}} \frac{\partial l_{-m}^{jj}}{\partial x_3}+S_{-m-1}^{(8)}(\mu,\sigma).\nonumber
      \end{eqnarray} Combing this and (\ref{19.19.4-6}) we have
      \begin{eqnarray} \label{19.12.7-6} &&\; \;\; l_{-m-1}^{jj} =\frac{1}{2\sqrt{\sum_{\alpha,\beta} g^{\alpha\beta} \xi_\alpha\xi_\beta}} \frac{\partial}{\partial x_3}\bigg\{\!
       \bigg(2\sqrt{\sum\limits_{\alpha,\beta}\! g^{\alpha\beta} \xi_\alpha\xi_\beta}\,\bigg)^{\!-m} \bigg[
\frac{\partial^{m+1} \mu}{\partial x_3^{m+1}} \bigg(\!\!- \!\frac{1}{2\mu} \\
&&\quad \;\; \qquad +
\frac{(-1)^{j-1}(g_{3-j,2}\xi_1-g_{3-j,1}\xi_2) \xi_j}{2\mu|g|\big(\sum_{\alpha,\beta}g^{\alpha\beta}\xi_\alpha\xi_\beta\big)}\!\bigg)\!\!+\! \frac{\sum_{\alpha}g^{j\beta} \xi_j\xi_\beta}{2\sigma \sum_{\alpha} g^{\alpha\beta} \xi_\alpha\xi_\beta }\frac{\partial^{m+1}\sigma}{\partial x_3^{r+1}}\!\bigg] \!\!+\!  S_{\!-m}^{(9)}(\mu,\sigma)\!
      \bigg\}\!\!+\!  S_{\!-m-1}^{(10)}(\mu\!,\sigma)\nonumber\\
   && \quad \qquad =
       \bigg(2\sqrt{\sum\limits_{\alpha,\beta}\! g^{\alpha\beta} \xi_\alpha\xi_\beta}\bigg)^{\!-m-1} \bigg\{\!
\frac{\partial^{m+2} \mu}{\partial x_3^{m+2}} \Big(\!\!- \!\frac{1}{2\mu}  +
\frac{(-1)^{j-1}(g_{3-j,2}\xi_1-g_{3-j,1}\xi_2) \xi_j}{2\mu|g|\big(\sum_{\alpha,\beta}g^{\alpha\beta}\xi_\alpha\xi_\beta\big)}\Big)\nonumber\\
&&\quad \;\,\; \qquad + \frac{\sum_{\alpha}g^{j\beta} \xi_j\xi_\beta}{2\sigma \sum_{\alpha} g^{\alpha\beta} \xi_\alpha\xi_\beta}  \; \frac{\partial^{m+2}\sigma}{\partial x_3^{m+2}}\!\bigg\} +  S_{-m-1}^{(11)}(\mu,\sigma) \quad \; \mbox{for} \;\; j=1,2.\nonumber   \end{eqnarray}
Solving  the above linear equation system, we get that $\frac{\partial^{m+2} \mu}{\partial x_3^{m+2}}$ and $\frac{\partial^{m+2} \sigma}{\partial x_3^{m+2}}$  can uniquely be determined on $\Gamma$ by $l_{-m-1}^{11}$ and  $l_{-m-1}^{11}$.
More precisely, we can solve for the first $(m+2)$-order normal derivatives of $\mu$ and $\sigma$ in terms of $l_{-k}^{jj}$ for $k=1,2,\cdots, m+1$ and $j=1,2$. \qed

\vskip 0.58 true cm

  \noindent  {\it Proof of Theorem 1.5.} \  Let $(x_1, x_{2})$ be any local coordinates for an open set $W\subset \Gamma$,  and let $\{\psi_j\}_{j\le 1} $ denote the full symbol of $\Lambda_{g,\Gamma}$ in these coordinates. Then for any $x_0\in W$, $\frac{\partial^{|K|} \mu}{\partial x^K}$ and $\frac{\partial^{|K|} \sigma}{\partial x^K}$  for all multi-indices $K=(k_1,k_2,k_3)$ with $|K|\ge 0$ at $x_0$ in boundary normal coordinates is given by explicit formula in terms of the matrix-valued functions $\{ \psi_j\}_{j\le 1}$
 and their tangential derivatives at $x_0$. This implies that we can determines the functions $\mu$ and $\sigma$ at a small neighborhood of $x_0$ by the real-analyticity of $\mu$ and $\sigma$ on $M\cup \Gamma$. Hence, by unique continuation of real analytic function (see, for example, p.$\,$65 in \cite{John}), we can uniquely determine $\mu$ and $\sigma$ on real-analytic manifold $(M,g)$. \qed

\vskip 0.33 true cm
\noindent{\bf Remark 5.1.} \  {\it It is easy to verify that Theorem 1.5 still holds for piece-wise real-analytic manifold $(M,g)$. }

\vskip 1.68 true cm

\section*{Acknowledgments}
\renewcommand{\thesection}{\arabic{section}}
\renewcommand{\theequation}{\thesection.\arabic{equation}}
\setcounter{equation}{0} \setcounter{maintheorem}{0}

\vskip 0.59 true cm
  I would like to thank Professor Gunther Uhlmann for his very useful comments and suggestions, which make me to improve an earlier result of Theorem 1.2. This research was supported by NNSF of China (11671033/A010802).

  \vskip 1.68 true cm


\begin{thebibliography}{9999}



\bibitem{ADN}  S. Agmon, A. Douglis and L. Nirenberg, \textsl{Estimates near the
boundary for solutions of elliptic partial differential equations
satisfying general boundary conditions. I},
 Comm. Pure Appl. Math.  12(1959), 623-727.




\bibitem {Ale} G. Alessandrini,  \textsl{Stable determination of conductivity by boundary measurements}.
Appl. Anal., (1-3), 27(1988), 153-172.


\bibitem {AP}  K. Astala  and L. P\"{a}iv\"{a}rinta,  \textsl{Calder\'{o}n's
  inverse conductivity problem in the plane}, Ann. Math., 163(2006), 265-299.


\bibitem{Bro} R. M. Brown,  \textsl{Global uniqueness in the impedance-imaging problem for less regular
conductivities}, SIAM J. Math. Anal.  no.4, 27(1996), 1049-1056.

\bibitem{BU} B. M. Brown, G. A. Uhlmann, \textsl{Uniqueness in the inverse conductivity problem with less regular conductivities
in two dimensions}, Comm. Partial Differential Equations, no. 5-6,
 22(1997), 1009-1027.


\bibitem{Br} R. M. Brown, . \textsl{Global uniqueness in the impedance-imaging problem for less regular
conductivities}, SIAM J. Math. Anal., no.4, 27(1996), 1049-1056.

\bibitem{BT}  R. M. Brown and R. H. Torres, \textsl{Unqueness in the inverse conductivity problem for conductivities with
$3/2$ derivatives in $L^p, p>2n$}, J. Fourier Anal. Appl., no.6, 9(3003), 563-574.

\bibitem{ABU} A. L. Bukhgeim and G. Uhlmann, \textsl{Recovering a potential from partial Cauchy data}, Comm.
Partial Differential Equations, No.3-4, 27(2002), 653-668.


\bibitem{Cald} A.P. Calder\'{o}n, \textsl{On an inverse boundary value problem}, Seminar on Numerical Analysis and
its Applications to Continuum Physics (Rio de Janeiro, 1980), pp. 65-73, Soc. Brasil. Mat.,
R\'{i}o de Janeiro, 1980.



\bibitem{Cha} S. Chanillo, \textsl{A problem in electrical prospection and ann-dimensional Borg-Levinson
theorem},  Proc. Amer. Math. Soc., no. 3, 108(1990), 761-767.

\bibitem{CR} P. Caro and K.M Rogers, \textsl{Glabal uniqueness for the Calder\'{o}n problem with Lipschitz conductivities}, Forum Math. Pi, e2, 4(2016).


\bibitem{COR}  P. Caro, P. Ola, and M. Salo, \textsl{Inverse boundary value problem for Maxwell equations with
   local data}, Comm. Partial Differential Equations, 34(2009), 1425-1464.


\bibitem{CZ} P. Caro and Ting Zhou,   \textsl{Global uniqueness for an IBVP for the time-harmonic Maxwell equations},
Anal. PDE, No.2,  7(2014), 375-405.

\bibitem{Chew} W. C. Chow,  \textsl{Waves and Fields in Inhomogeneous Media}, Institute of Electrical and Electronics Engineers, Inc., New York, 1995.


\bibitem{CK2} D. Colton and R. Kress, \textsl{Inverse Acoustic and Electromagnetic Scattering Theory}, Springer-Verlag, New York, 1992.

\bibitem{Coo}  D. M. Cook, \textsl{The Theory of the Electromagnetic Field}, Mineola NY: Courier Dover Publications, 2002.


\bibitem{COST} F. J. Chung, P. Ola, M. Salo, and L. Tzou, \textsl{Partial Data Inverse Problems for Maxwell
Equations via Carleman Estimates}, ArXiv e-prints, February 2, 2015.


\bibitem{Esk} G. I. Eskin, \textsl{Boundary Value Problemsfor Elliptic Pseudodigerential Equations} (translated from
Russian by S. Smith), Ann. Math. Society Translation of Mathematical Monographs, Vol. 52,
Providence, R.I., 1981.

\bibitem{FKSU} D. S. Ferreira, C. E. Kenig, M. Salo and G. Uhlmann, \textsl{Limiting Carleman weights and anisotropic inverse
problems}, Invent. Math., 178(2009), 119-171.


\bibitem{Frie} A. Friedman, \textsl{Partial Diferential Equations of Parabolic Type}, Prentice Hall, Englewood Cliffs, NJ, 1964.

\bibitem{GT}  D. Gilbarg  and N. Trudinger, \textsl{Elliptic Partial
Differential Equations of Second Order}, Reprint of the 1998
edition, Classics Math. Springer-Verlag, Berlin, 2001.


\bibitem{Gil2} P. Gilkey, \textsl{Invariance Theory, the Heat Equation and the Atiyah-Singer
Index Theorem}, CRC Press, Boca Raton, 1995.


\bibitem{GKLU}  A. Greenleaf,  Y. Kurylev, M. Lassa and G. Uhlmann, \textsl{Invisibility and inverse problems},
 Bulletin of American Mathematical Society,
   46(2009), 55-97.


\bibitem{GLU} A. Greenleaf, M. Lassas, G. Uhlmann,  \textsl{The Calder¨®n problem for conormal
potentials,  I. Global uniqueness and reconstruction}, Comm. Pure Appl. Math., no. 3,  56(2003),
328-352.


\bibitem{Gr}  G. Grubb, \textsl{Functional Calculus of Pseudo-differential Boundary Problems}, Birkh\"{a}user, Boston, 1986.



\bibitem{Ho3} L. H\"{o}rmander, \textsl{The Analysis of Partial Differential Operators} III,
  Springer-Verlag, Berlin Heidelberg New York, 1985.


\bibitem{Isak} V. Isakov, \textsl{Inverse Problems for Partial Differential Equations}, 2nd ed., Springer, New York, 2006.


\bibitem{Isako} V. Isakov, \textsl{On uniqueness in the inverse conductivity problem with local data}, Inverse Probl. Imaging,
No.1, 1(2007), 95-105.


 \bibitem{John} F. John, \textsl{Partial Differential Equations}, fourth ed., Springer-Verlag, New York Inc., 1982.


\bibitem{JoMcD} M. S. Joshi and S. R. McDowall, \textsl{Total determination of material parameters from electromagnetic boundary information}, Pacific Journal of Mathematics, 193(2000), 107-129.

\bibitem{KeSU} C. E. Kenig, M. Salo and G. Uhlmann, \textsl{Inverse problems for the anisotropic Maxwell equations}, Duke Math. J., 157(2011), 369-419.


\bibitem{KeSjU} C. E. Kenig, J. Sj\"{o}strand and G. Uhlmann, \textsl{The Calder¡äon problem with partial data},
Annals of Mathematics, No.2, 165(2007), 567-591.



\bibitem{Kir} A. Kirsch, \textsl{An introduction to mathematical theory of inverse problems}, Second Edition, Springer Science+Business Media, LLC, 2011.

\bibitem{KN} J.  Kohn  and L. Nirenberg, \textsl{An algebra of pseudo-differential operators}. Comm. Pure Appl. Math., 18(1965), 269-305.



\bibitem{KT} K. Knudsen, A. Tamasan, \textsl{Reconstruction of less regular conductivities in the plane},
Comm. Partial Differential Equations, no. 3-4, 29(2004), 361-381.

\bibitem{LTU} M. Lassas, M. Taylor  and G.  Uhlmann,  \textsl{The dirichlet-to-neumann map for complete Riemannian manifolds with boundary},  Comm. Geom. Anal.,  11(2003), 207-222.

\bibitem{LU}  J.  Lee and G. Uhlmann, \textsl{Determing anisotropic real-analytic conductivities by boundary measurements},
Comm. Pure Appl. Math., 42(1989), 1097-1112.


\bibitem{Liu} G. Q. Liu, \textsl{Determination of isometric real-analytic metric and spectral invariants for elastic Dirichlet-to-Neumann map on Riemannian manifolds}, arXiv: 1908.05096 [math.AP].


\bibitem{McD} S. R. McDowall, \textsl{Boundary determination of material parameters from electromagnetic boundary information}, Inverse Problems 13(1997), 153-163.


\bibitem{Myers} S. D. Myers,  \textsl{Riemannian manifolh in the large}, Duke Math. J.  1(1935), 39-49.

\bibitem{Nac1} A. R. Nachman,  \textsl{Global uniqueness for a two-dimensional inverse boundary value
problem}, Ann. of Math., no. 1, 143(1996), 71-96.


\bibitem{NSU}  A. I. Nachman, J. Sylvester and G. A.  Uhlmann, \textsl{An n-dimensional Borg-Levinson
theorem},  Comm. Math. Phys., no. 4, 115(1988), 595-605.



\bibitem{Nov} R. G. Novokov, \textsl{A multidimensional inverse spectral problem for the equation $-\Delta \phi +(v(x) -Eu(x)) \psi=0$}, Funksional  Annal. i Prilozhen, no.4, 22(1988), 11-22, Translation in Func. Anal. Appl., 22(1988), 263-272.

\bibitem{Pich} M. Pichler, \textsl{An inverse problem for Maxwell's equations with Lipschitz parameters},
  Inverse Problems, 34(2018), 1-21.


\bibitem{OPS} P. Ola, L. P\"{a}iv\"{a}rinta and E. Somersalo, \textsl{An inverse boundary problem in electrodynamics}, Duke Math. J.
70(1993), 617-653.

\bibitem{OPS2} P. Ola, L. P\"{a}iv\"{a}rinta and E. Somersalo, \textsl{Inverse problems for time harmonic electrodynamics},
In {\it Inside out: inverse problems and applications}, volume 47 of Math. Sci. Res. Inst. Publ., pages 169-191.
Cambridge Univ. Press, Cambridge, 2003.


\bibitem{OS} P. Ola,  E. Somersalo, \textsl{Electromagnetic inverse problems and generalized Sommerfeld potentials}.
SIAM J. Appl. Math., No.4, 56(1996), 1129-1145.

\bibitem{PPU} L. P\"{a}iv\"{a}rinta, A. Panchenko and G. A. Uhlmann, \textsl{Complex geometrical optics solutions
for Lipschitz conductivities}, Rev. Mat. Iberoamericana  no.1, 19(2003), 56-72.



\bibitem{SIC}  E. Somersalo,  D. Isaacson and M. Cheney,  \textsl{A linearized inverse boundary value problem for Maxwell¡¯s
equations},  J. Comput. Appl. Math. 42(1992),  123-36.



\bibitem{Som} E. Somersalo,  \textsl{Layer stripping for time harmonic Maxwell¡¯s equations with high frequency}, Inverse
Problems 10(1994), 449-66.


\bibitem{SunU} Z. Sun and G. Uhlmann,  \textsl{An inverse boundary value problem for Maxwell's equations}, Arch.Rat. Mech. Anal. 119(1992), 71-93.


\bibitem{SU} Z. Sun and G. Uhlmann,  \textsl{Anisotropic inverse problems in two dimensions},
 Inverse Problems, 19(2003),1001-1010.


 \bibitem{SU1} J. Sylvcster and G. Uhlmaun,  \textsl{A global uniqueness theorem for an inverse boundary
value problem}, Ann. Math., 125 (1987), 153-169.



\bibitem{SU2} J. Sylvester and G. Uhlmann, \textsl{The Dirichlet to
 Neumann map and applications},
  in: Inverse problems in partial differential equations, Edited by
  David Colton, the Society for Industrial and Applications, 1990.


\bibitem{Ta2} M. E. Taylor, \textsl{Partial Differential
Equations II}, Appl. Math. Sci., vol. 116, Springer-Verlag, New York, 1996.

\bibitem{Ta3} M. E. Taylor, \textsl{Partial Differential
Equations III}, Appl. Math. Sci., vol. 117, Springer-Verlag, New York, 1996.

\bibitem{Tre} F. Treves, \textsl{Introduction to pseudodifferential and Fourier integral operator}, Plenum Press, New York, 1980.


\bibitem{Whi} G. Whitehead,  \textsl{Elements of Homotopy theory}, Springer-Verlag, New York, 1978.

\bibitem{Youn} E. C. Young, \textsl{Vector and tensor analysis},  Second Edition, Marcel Dekker, Inc., New York, Basel, Hong Kong, 1993.




\end{thebibliography}
\end{document}